\documentclass[10pt]{article}

\usepackage[left=1in,top=1in,right=1in,bottom=1in,letterpaper]{geometry}

\usepackage[utf8]{inputenc} 
\usepackage{url}    

\usepackage{booktabs}      
\usepackage{amsfonts}       
\usepackage{microtype}      
\usepackage[table]{xcolor}         

\usepackage{amsmath}
\usepackage{amsthm}
\usepackage{amssymb}

\usepackage{multirow}
\usepackage{enumitem}
\usepackage{dsfont}

\usepackage{graphicx}

\usepackage{algorithm}
\usepackage{algorithmicx}
\usepackage{algpseudocode}

\usepackage{threeparttable}
\usepackage{mathtools}
\usepackage{float}
\usepackage{placeins}
\usepackage{chngcntr}
\usepackage[colorlinks,linkcolor=blue,filecolor=blue,citecolor=blue,urlcolor=blue,pdfusetitle]{hyperref}

\restylefloat{algorithm} 
\newtheorem{assumption}{Assumption}
\newtheorem{theorem}{Theorem}

\newtheorem{lemma}{Lemma}
\newtheorem{example}{Example}

\newlength{\wrapBottomVspace}
\numberwithin{equation}{section}
\counterwithin{figure}{section}
\counterwithin{table}{section}
\counterwithin{algorithm}{section}

\newcommand{\vg}{\mathbf{g}}

\newcommand{\vF}{{\mathbf{F}}}

\newcommand{\vW}{{\mathbf{W}}}

\newcommand{\bb}{{{b}}}

\newcommand{\bd}{{{d}}}

\newcommand{\bg}{{{g}}}

\newcommand{\bu}{{{u}}}
\newcommand{\bv}{{{v}}}

\newcommand{\bx}{{{x}}}

\newcommand{\bz}{{{z}}}

\newcommand{\bB}{{{B}}}

\newcommand{\bF}{{{F}}}
\newcommand{\bG}{{{G}}}

\newcommand{\bM}{{{M}}}

\newcommand{\bP}{{{P}}}

\newcommand{\bR}{{{R}}}

\newcommand{\bU}{{{U}}}
\newcommand{\bV}{{{V}}}
\newcommand{\bW}{{{W}}}

\newcommand{\bY}{{{Y}}}

\newcommand{\vvvert}{{\vert\kern-0.25ex\vert\kern-0.25ex\vert}}

\DeclareMathOperator*{\argmin}{argmin}

\makeatletter
\let\c@lofdepth\relax
\let\c@lotdepth\relax

\makeatother
\usepackage{tocloft}
\usepackage{setspace}

\allowdisplaybreaks

\title{A Primer on Matrix Subspace Optimization}

\author{%
  Yihan Jin\textsuperscript{1}\thanks{Equal contribution. Emails:
    \href{mailto:jyh2714225423@stu.pku.edu.cn}{\nolinkurl{jyh2714225423@stu.pku.edu.cn}} and
    \href{mailto:zhanghengrui26@stu.pku.edu.cn}{\nolinkurl{zhanghengrui26@stu.pku.edu.cn}}.}
  \enspace Hengrui Zhang\textsuperscript{1}\footnotemark[1]
  \enspace Wenxuan Wu\textsuperscript{2}\thanks{Work done during an internship at Peking University. Email:
    \href{mailto:tico2214215527@stu.xjtu.edu.cn}{\nolinkurl{tico2214215527@stu.xjtu.edu.cn}}.}
  \enspace Jiaqi Zeng\textsuperscript{1}\thanks{Emails:
    Jiaqi Zeng (\href{mailto:jiaqi.zeng@stu.pku.edu.cn}{\nolinkurl{jiaqi.zeng@stu.pku.edu.cn}})
    and Zixin Zhang (\href{mailto:zixinzhang682@gmail.com}{\nolinkurl{zixinzhang682@gmail.com}}).}
  \enspace Zixin Zhang\textsuperscript{1}\footnotemark[3]
  \enspace Kun Yuan\textsuperscript{1}\thanks{Corresponding author. Email:
    \href{mailto:kunyuan@pku.edu.cn}{\nolinkurl{kunyuan@pku.edu.cn}}.}\\[0.5em]
  \normalsize\textsuperscript{1}Peking University
  \qquad
  \normalsize\textsuperscript{2}Xi'an Jiaotong University
}
\date{}
\hypersetup{pdfauthor={Yihan Jin, Hengrui Zhang, Wenxuan Wu, Jiaqi Zeng, Zixin Zhang, Kun Yuan}}
\begin{document}

\makeatletter
\begingroup
\let\originalauthorsymbol\@fnsymbol
\renewcommand*{\@fnsymbol}[1]{%
  \ifnum#1=2 \originalauthorsymbol{4}%
  \else\ifnum#1=4 \originalauthorsymbol{2}%
  \else\originalauthorsymbol{#1}%
  \fi\fi
}
\maketitle
\endgroup
\makeatother

\begin{abstract}
Modern large-scale optimization problems increasingly involve matrix- or tensor-valued parameters whose sheer dimension makes full-space optimization computationally impractical. This primer develops a matrix-oriented subspace optimization framework that updates a high-dimensional parameter matrix through low-dimensional updates: at each iteration, the framework selects a low-dimensional subspace, solves a small subproblem within it (exactly or approximately), and lifts the resulting low-dimensional update back to the full parameter matrix. By treating subspace selection as a modular component and decoupling it from the downstream optimizer that solves the reduced subproblem, the framework accommodates random, greedy, power-iteration, and coordinate-aligned subspaces while preserving enough of the original problem structure to reuse established solvers. We instantiate this perspective across stochastic, zeroth-order, second-order, constrained, federated, decentralized, bilevel, and large-language-model training settings. Through derivations and numerical examples, we show that low-dimensional matrix subspace updates can reduce computation, lower the memory footprint of large-scale model training, and cut communication overhead in distributed optimization, all while attaining convergence behavior comparable to full-rank baselines and substantially outperforming coordinate-descent-type variants. The primer is intended as an accessible guide for practitioners seeking computation-efficient and flexible tools for large-scale matrix optimization.
\end{abstract}

{
\singlespacing
\tableofcontents
}

\section{Introduction}
\subsection{Motivation}

Subspace optimization is a classical and well-established technique for reducing the computational cost of solving large-scale problems. By restricting the search variables to a low-dimensional subspace, it replaces each original high-dimensional update with a much smaller subproblem. Reference \cite{liu2021subspace} surveys a broad range of subspace optimization methods, including block coordinate descent, active-set methods, limited-memory schemes, and sketching-based approaches. It also demonstrates their applications across diverse settings such as unconstrained optimization, stochastic optimization, sparse optimization, and constrained optimization, highlighting the richness and broad applicability of the subspace optimization paradigm.


However, in the era of artificial intelligence and big data, subspace optimization faces several new challenges that are less prominent in its classical settings. (\romannumeral 1) Modern deep neural networks, especially large language models, involve extremely large matrix- and tensor-valued parameters rather than moderate-dimensional vector variables. At this scale, constructing subspaces from historical iterates, gradients, or optimizer states can introduce a prohibitive memory burden, making classical subspace construction strategies difficult to apply directly. (\romannumeral 2) The objective functions arising in modern machine learning are typically highly nonconvex. Consequently, convergence analyses can no longer rely on convexity, global optimality of reduced subproblems, or deterministic descent properties that are often available in classical settings. Instead, they must establish stationarity-type guarantees under weaker assumptions such as smoothness, lower-boundedness, and controlled stochastic errors. (\romannumeral 3) Large-scale learning problems are usually stochastic rather than deterministic, and exact gradient information is often unavailable. Therefore, both the subspace construction and the reduced optimization procedure must be robust to gradient noise, sampling errors, and mini-batch effects. (\romannumeral 4) Downstream optimization algorithms in modern machine learning have evolved from simple gradient-type methods to more sophisticated momentum-based and adaptive algorithms. Although subspace ideas can in principle be combined with such mechanisms, a systematic framework that supports stochastic momentum or adaptive solvers under matrix-structured, memory-constrained, and nonconvex regimes remains underdeveloped.

Among these challenges, the transition from vector-valued optimization to matrix- and tensor-valued optimization represents a particularly fundamental shift. Classical subspace methods are often formulated for vector variables, whereas modern machine learning models naturally organize their trainable parameters as structured matrices and tensors, such as embedding matrices, attention projection matrices, feed-forward weight matrices, and higher-order tensors. Although these parameters can always be vectorized mathematically, a flat-vector representation may obscure their row-column organization, spectral patterns, low-rank structures, and layer-wise relationships. This shift is further amplified by the rapid growth of large language models, where memory consumption, computational cost, and communication overhead have become central bottlenecks in training and fine-tuning. Motivated by these observations, this primer adopts matrix variables as the basic modeling unit, while viewing tensor-valued parameters as natural extensions. We advocate a matrix-oriented view of subspace optimization, in which optimization problems are formulated directly with matrix variables and low-dimensional updates are performed within structure-preserving matrix subspaces.

At the same time, many existing subspace methods are designed by tailoring the subspace construction to a particular problem structure or a specific downstream algorithm. Such specialization is often powerful, but it can also make the underlying algorithmic ideas difficult to transfer to other settings. This primer therefore presents subspace optimization as a reusable design pattern for incorporating dimension reduction into downstream optimization algorithms. Building on the randomized subspace formulation developed in \cite{chen2025memory} for memory-efficient LLM pre-training, we adopt a broader modular view in which subspace construction and the downstream solver are treated as separate design choices. At each iteration, one selects or samples a low-dimensional subspace, formulates a reduced problem within this subspace, solves it using an appropriate domain-specific optimization method, and lifts the resulting update back to the original high-dimensional space. Different instantiations of this recipe mainly differ in how the subspace is constructed, how the reduced problem is defined, and which solver is used. This modular perspective allows existing optimization tools to be reused across different domains while reducing the effective dimension, memory footprint, computational cost, or communication overhead of the original problem. Through this unified recipe, our goal is to help readers quickly understand the core idea of subspace optimization, master its basic theoretical principles, and apply it flexibly to a broad range of downstream optimization problems.

This paper is intended as a primer for engineers, scientists, and students who may not have specialized training in optimization. Rather than developing a separate convergence theory for every possible combination of subspace selection and downstream solver, we focus on explaining the main ideas, practical principles, and empirical behavior of subspace-based optimization in an accessible manner. We illustrate how this idea can be adapted to zeroth-order optimization, second-order optimization, nonsmooth optimization, constrained optimization, bilevel optimization, stochastic optimization, distributed optimization, and federated learning. Through these examples, we aim to help readers understand how to choose or sample a subspace, formulate the corresponding reduced problem, select a compatible downstream solver, and evaluate the resulting method in practice.

\subsection{Problem Formulations}
Consider the following unconstrained, matrix-parameterized optimization problem:
\[
    \min_{W \in \mathbb{R}^{m \times n}} f(W),
\]
where $W \in \mathbb{R}^{m \times n}$ denotes the high-dimensional model parameter matrix, such as the weight matrix of a hidden layer in a neural network. When $n=1$, $W$ reduces to a vector, thereby recovering the standard vector-parameterized setting. Here, $f: \mathbb{R}^{m \times n} \to \mathbb{R}$ is a continuous objective function. Unlike flattening $W$ into a vector of dimension $mn$, the matrix formulation preserves its row-column organization, spectral patterns, and potential low-rank structure.

In particular, this problem formulation is compatible with modern deep learning problems under minor modifications. Specifically, the parameter is extended from a single matrix to a collection of layer-wise matrices $\boldsymbol{W} = \{W_{\ell}\}_{\ell=1}^L$, with each $W_{\ell} \in \mathbb{R}^{m_{\ell} \times n_{\ell}}$ representing the parameter matrix of the $\ell$-th layer. $f(\boldsymbol{W}) = \mathbb{E}_{\xi \sim \mathcal{D}}[F(\boldsymbol{W}; \xi)]$ is a stochastic objective, where $F(\boldsymbol{W}; \xi)$ is the loss function for a data sample $\xi$ drawn from a distribution $\mathcal{D}$.

\subsection{Subspace Optimization Framework}
To address the memory challenges arising in high-dimensional optimization problems, especially those encountered in large-scale model training, \cite{chen2025memory} proposed the random subspace optimization (RSO) framework. Motivated by this line of work, we consider a more general \emph{matrix subspace optimization} framework. Unlike the framework in \cite{chen2025memory}, which focuses specifically on random subspaces, the present framework does not restrict the subspace selection strategy to the random setting. Instead, at each iteration, the subspace may be either sampled or computed according to a prescribed selection mechanism. This generalization decouples the subspace selection mechanism from the downstream optimization solver, thereby substantially enlarging the space of admissible algorithmic combinations and making the framework more flexible for different problem structures. In particular, instead of updating the full high-dimensional parameter matrix $W$ directly at each iteration, the framework decomposes the original problem into a sequence of lower-dimensional subproblems, which can then be solved by existing optimization algorithms in the reduced space. This structured dimensionality reduction is the key to its efficiency.

The general iterative update rule for the subspace optimization framework is as follows:
\begin{subequations}
    \label{eq:so-general}
    \begin{align}
       &B^k \approx \underset{B}{\arg\min} \left\{ f(W^k + P^k B) + \frac{1}{2 \eta^k} \|B\|^2_F \right\}, \label{eq:subproblem} \\
       &W^{k+1} = W^k + P^k B^k, \label{eq:update}
    \end{align}
\end{subequations}
where:
\begin{itemize}
    \item $W^k \in \mathbb{R}^{m \times n}$ is the full parameter matrix at iteration $k$.
    \item $P^k \in \mathbb{R}^{m \times r}$ is a columnwise-orthogonal matrix whose columns span the selected subspace at iteration $k$, where the subspace dimension satisfies $r \ll m$. Here, $P^k$ may be randomly sampled or deterministically computed according to a given subspace selection rule.
    \item $B \in \mathbb{R}^{r \times n}$ denotes the coordinates of the optimization variable in the selected subspace, and the optimization is carried out over this much smaller matrix variable.
    \item The subproblem \eqref{eq:subproblem} is solved for $B^k$. Notably, an exact solution is not required; an approximate solution, often obtained by running a few steps of a standard optimizer, is sufficient.
    \item The term $\frac{1}{2 \eta^k} \|B\|^2_F$ is a proximal term that regularizes the subproblem and helps ensure convergence.
\end{itemize}

\begin{algorithm}[t]
    \caption{Subspace Optimization with Gradient Descent (SO-GD)}
    \label{alg:so-gd}
    \textbf{Input:} $W^{0}$; outer steps $K$; proximal coefficients $\{\eta_k\}$; maximum inner steps $T_{\max}$; inner step sizes $\{\alpha_t\}$.\\
    \textbf{Output:} $W^{K}$.
    \begin{algorithmic}[1]
    \For{$k = 0,1,\dots,K-1$}
        \State Sample or compute a subspace matrix $P^{k}$.
        \State Initialize $B \gets 0$.
        \For{$t = 0,1,\dots,T_{\max}-1$}
            \State Compute gradient at current $B$: $g_B \gets {P^{k}}^{\!\top}\,\nabla f\!\left(W^{k}+P^{k}B\right) \;+\; \tfrac{1}{\eta_k} B$
            \State Update $B$ by $B \gets B - \alpha_t \, g_B$
            \State \textbf{if} {$B$ converges} \textbf{then} $\tilde{B}^{k} \gets B$ and \textbf{break}
        \EndFor
        \State Update $W$ by $W^{k+1} \gets W^{k} + P^{k}\tilde{B}^{k}$
    \EndFor
    \end{algorithmic}
\end{algorithm}

The subspace optimization framework is highly flexible and can be viewed  as a meta-framework that interfaces subspace selection strategies with downstream optimization algorithms for solving the reduced subproblem \eqref{eq:subproblem}. Because the choice of subspace is decoupled from the choice of solver, one may combine different subspace construction mechanisms with different optimization dynamics according to the structure of the target problem. Below, we present two simple instantiations of this framework.

\begin{example}
    Subspace Optimization with Gradient Descent (GD) is a simple instantiation that uses gradient descent to solve the reduced subproblem \eqref{eq:subproblem}. The algorithm is given in Algorithm \ref{alg:so-gd}.

    Compared with full GD performed in the original high-dimensional parameter space, this subspace-based variant only needs to maintain the low-dimensional variable $B$ together with the subspace matrix $P^k$, rather than a full gradient in the ambient space. As a result, the total storage associated with the reduced representation is $O(r(m+n))$, which is substantially more memory-efficient than the $O(mn)$ cost in the full-space setting. More importantly, this example illustrates that once a subspace has been selected, a standard downstream optimizer can be directly applied without modifying the overall framework.
\end{example}

\begin{example}
    Nesterov's Accelerated Gradient (NAG) \cite{nesterov1983method} is a technique for accelerating gradient-based optimization. Subspace Optimization with NAG is obtained by applying Nesterov's acceleration method to solve the reduced subproblem \eqref{eq:subproblem}. The algorithm is given in Algorithm \ref{alg:so-nag}.

    \begin{algorithm}[t]
        \caption{Subspace Optimization with Nesterov's Accelerated Gradient (SO-NAG)}
        \label{alg:so-nag}
        \textbf{Input:} $W^{0}$; outer steps $K$; proximal coefficients $\{\eta_k\}$; maximum inner steps $T_{\max}$.\\
        \textbf{Output:} $W^{K}$.
        \begin{algorithmic}[1]
        \For{$k = 0,1,\dots,K-1$}
            \State Sample or compute a subspace matrix $P^{k}$.
            \State Initialize $B \gets 0$, $Y \gets B$.
            \For{$t = 0,1,\dots,T_{\max}-1$}
                \State Compute gradient at the lookahead point $Y$: $g_Y \gets {P^{k}}^{\!\top}\,\nabla f\!\left(W^{k}+P^{k}Y\right) \;+\; \tfrac{1}{\eta_k} Y$
                \State $B_{\text{new}} \gets Y - \alpha \, g_Y$
                \State $Y \gets B_{\text{new}} + \beta \left(B_{\text{new}} - B\right)$
                \State $B \gets B_{\text{new}}$
                \State \textbf{if} {$B$ converges} \textbf{then} $\tilde{B}^{k} \gets B$ and \textbf{break}
            \EndFor
            \State Update $W$ by $W^{k+1} \gets W^{k} + P^{k}\tilde{B}^{k}$
        \EndFor
        \end{algorithmic}
    \end{algorithm}

Compared with full NAG in the original high-dimensional parameter space, the subspace-based version only needs to maintain the low-dimensional auxiliary variable $Y \in \mathbb{R}^{r \times n}$ as the optimizer state, which is more memory-efficient than storing optimizer states in the full space. This example further demonstrates the flexibility of the framework: the same subspace selection mechanism can be paired with different downstream solvers, while different subspace construction strategies can also be integrated with the same solver.

\end{example}

These examples highlight that the proposed subspace optimization framework can accommodate more sophisticated optimization dynamics while retaining the computational and memory advantages of reduced-dimensional optimization. More importantly, by decoupling subspace selection from the downstream solver, the framework provides a unified and extensible interface for combining diverse subspace construction strategies with a broad range of optimization algorithms.

\section{Related Work}
\subsection{Direct Subspace in Traditional Optimization}
Subspace methods in traditional optimization long predate recent low-rank approaches. Their core idea is not generic dimension reduction, but to restrict each update to a carefully chosen low-dimensional space while preserving the original problem class so that a specialized solver can still be reused. This perspective is closely aligned with our Subspace Matrix Optimization.

A first line of work constructs subspaces from historical iterates, gradients, residuals, or search directions. Early analyses of conjugate-gradient-type methods showed that their updates can be understood as minimization over a very small subspace spanned by the current gradient and previous search information, thereby revealing that many algorithmic variants differ primarily in how coefficients are chosen within this restricted space \cite{yuan1995subspace}. Subsequent developments enlarged this space by incorporating additional gradients, iterate differences, and historical directions, and then performing local quadratic or conic minimization in the resulting low-dimensional model space \cite{karimi2012conjugate}. Sequential subspace optimization methods further advanced this idea by combining gradient information, step differences, and optional second-order directions, yielding memory-efficient yet effective schemes for large-scale optimization \cite{narkiss2005sesop,zibulevsky2008sesop}. Closely related extrapolation-based methods replaced purely Markovian updates by structured combinations of previous solutions or residuals~\cite{pulay1980convergence}, showing that historical information can be exploited through explicit coefficient fitting within an iterative subspace \cite{csaszar1984geometry,chen2019convergence}. Similar principles also appear in acceleration methods that determine coefficients by minimizing residuals over history-based subspaces; in linear settings, these methods exhibit strong connections to Krylov-space algorithms such as GMRES \cite{walker2011anderson}. More broadly, this line of research provides a unifying viewpoint in which conjugate-gradient and limited-memory quasi-Newton methods can both be interpreted as iterated-subspace minimization procedures \cite{conn1996iterated}.

A second line of work uses low-dimensional subspaces to preserve local model or curvature information. Sketching-based methods, for example, project the search or constraint structure into a smaller space~\cite{gower2015randomized,gower2019adaptive} while retaining essential first- or second-order information \cite{necoara2021randomized}.
Trust-region and reduced-Hessian methods follow a similar philosophy. Rather than solving the full-dimensional model directly, they restrict the step computation to a carefully chosen low-dimensional subspace \cite{byrd1988approximate}, typically spanned by gradients, displacement vectors, gradient differences, or approximate curvature directions \cite{wang2004subspace,wang2006subspace}.
Limited-memory and reduced-Hessian techniques further suggest that useful curvature information can often be captured in a compact subspace \cite{liu1989limited,gill2001reduced}. Later variants make this idea more scalable by limiting the reduced-Hessian dimension or combining compact quasi-Newton models with trust-region techniques, while treating the remaining directions with simpler approximations \cite{gill2003limited,burdakov2017efficiently}.
These developments show that subspace methods are not merely computational shortcuts; rather they provide a way to preserve the essential local optimization structure in a lower-dimensional representation.

The role of subspaces becomes more explicit in constrained optimization, where they must encode not only descent information but also feasibility and local geometry. For bound-constrained problems, subspaces have been constructed from scaled gradient and approximate Newton directions so that the update process respects boundary geometry while remaining computationally tractable \cite{branch1999subspace}. In composite-step sequential quadratic programming, a quasi-normal step improves feasibility and a tangential step reduces the objective model \cite{heinkenschloss2014matrix}. Active-set methods follow a related principle by first identifying an approximately active constraint set and then solving an equality-constrained quadratic subproblem within the corresponding reduced space \cite{byrd2011active}. These examples make clear that, in constrained settings, the subspace is not merely a tool for efficiency, but also a mechanism for encoding feasibility, geometry and problem structure.

Compared with these classical subspace methods, Subspace Matrix Optimization is motivated by a different level of abstraction. Traditional methods typically design subspaces for a specific solver family, such as conjugate-gradient, trust-region, SQP, and Newton methods, etc. By contrast, our framework targets matrix-variable optimization problems more broadly and emphasizes a structure-preserving reformulation: the reduced problem remains in the same optimization domain after replacing the update by $X + PB$. This makes it possible to reuse existing domain-specific solvers rather than designing a new one from scratch. Our distinctive contribution is therefore not only low-dimensional search, but a general matrix-oriented interface between subspace selection and structured optimization.

\subsection{Matrix-Based Optimization}
While the classical subspace constructions largely developed in vector-space settings, modern machine learning problems increasingly require a matrix-oriented viewpoint.
Unlike approaches that flatten all parameters into a single vector and then apply generic Euclidean updates, matrix-based optimization treats the matrix itself as the primitive optimization object. The following will explain in three perspectives why it is fruitful. 

First, many matrix-parameterized learning problems exhibit an intrinsic low-rank or low-dimensional subspace bias. In deep matrix factorization, increasing depth strengthens gradient descent's implicit preference for low-rank solutions \cite{arora2019implicit}. Similar low-rank simplicity bias has also been observed for nonlinear deep networks through the effective rank of learned feature Gram matrices \cite{huh2021lowrank}. At the level of optimization dynamics, the first gradient step in certain two-layer models already concentrates on a rank-one spike \cite{ba2022highdim}, while more realistic spiked data models lead to gradient matrices that are well approximated by a small number of dominant spectral components \cite{sonthalia2025low}. More recently, Jaiswal et al.\ \cite{jaiswal2025from} connect the stabilization of gradient subspaces to the eventual low-rankness of different LLM weight matrices. Taken together, these results support the view that low-rank or subspace updates are not merely heuristic compression devices, but often align with the intrinsic optimization dynamics of matrix variables.

Second, matrix-variable problems often possess richer non-Euclidean structures, invariances, and quotient geometries that disappear after vectorization. A canonical example is low-rank matrix completion: the factorization $X=UW$ is invariant under the transformation $(U,W)\mapsto (UM,M^{-1}W)$ for any invertible matrix $M$, indicating that the true search space is better described by subspaces or fixed-rank manifolds than by a flat Euclidean parameterization. This observation has motivated Grassmannian formulations for matrix completion \cite{dai2011geometric,keshavan2009gradient,boumal2011rtrmc}, as well as optimization directly on the manifold of fixed-rank matrices \cite{vandereycken2013low}. Related ideas also appear beyond matrix completion, for example in optimization over fixed-rank positive semidefinite matrices for metric learning and related problems \cite{meyer2011regression}. These works show that once the variable is recognized as a matrix, the relevant optimization structure is often governed by row and column spaces, rank constraints, and positive-semidefinite geometry rather than by independent coordinates.

Third, these matrix-specific structures induce a richer family of matrix-native optimizers. On constrained matrix spaces, one can design geometry-preserving optimizers that update in tangent spaces and retract back to the manifold, ranging from classical Riemannian first- and second-order methods \cite{absil2008optimization} to more specialized variants such as Cayley-transform-based stochastic optimizers on the Stiefel manifold \cite{li2020efficient}, Riemannian adaptive methods \cite{becigneul2019riemannian}, and Riemannian coordinate descent on matrix manifolds \cite{han2024riemannian}. Even when the matrix variable is unconstrained, the optimizer can still exploit its two-sided row/column structure. For example, K-FAC approximates layerwise curvature by Kronecker factors \cite{martens2015optimizing}, Shampoo builds matrix- or tensor-wise preconditioners along each mode \cite{gupta2018shampoo}, and Adafactor compresses second-moment statistics into row and column accumulators \cite{shazeer2018adafactor}. More recently, Muon orthogonalizes matrix updates so as to directly manipulate their singular structure instead of treating them as unrelated scalar coordinates \cite{jordan2024muon}. Classical spectral methods such as singular value thresholding and singular value projection \cite{cai2010singular,jain2010guaranteed} further demonstrate that, for matrix problems, the optimizer itself may naturally act on singular values or low-rank projections. In this sense, matrix-based optimization is not simply vector optimization in disguise: the matrix representation changes both the geometry of the search space and the design space of effective optimizers. 

These observations naturally raise a further question: how should one organize subspace updates once the optimization variable is treated natively as a matrix rather than as a flattened vector?
Compared with the methods above, the proposed subspace optimization framework is best viewed as a complementary matrix-based meta-framework that does not commit to a single geometry, factorization, or preconditioner. Instead, it uses low-dimensional matrix subspaces as a generic computational primitive and emphasizes that the subspace selection mechanism can be decoupled from the downstream optimizer.

\subsection{Subspace Methods in Diverse Optimization Settings}
Building on the classical principles reviewed above, subspace methods have been developed in a broad range of modern optimization settings. These works do not necessarily share a single update formula or a common algorithmic template. In particular, some methods perform additive updates in a reduced space, while others use subspaces to define local models, preserve constraints, construct Krylov approximations, or compress data and curvature information. Therefore, we review existing subspace-related methods according to how the relevant subspace is selected or constructed.

A first class of methods adopts randomized subspaces, where the subspace is generated through stochastic mechanisms. 
Such methods use randomness to reduce the effective search dimension, improve scalability, and sometimes provide additional exploratory behavior in high-dimensional optimization.
In first-order optimization, randomized or adaptive subspace sampling has been used to reduce the effective dimension of the update while retaining the basic structure of gradient-based methods. For example, proximal gradient methods with adaptive subspace sampling demonstrate that restricting updates to low-dimensional random subspaces can lead to effective and scalable first-order schemes \cite{grishchenko2021proximal}.
In general nonsmooth or convex zeroth-order settings, random subspace methods reduce the number of directions that need to be queried and thereby mitigate the dimensional dependence of gradient-free optimization \cite{nozawa2025zeroth}. The concept of subspace zeroth-order optimization also transfers to language-model optimization, where the parameter dimension is much higher. Methods such as LOZO and SubZero introduce low-dimensional random or structured projections to reduce the cost of zeroth-order fine-tuning for large models \cite{chen2024enhancing,yu2025zeroth}. 
Derivative-free subspace methods follow a similar template: each iteration restricts the local model to a randomly selected low-dimensional space, solves the resulting subproblem there, and maps the step back to the original variables. Cartis and Roberts instantiate this idea for nonlinear least-squares problems using scalable Gauss--Newton-type subspace models \cite{cartis2023scalable}, while their later randomized subspace DFO framework uses quadratic models to obtain second-order convergence guarantees \cite{cartis2026randomized}.
Recent randomized subspace DFO methods extend this paradigm to constrained and stochastic settings. CLARSTA uses random-subspace trust-region models for convex-constrained derivative-free optimization \cite{chen2025clarsta}, while STARS and StoDARS provide stochastic trust-region and direct-search variants with convergence and expected-complexity guarantees \cite{dzahini2024stochastic,dzahini2024direct}.
Randomized subspaces have also been explored in second-order and nonconvex optimization. Randomized subspace Newton methods reduce the cost of Newton-type updates by solving curvature-informed subproblems in low-dimensional spaces \cite{gower2019rsn}. 
Randomized second-order methods reduce Hessian costs by computing regularized Newton steps \cite{fuji2025theoretical} or cubic-regularized models \cite{cartis2025random} within low-dimensional subspaces.
Recent random subspace second-order algorithms also study nonconvex settings and improve convergence guarantees by combining random projections with trust-region or regularized second-order models \cite{higuchi2024improving}.
In constrained optimization, the selected subspace must be compatible with feasibility requirements, so random subspace updates are combined with projections or constraint-aware mechanisms \cite{nozawa2023randomized}. 
In distributed or federated environments, randomized subspace updates can be used to reduce communication and memory costs, as in FedKRSO and related methods \cite{yang2026fedkrso}.

A second class of methods relies on structure-aware or deterministically constructed subspaces, where the subspace is derived from geometry, constraints, or the algebraic structure of the problem. 
Krylov and Lanczos methods build iterative subspaces for cubic-regularized Newton \cite{jiang2024krylov} and trust-region or norm-regularized subproblems \cite{al2025extended}, while Lanczos-based bilevel methods use \cite{yang2024lancbio} or recycle \cite{ehrhardt2024efficient} Krylov subspaces to approximate Hessian-inverse hypergradient computations.
In manifold optimization, subspaces are often dictated by local geometry: sequential methods build Stiefel-compatible reduced models \cite{chen2024sequential}, randomized methods exploit fixed-rank manifold structure \cite{gutman2026randomized}, and Riemannian SQO uses tangent-space models with geometry-aware mappings to respect intrinsic constraints \cite{obara2022sequential}.
Composite-step SQP methods use local constraint and gradient information to compute separate steps for feasibility improvement and objective model reduction \cite{heinkenschloss2014matrix}. These methods require carefully selected subspaces to align with specific geometric or structural properties of the problem, often leading to more targeted updates.

A third class of methods learns data- or curvature-adaptive subspaces from sketches, stochastic gradients, or Hessian samples. Double-sketching subspace Newton methods build compact Newton models from sketched data and curvature information \cite{shang2024accelerated}, while subsampled Newton variants use adaptive variance reduction \cite{zhang2022adaptive} or Hessian subsampling \cite{roosta2019sub} to reduce per-iteration cost.

Motivated by these observations, our work abstracts a large portion of these ideas into a matrix-oriented subspace optimization framework.
This framework is inspired partly by the recent Random Subspace Optimization (RSO) framework \cite{chen2025memory}, which provides a compelling perspective for reducing high-dimensional matrix optimization problems to low-dimensional random subspace updates. While inheriting this basic subspace-reduction principle, our formulation does not restrict the subspace selection mechanism to random sampling. Instead, the subspace may be random, structure-aware, or data-adaptive, depending on the problem structure and algorithmic requirements.
In particular, we identify a broad and practically important family of subspace-based algorithms that can be represented, generalized, or approximated through an additive matrix update of the form $W^{k+1}=W^k+P^k B^k$. The key distinction is that our framework explicitly decouples two components that are often intertwined in existing methods: the selection of the subspace and the solution of the reduced subproblem. This decoupling enlarges the algorithmic design space and enables flexible combinations of different subspace construction strategies with different domain-specific solvers.

\subsection{Efficient LLM Training via Subspace Projections}

The broad use of subspace ideas becomes particularly visible in modern large language model training. While the preceding subsection reviews subspace methods across general optimization settings, LLM training provides a concrete and practically important domain in which dimensionality reduction, low-rank structure, and subspace projections are directly tied to memory and computation bottlenecks. A rapidly growing body of recent work shares a common philosophy: exploiting low-dimensional structure in parameters, gradients, blocks, or activations to make large-scale training and fine-tuning more efficient.

One prominent line of work focuses on parameter subspaces.
LoRA \cite{hu2022lora} posits that the update $\Delta W$ during adaptation has low intrinsic rank and parameterizes it as a product of two low-rank matrices.
Subsequent methods such as SLTrain \cite{han2024sltrain} and ReLoRA \cite{lialin2023relora} seek to recover tail information beyond the principal low-rank component, while LORO \cite{mo2025parameter} treats the low-rank parameter as a point on a low-rank matrix manifold and performs Riemannian optimization on that structured space.
These methods illustrate how low-dimensional parameterizations can significantly reduce adaptation cost while remaining tied to matrix structure.

Another important line of work constructs subspaces directly from gradients.
GaLore \cite{zhao2024galore} projects gradients onto a low-rank principal subspace extracted by SVD, thereby compressing optimizer states. Fira \cite{chen2024fira} studies how low-rank updates can recover full-rank behavior through suitable scaling corrections. FRUGAL \cite{zmushko2024frugal} augments a principal low-dimensional gradient subspace with a lightweight treatment of the orthogonal complement, while COAP \cite{xiao2025coap} formulates the projection matrix itself as a dynamic optimization variable. GoLore \cite{he2024subspace} further shows that stochastic subspace selection can be used in place of deterministic spectral projection, thereby providing a concrete example in which random subspace selection is one useful instantiation of a broader subspace-optimization philosophy.

A further line of work operates on activation subspaces. CoLA \cite{liu2025cola} stores low-rank approximations of activations, while CompAct \cite{shamshoum2025compact} and VeLoRA \cite{miles2024velora} use low-rank activations to compute gradients and optimizer states before restoring full-rank quantities only when needed. CR-Net \cite{kong2025cr} goes beyond intra-layer compression by exploiting low-rank cross-layer activation residuals. These methods suggest that the subspace idea is not restricted to parameters or gradients, but can appear throughout the training pipeline.

A related theoretical question concerns how the relevant subspaces evolve over the training trajectory. Flexora \cite{wei2025flexora} and dEBORA \cite{zangrando2025debora} study adaptive rank allocation, SubTrack++ \cite{rajabi2025subtrack++} treats the subspace itself as a trajectory on the Grassmannian manifold, and I3S \cite{zhang2025i3s} develops importance-sampling strategies for dynamic subspace selection during pre-training. Collectively, these works suggest that subspace selection should not be viewed as static: in large-scale training, the quality of a low-dimensional method depends crucially on how the subspace adapts over time.

Overall, existing LLM methods have produced highly effective specialized algorithms, but most remain closely tied to particular structural assumptions, such as low-rank parameterization, spectral gradient compression, activation compression, or coordinate-aligned block updates. However, the proposed matrix-oriented subspace optimization framework is intended as a more general interface: it does not commit to a single form of structure, but instead decouples subspace selection from the downstream optimizer. This perspective helps unify a broad range of existing methods while also leaving room for novel algorithm design.

\section{Comparison with Coordinate Descent}
\label{sec3}

Block Coordinate Descent (BCD) is one of the simplest and most classical instances of subspace-based optimization, and it also serves as one of the most natural baselines for our numerical experiments. From the perspective of the proposed subspace optimization framework, BCD corresponds to the special case in which the selected subspace is aligned with coordinate blocks rather than with more general low-dimensional linear combinations of variables. Therefore, comparing our framework with BCD is useful for clarifying both their connections and their differences.

This relation can be made explicit at the level of the update formula. Let $S_k = \{i_1,\dots,i_r\} \subset \{1,\dots,m\}$ be a selected row index set, and let $E_{S_k} = [e_{i_1},\dots,e_{i_r}] \in \mathbb{R}^{m \times r}$
be the corresponding row-selection matrix. Then the reduced subproblem
\[
B^k \approx \arg\min_{B \in \mathbb{R}^{r \times n}}
\left\{
f\!\left(W^k + E_{S_k} B\right) + \frac{1}{2\eta_k}\|B\|_F^2
\right\},
\]
together with the update
\[
W^{k+1} = W^k + E_{S_k} B^k,
\]
is precisely a proximal row-block update on the rows indexed by $S_k$. Symmetrically, if $S_k \subset \{1,\dots,n\}$ is a selected column index set and $E_{S_k}\in\mathbb{R}^{n\times r}$ is the corresponding column-selection matrix, then
\[
B^k \approx \arg\min_{B \in \mathbb{R}^{m \times r}}
\left\{
f\!\left(W^k + B E_{S_k}^{\top}\right) + \frac{1}{2\eta_k}\|B\|_F^2
\right\},
\]
with the update
\[
W^{k+1} = W^k + B^k E_{S_k}^{\top},
\]
corresponds to a proximal column-block update. Therefore, by choosing a coordinate-aligned subspace matrix, the proposed subspace optimization framework reduces to row-wise or column-wise BCD.

Under this unified viewpoint, the distinctive feature of BCD is not that it is unrelated to subspace optimization, but that it uses a particularly simple and highly structured subspace family: at each iteration, the subspace is spanned by a subset of coordinate directions or blocks, and the reduced subproblem is correspondingly simple. This simplicity gives BCD strong theoretical analyzability and engineering scalability, which explains why it has become one of the foundational tools in large-scale optimization.

A central design issue in BCD is how to select blocks. The classical progression is from cyclic selection, to random selection, to more structured adaptive or greedy rules. Cyclic BCD is easy to implement but can be bottlenecked by the slowest block and may fail to exploit problem structure. Randomized BCD improves robustness and enjoys clean convergence guarantees; early works such as \cite{nesterov2012efficiency,richtarik2014iteration} showed that uniformly random block selection can already achieve strong convergence behavior. Building on this, importance sampling \cite{flamary2015importance} biases selection probabilities according to block sensitivity, while later works introduced adaptive schemes based on current gradients \cite{perekrestenko2017faster}, volume sampling \cite{rodomanov2020randomized}, and safe adaptive importance sampling \cite{stich2017safe}. Markov-chain BCD \cite{sun2020markov} further extend the selection rule to dependent stochastic processes. In contrast, greedy rules like Gauss-Southwell select the block expected to produce the largest immediate improvement, and Nutini et al. \cite{nutini2015coordinate} showed that such rules can converge strictly faster than uniform random selection. More recently, Peng and Yin \cite{peng2024block} demonstrated that properly optimized step sizes can accelerate block gradient descent even without explicit momentum. Taken together, these works show that even within the coordinate-aligned regime, there is already a rich design space around block selection.

Because of flexibility and scalability, BCD has expanded far beyond its original role in separable convex optimization. Classical examples include Platt's SMO algorithm for SVM training \cite{platt1998sequential}, glmnet for sparse generalized linear models \cite{friedman2010regularization}, stochastic dual coordinate ascent \cite{shalev2013stochastic}, and multiconvex optimization for tensor factorization and completion \cite{xu2013block}. The same basic idea has also been extended to constrained optimization \cite{birgin2022block}, manifold optimization \cite{cheung2024randomized,han2024riemannian,peng2023block}, zeroth-order optimization \cite{cai2021zeroth,jing2024asynchronous,park2025unraveling}, and bilevel or multilevel optimization \cite{david2024locally,gupta2023efficient,lauga2025flexible}. In large-model training, the notion of a ``block'' has further evolved from coordinate groups to layers or modules: BAdam \cite{luo2024badam}, BREAD \cite{luoaccelerating}, and BLOCKLLM \cite{ramesh2024blockllm} all use block-wise or layer-wise updates to reduce memory and computation during fine-tuning. This broad applicability explains why BCD is both a strong baseline and a useful reference point for evaluating more general subspace-based methods.

At the same time, the proposed framework differs from BCD in several important ways. First, BCD is restricted to coordinate-aligned or block-aligned subspaces, whereas the present framework allows a much richer family of subspace choices, including random, adaptive, structure-driven, or otherwise computed low-dimensional spaces. Once the subspace family is no longer restricted to coordinate blocks, the method can capture richer linear combinations of variables and cross-block interactions that are invisible to purely coordinate-wise updates. This is particularly important when the optimization problem exhibits low-rank, spectral, or other non-coordinate structures.

Second, BCD is most natural when the objective already has coordinate- or block-separable structure, because cheap block updates can then exploit sparsity and locality directly. The proposed framework is broader: it does not require the useful low-dimensional structure to be aligned with coordinate blocks. Instead, it allows one to choose a subspace that reflects the structure of the problem and then solve the reduced problem with an appropriate downstream optimizer. In this sense, the framework is designed not only to exploit separability, but also to exploit more general low-dimensional linear structure.

Third, our emphasis is on the decoupling between subspace selection and the downstream solver. In BCD, the coordinate-aligned subspace family and the block-update mechanism are often discussed together as a single algorithmic object. In the proposed framework, these two roles are separated: one first chooses a subspace, and then independently chooses how to solve the reduced subproblem. This decoupling enlarges the range of admissible combinations. BCD therefore appears as a special case and baseline within a broader subspace optimization framework, rather than as the conceptual endpoint of subspace-based optimization.

\section{Stochastic Optimization}

Consider a machine learning task of minimizing a loss function over a data distribution. Let \( x \) denote the model parameter that belongs to the feasible set $\mathcal{X}$, and let \( \ell(x; \xi) \) denote the loss on a data sample \( \xi \) drawn from distribution \( \mathcal{D} \). The global objective is to minimize the expected loss:
\begin{equation}
\min_{x \in \mathcal{X}} f(x) := \mathbb{E}_{\xi \sim \mathcal{D}}[\ell(x; \xi)].
\label{eq:expected_loss}
\end{equation}
In practice, we typically have access to a finite set of i.i.d. samples \( \{\xi_i\}_{i=1}^N \) drawn randomly from the true distribution \( \mathcal{D} \), which leads to the empirical risk minimization problem:
\begin{equation}
\min_{x \in \mathcal{X}} f(x) := \frac{1}{N} \sum_{i=1}^N \ell(x; \xi_i).
\label{eq:empirical_loss}
\end{equation}
The traditional gradient descent method requires computing the full exact loss gradient. In the empirical risk minimization problem \eqref{eq:empirical_loss}, the full gradient is given by:
\begin{equation}\label{eq:gradient_average}
\begin{aligned}
& \nabla f(x) = \frac{1}{N} \sum_{i=1}^N \nabla \ell(x; \xi_i),
\end{aligned}
\end{equation}
When \( N \) is large, evaluating the full gradient \( \nabla f(x) \) becomes computationally prohibitive. Stochastic optimization methods address this issue by replacing the full gradient with a stochastic estimate computed from sampled data. The simplest form is Stochastic Gradient Descent (SGD), where at each iteration a single sample \(\xi_i\) is selected and the gradient \(\nabla \ell(x;\xi_i)\) is used as an unbiased estimate of the full gradient.

In practice, a more stable and widely used variant is mini-batch SGD. Given a randomly sampled mini-batch B $\subseteq \{1, 2, \dots, N\}$, the stochastic gradient is computed as :
\begin{equation}\label{eq:batch_gradient}
\begin{aligned}
& g(x) = \frac{1}{|B|} \sum_{i \in B} \nabla l(x; \xi_i).
\end{aligned}
\end{equation}
This estimator reduces the computational cost compared with the full gradient while usually having lower variance than the single-sample estimator.
This principle leads to SGD \cite{robbins1951stochastic, kiefer1952stochastic}, which serves as a foundational approach for large-scale machine learning.
Building on this foundation, several influential variants have been proposed to improve convergence and training stability. Momentum SGD (SGDM) \cite{polyak1964some} accelerates convergence by introducing a momentum term into the parameter updates. Adaptive Gradient Algorithm (AdaGrad) \cite{duchi2011adaptive} adjusts the learning rate according to the accumulated historical gradients, which is particularly effective for sparse or unevenly scaled gradients. Root Mean Square Propagation (RMSProp) \cite{hinton2012divide} follows a similar adaptive learning rate strategy to AdaGrad, but replaces the cumulative gradient statistics with an exponentially weighted moving average of squared gradients, thereby assigning greater weight to recent gradient information. Adaptive Moment Estimation (Adam) \cite{kingma2014adam} combines the momentum mechanism with adaptive per-parameter learning rates, enabling efficient and robust optimization in many deep learning tasks. Stochastic Variance Reduced Gradient (SVRG) \cite{johnson2013accelerating} periodically uses the full gradient to correct the gradient estimate, thereby reducing variance while maintaining the efficiency of stochastic updates. These methods strike different balances between computational efficiency and convergence speed, collectively forming an important foundation for modern large-scale machine learning.

\subsection{Stochastic Gradient Descent (SGD)}
\subsubsection{Preliminary}

The basic update formula of single-sample SGD is given by: 
\begin{equation}\label{eq:sgd_formulation1}
\begin{aligned}
W^{k+1} = W^k - \alpha_k \cdot \nabla f_i(W^k)
\end{aligned}
\end{equation}
where \( W \) denotes the model parameter to be optimized, \( k \) is the current iteration number, \( \alpha_k \) is the learning rate, and \( f_i(\cdot) \) denotes the loss function associated with the \(i\)-th randomly sampled data point.

In practical implementation, the single-sample gradient is often replaced by a mini-batch gradient. For a mini-batch \(B_k\) sampled at iteration \(k\), the update becomes
\begin{equation}
\begin{aligned}
    W^{k+1}=W^k-\alpha_k \frac{1}{|B_k|}\sum_{i\in B_k}\nabla f_i(W^k).
\end{aligned}
\end{equation}
This formulation preserves the stochastic nature of SGD while reducing the variance of the gradient estimate.

\subsubsection{Update Derivation}

To apply \emph{Matrix Subspace Optimization} to the SGD update rule, we consider the subproblem given in \eqref{eq:subproblem}. In this subsection, we use the mini-batch form of SGD to solve the low-dimensional subproblem. Let \( B^{k, 0} = 0 \). For inner iteration \( t \), we draw a random mini-batch \( \xi_{k, t} \) with size \(b\), and compute a stochastic gradient of the objective in \eqref{eq:subproblem} with respect to \( B \). By the chain rule, the stochastic gradient is given by
\begin{subequations}
\begin{align}
g_B^{k, t} &= \frac{1}{b} (P^k)^\top \sum_{i \in \xi_{k, t}} \nabla f_i(W^k + P^k B^{k, t}) + \frac{1}{\eta^k} B^{k, t} \label{eq:rso-gradient}\\
B^{k, t+1} &= B^{k, t} - \alpha_t g_B^{k, t} \label{eq:rso-inner-update}
\end{align}
\end{subequations}

where \(\alpha_t\) is the step size of the inner SGD update, and all other notations are consistent with the definitions above. Upon completing the optimization of the subproblem, we denote the solution obtained from the iteration process as \( \tilde{B}^k \) and update the full-dimensional variable:  \(W^{k+1} = W^k + P^k \tilde{B}^k\). 
The complete procedure is summarized in Algorithm ~\ref{alg:rso-sgd}.

\begin{algorithm}[!t]
        \caption{Subspace Optimization with Stochastic Gradient Descent (SO-SGD)}
        \label{alg:rso-sgd}
        \textbf{Input:} $W^{0}$; batch size b; outer steps $K$; proximal coefficients $\{\eta^k\}$; maximum inner steps $T_{\max}$; inner step sizes $\{\alpha_t\}$.\\
        \textbf{Output:} $W^{K}$.
        \begin{algorithmic}[1]
        \For{$k = 0,1,\dots,K-1$}
            \State Sample or Compute projection matrix $P^k$.
            \State Initialize $B \gets 0$.
            \For{$t = 0,1,\dots,T_{\max}-1$}
                \State Sample mini-batch $\xi$.
                \State Compute stochastic gradient at current $B$: $g_B \gets \frac{1}{b}(P^{k})^{\!\top}\,\sum_{i \in \xi}\nabla F_{i}\!\left(W^{k}+P^{k}B\right) \;+\; \tfrac{1}{\eta^k} B$.
                \State Update $B$ by $B \gets B - \alpha_t \, g_B$.
                \State \textbf{if} {$B$ converges} \textbf{then} $\tilde{B}^{k} \gets B$ and \textbf{break}
            \EndFor
            \State Update $W$ by $W^{k+1} \gets W^{k} + P^{k}\tilde{B}^{k}$.
        \EndFor
        \end{algorithmic}
    \end{algorithm}

\subsubsection{Numerical Example}
We conduct numerical experiments on the MNIST multi-class classification task to evaluate the proposed Subspace SGD method. The experiments are designed to assess two aspects: the influence of the subspace rank on optimization performance, and the comparative behavior of different subspace selection strategies relative to block coordinate descent (BCD) methods.

\begin{figure}[t]
    \center
    \includegraphics{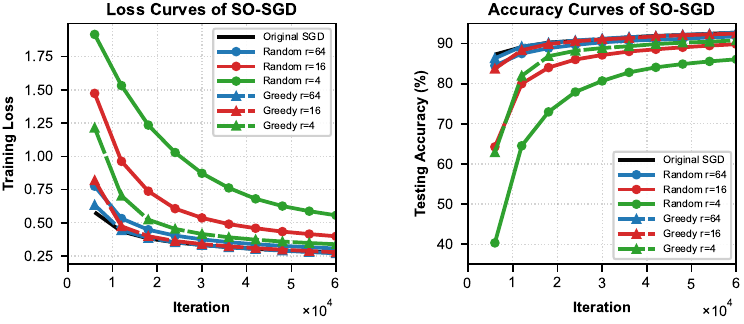}
    \caption{\small Subspace SGD with different ranks. Left: Training loss. Right: Test classification accuracy.
    }
    \label{fig:rso_sgd_ablation}
\end{figure}

\begin{figure}[t]
    \center
    \includegraphics{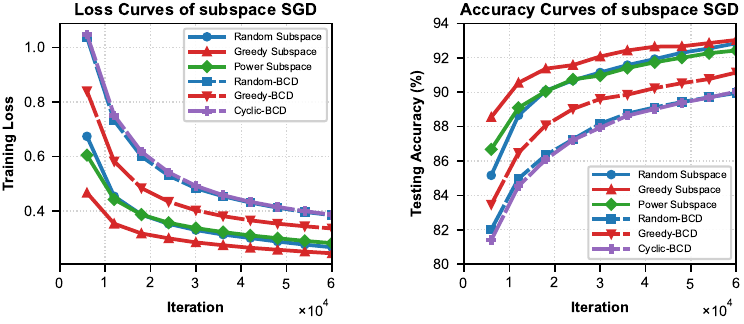}
    \caption{\small Comparison of Subspace Selection Strategy. Left: Training loss. Right: Test classification accuracy. All subspace methods (random subspace, greedy subspace, power subspace) outperform BCD-based methods, since subspace projection better captures low-rank structures and global dimensional dependencies.}
    \label{fig:sgd_rank16}
\end{figure}

\phantomsection
\label{MLP-setting}
The classification model is a two-layer bias-free multi-layer perceptron (MLP). Each \(28 \times 28\) input image is flattened into a 784-dimensional vector. The hidden layer has 128 units, and the output layer has 10 units, matching the 10 classes in MNIST. The first-layer weight matrix therefore has size \(128 \times 784\). For all methods, the mini-batch size is fixed at 10. Unless otherwise specified, all algorithmic hyperparameters are chosen to be the same as those in the corresponding original methods. 
In our settings, we only apply the proposed subspace method to the first-layer parameter matrix.
We report two evaluation metrics throughout the experiments: training loss, which reflects optimization efficiency, and test accuracy, which measures generalization performance. 
The corresponding results can be found in \ref{fig:rso_sgd_ablation} and \ref{fig:sgd_rank16}.

\paragraph{Effect of Subspace Rank}
We first investigate the effect of the subspace rank on the convergence behavior. To this end, we compare Greedy Subspace (SVD) SGD and Random Subspace SGD under different subspace ranks \( r \in \{64, 16, 4\} \). As reported in Figure \ref{fig:rso_sgd_ablation}, the results exhibit a clear monotonic trend: larger subspace ranks consistently lead to improved optimization performance. Notably, when $r = 64$, corresponding to half of the full rank of 128, SO-SGD attains performance close to that of the original full-rank algorithm. This observation suggests that the proposed subspace approach is able to preserve the dominant update directions and retain most of the effective optimization information, even when operating in a substantially reduced-dimensional subspace. Moreover, under the same rank setting, the random subspace strategy generally underperforms the greedy SVD-based method, as the latter constructs the projection subspace by explicitly leveraging the structural information of the gradient matrix.

\paragraph{Ablation on Projection Strategy}
We next compare the proposed subspace-based optimization methods with traditional block coordinate descent (BCD) methods. The core purpose of this experiment is to examine whether subspace methods can achieve better optimization performance than standard row-wise block updates under the exact same per-iteration computational budget.

To ensure a fair comparison, we match the number of free variables updated at each iteration. Specifically, for a rank-$r$ subspace update, the corresponding subproblem is parameterized by a matrix in $\mathbb{R}^{r \times n}$, so the number of decision variables is $rn$. Therefore, in the BCD baselines, we also update exactly $r$ rows at each iteration. Since each selected row contains $n$ variables, the total number of updated variables is likewise $rn$. In this way, the subspace methods and the BCD methods operate under the same update budget.

Given a fixed rank $r$, the descriptions of the six comparison methods used in our experiment are as follows:
\begin{enumerate}\label{subspace-choice}
    \item \textbf{Random Subspace:} 
    The subspace matrices are generated by random orthonormal sampling. Concretely, we first draw Gaussian random matrices and then orthogonalize their columns.

    \item \textbf{Greedy Subspace:} 
    The subspace matrices are chosen as the top-$r$ left singular vectors of the exact gradients of the optimization variable $\bW$. Equivalently, this strategy selects the rank-$r$ subspace that captures the dominant gradient energy at the current iterate.

    \item \textbf{Power Subspace:} 
    This method approximates the top-$r$ singular subspace of the gradient by power iteration. Given a gradient matrix $G \in \mathbb{R}^{m \times n}$ of $\bW$, we start from a random matrix $\Omega \in \mathbb{R}^{n \times r}$, form $Y = G\Omega$, and compute an orthonormal basis of $Y$ via QR factorization, which is the subspace matrix.

    \item \textbf{Greedy-BCD:} 
    This method selects the $r$ rows whose exact gradient row $\ell_2$-norms are the largest. Equivalently, it chooses the corresponding $r$ columns of the identity matrix as the block basis and updates the rows with the strongest first-order signal.

    \item \textbf{Random-BCD:} 
    This method randomly selects $r$ rows to update. Equivalently, it randomly chooses $r$ columns of the identity matrix as the block basis.

    \item \textbf{Cyclic-BCD:} 
    This method cyclically selects $r$ rows to update. Equivalently, it cyclically chooses $r$ columns of the identity matrix as the block basis.
\end{enumerate}
Specifically, we fix the subspace rank to \( r = 16 \) in this experiment. We still report training loss and test accuracy as evaluation metrics throughout this experiment.

As shown in Figure \ref{fig:sgd_rank16}, all projection-based subspace methods, including random subspace, greedy subspace, and power subspace methods, consistently outperform the BCD-based methods across all evaluation metrics. This result suggests that constructing a low-dimensional subspace through projection is more effective than directly selecting a subset of coordinates. One possible explanation is that low-rank projections are better able to capture the intrinsic low-rank structure of the problem. In particular, projection-based subspaces incorporate information from all dimensions in a compressed representation, whereas BCD-based methods update only a subset of coordinates at each iteration. As a result, the latter may discard substantial information from the remaining dimensions and fail to adequately capture the coupling relationships among different coordinates, leading to less effective optimization performance.

\subsection{SGD with momentum (SGDM)}
\subsubsection{Preliminary}

SGDM introduces a momentum mechanism on the basis of SGD, so that the parameter update is governed not only by the current stochastic gradient but also by the accumulated gradient information from previous iterations. In its single-sample form, SGDM can be written as \begin{subequations}
\begin{align}
M^k &= \beta M^{k-1} + \nabla f_i(W^k) \label{eq:momentum_m}\\
W^{k+1} &= W^k - \alpha_k \cdot M^k \label{eq:momentum_theta}
\end{align}
\end{subequations}
where $i$ denotes the index of the randomly sampled data point, \(\beta\) is the decay factor, and \(\alpha_k\) is the learning rate. The momentum term is a linear combination of the previous momentum and the current gradient, and can be interpreted as a moving average of historical gradient information.

In practical implementations, the single-sample gradient is often replaced by a mini-batch gradient. Given a mini-batch \(\mathcal{I}_k\) sampled at iteration \(k\), the mini-batch SGDM update is given by
\begin{subequations}
\begin{align}
M^k
&=
\beta M^{k-1}
+
\frac{1}{|\mathcal{I}_k|}
\sum_{i\in \mathcal{I}_k}
\nabla f_i(W^k), \\
W^{k+1}
&=
W^k-\alpha_k M^k.\label{eq:momentum_theta_minibatch}
\end{align}
\end{subequations}
Compared with the single-sample formulation, the mini-batch version reduces the variance of the stochastic gradient estimate while retaining the computational efficiency of stochastic optimization. In the following algorithm, we adopt this mini-batch formulation of SGDM.

\subsubsection{Update Derivation}
To apply \emph{Matrix Subspace Optimization} in SGDM, we retain the same subspace construction method and subproblem formulation as those in SO-SGD, while replacing the inner optimizer with mini-batch SGDM. Specifically, at each inner iteration, the stochastic gradient of the subproblem is estimated using a randomly sampled mini-batch, and the subspace variable is then updated according to the SGDM rule. The complete procedure for integrating SO with mini-batch SGDM is presented in Algorithm \ref{alg:rso-sgdm}. 

    \begin{algorithm}[t]
        \caption{Subspace Optimization with Momentum SGD (SO-SGDM)}
        \label{alg:rso-sgdm}
        \textbf{Input:} $W^{0}$; batch size b; outer steps $K$; proximal coefficients $\{\eta^k\}$; momentum coefficient $\beta$; maximum inner steps $T_{\max}$; inner step sizes $\{\alpha_t\}$.\\
        \textbf{Output:} $W^{K}$.
        \begin{algorithmic}[1]
        \For{$k = 0,1,\dots,K-1$}
            \State Sample or Compute projection matrix $P^k$.
            \State Initialize $B,M^{0} \gets 0$.
            \For{$t = 0,1,\dots,T_{\max}-1$}
                \State Sample mini-batch $\xi$.
                \State Compute stochastic gradient at current $B$: $g_B \gets \frac{1}{b}(P^{k})^{\!\top}\,\sum_{i \in \xi}\nabla F_{i}\!\left(W^{k}+P^{k}B\right) \;+\; \tfrac{1}{\eta^k} B$.
                \State Update $B$ by $M^{t+1} \gets \beta \cdot M^t + g_B$,\quad $B \gets B - \alpha_t \, M^{t+1}$.
                \State \textbf{if} {$B$ converges} \textbf{then} $\tilde{B}^{k} \gets B$ and \textbf{break}
            \EndFor
            \State Update $W$ by $W^{k+1} \gets W^{k} + P^{k}\tilde{B}^{k}$.
        \EndFor
        \end{algorithmic}
    \end{algorithm}

\subsubsection{Numerical Example}

All numerical experimental settings for SO-SGDM are kept fully consistent with those adopted in \ref{MLP-setting}.

\begin{figure}[t]
    \center
    \includegraphics{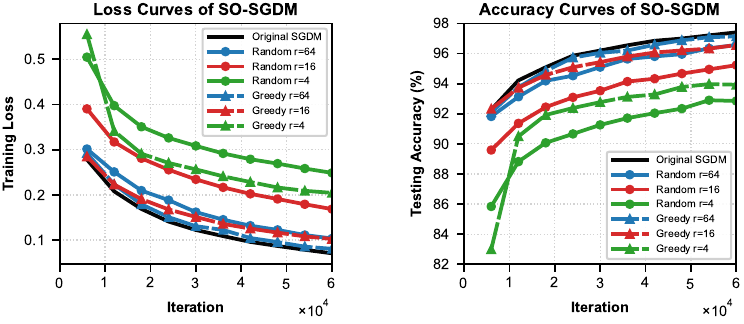}
    \caption{\small Subspace SGDM with different ranks. Left: Training loss. Right: Test classification accuracy.}
    \label{fig:rso_sgdm_ablation}
\end{figure}

\begin{figure}[t]
    \center
    \includegraphics{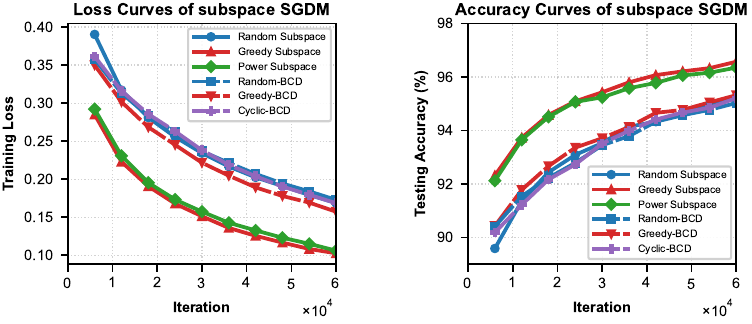}
    \caption{\small Comparison of Subspace Selection Strategy. Left: Training loss. Right: Test classification accuracy. Greedy subspace method and power subspace achieve overall the best performance, while random subspace method performs similarly to Random-BCD and Cyclic-BCD.}
    \label{fig:sgdm_rank16}
\end{figure}

\paragraph{Effect of Subspace Rank}
We first validate the influence of the rank of the subspace on convergence. We compare the Greedy Subspace (SVD) SGDM and Random Subspace SGDM with subspace ranks \( r \in \{64, 16, 4\} \) respectively. Results and insights from Figure \ref {fig:rso_sgdm_ablation} align with those of SO-SGD rank ablation.

\paragraph{Ablation on Projection Strategy}
We further investigate the nuanced differences between specific subspace selection mechanisms and traditional coordinate-wise updates. 
We compare six ``subspace'' selection strategies as defined in \ref{subspace-choice}, and the subspace rank is also fixed to r = 16 in this experiment.
A prominent conclusion from Figure \ref{fig:sgdm_rank16} is that random subspace method achieves comparable performance to BCD baselines, while being outperformed by greedy subspace method and power subspace strategy. This verifies that low-rank projection is superior, which aligns with the conclusions and insights from SO-SGD.

\subsection{Adaptive Moment Estimation (Adam)}
\subsubsection{Preliminary}

Adam incorporates both first-order and second-order moments, which enables it to account for both the direction and the magnitude of the gradient, thereby achieving more efficient optimization.
Based on these two moment estimates, the Adam update can be expressed as follows:
\begin{subequations}
\begin{align}
M^k &= \beta_1 M^{k-1} + (1 - \beta_1) \nabla f_i(W^k) \label{eq:adam_m}\\
V^k &= \beta_2 V^{k-1} + (1 - \beta_2) \left(\nabla f_i(W^k)\right)^2 \label{eq:adam_v}\\
W^{k+1} &= W^k - \alpha_k \cdot \frac{M^k/(1 - \beta_1^k)}{\sqrt{V^k/(1 - \beta_2^k)} + \varepsilon} \label{eq:adam_update}
\end{align}
\end{subequations}

where $M$ and $V$ denote the first-order momentum and the second-order momentum.

\subsubsection{Update Derivation}
For integrating \emph{Matrix Subspace Optimization} with Adam, we keep consistent subspace construction and subproblem formulation with SO-SGD, while only adopting Adam as the optimizer when solving each subproblem.
The complete algorithm for SO-Adam is summarized in Algorithm \ref{alg:rso-adam}.

    \begin{algorithm}[!t]
        \caption{Subspace Optimization with Adam (SO-Adam)}
        \label{alg:rso-adam}
        \textbf{Input:} $W^{0}$; batch size b; outer steps $K$; proximal coefficients $\{\eta^k\}$; first moment coefficient $\beta_1$; second moment coefficient $\beta_2$; small constant $\varepsilon$; maximum inner steps $T_{\max}$; inner step sizes $\{\alpha_t\}$.\\
        \textbf{Output:} $W^{K}$.
        \begin{algorithmic}[1]
        \For{$k = 0,1,\dots,K-1$}
            \State Sample or Compute projection matrix $P^k$.
            \State Initialize $B,M^{0},V^{0} \gets 0$.
            \For{$t = 0,1,\dots,T_{\max}-1$}
                \State Sample mini-batch $\xi$.
                \State Compute stochastic gradient at current $B$: $g_B \gets \frac{1}{b}(P^{k})^{\!\top}\,\sum_{i \in \xi}\nabla F_{i}\!\left(W^{k}+P^{k}B\right) \;+\; \tfrac{1}{\eta^k} B$.
                \State Update $B$ by:
                \State $M^{t+1} \gets \beta_1 \cdot M^t + (1 - \beta_1) \cdot g_B,\quad V^{t+1} \gets \beta_2 \cdot V^t + (1 - \beta_2) \cdot (g_B \odot g_B)$.
                \State $B \gets B - \alpha_t \cdot \left( M^{t+1}/(1-\beta_1^{t+1}) \right) \bigg/ \left( \sqrt{V^{t+1}/(1-\beta_2^{t+1})} + \varepsilon \right)$.
                \State \textbf{if} {$B$ converges} \textbf{then} $\tilde{B}^{k} \gets B$ and \textbf{break}
            \EndFor
            \State Update $W$ by $W^{k+1} \gets W^{k} + P^{k}\tilde{B}^{k}$.
        \EndFor
        \end{algorithmic}
    \end{algorithm}

\subsubsection{Numerical Example}
All numerical experimental settings for SO-Adam are kept consistent with those adopted in \ref{MLP-setting}.

\begin{figure}[t]
    \center
    \includegraphics{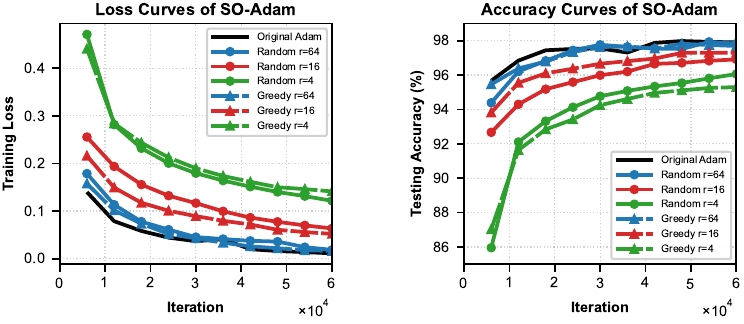}
    \caption{\small Subspace Adam with different ranks. Left: Training loss. Right: Test classification accuracy. }
    \label{fig:rso_adam_ablation}
\end{figure}

\begin{figure}[t]  
    \center 
    \includegraphics{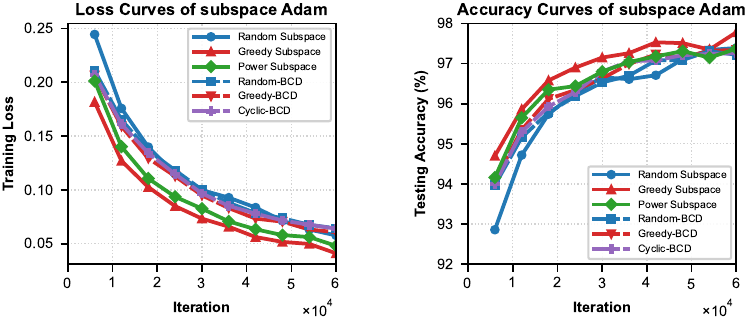}
    \caption{\small Comparison of Subspace Selection Strategy. Left: Training loss. Right: Test classification accuracy. Greedy subspace method and power subspace consistently achieve the best performance, while random subspace method performs comparably to all the BCD variants.}  
    \label{fig:adam_rank16}  
\end{figure}

\paragraph{Effect of Subspace Rank}
We first validate the influence of the rank of the subspace on convergence. We compare the greedy subspace Adam and random subspace Adam with ranks \( r \in \{64, 16, 4\} \) respectively. Results and insights from Figure \ref {fig:rso_adam_ablation} align with those of SO-SGD rank ablation.

\paragraph{Ablation on Projection Strategy}
We further investigate the nuanced differences between specific subspace selection mechanisms and traditional coordinate-wise updates. 
We compare six ``subspace'' selection strategies as defined in \ref{subspace-choice}, and the subspace rank is fixed to r = 16 in this experiment. Similar to the earlier experimental observations for SO-SGDM, it can be observed from Figure \ref{fig:adam_rank16} that the random subspace method attains performance comparable to all the BCD variants, yet is outperformed by greedy subspace and power subspace method, which construct effective subspaces by fully exploiting gradient information.

\subsection{Stochastic Variance Reduced Gradient (SVRG)}
\subsubsection{Preliminary}
SVRG leverages variance reduction by periodically computing the full gradient at a reference point, which is used to correct the noisy stochastic gradient estimate. The algorithm consists of two-layer loops. In the outer loop, the current parameter is saved as a snapshot \(\tilde{W}\), and the full gradient is computed as: \( \nabla f(\tilde{W}) = \frac{1}{N} \sum_{i=1}^N \nabla f_i(\tilde{W})\).
In the inner loop, we randomly select a sample \(i\) to construct the corrected gradient and perform the corresponding parameter update as follows:
\begin{subequations}
\begin{align}
\tilde{\nabla} f_i(W^k) &= \nabla f_i(W^k) - \nabla f_i(\tilde{W}) + \nabla f(\tilde{W}) \label{eq:corrected_gradient}\\
W^{k+1} &= W^k - \alpha_k \cdot \tilde{\nabla} f_i(W^k) \label{eq:update_rule}
\end{align}
\end{subequations}

\subsubsection{Update Derivation}
We maintain the same subspace construction strategy and subproblem formulation as in SO-SGD, and only replace the optimizer used to solve the inner subproblem with the variance-reduced SVRG method. The full procedure of the proposed SO-SVRG is given in Algorithm \ref{alg:rso-svrg}.

    \begin{algorithm}[!t]
        \caption{Subspace Optimization with Stochastic Variance Reduced Gradient (SO-SVRG)}
        \label{alg:rso-svrg}
        \textbf{Input:} Initialization $W^{0}$; batch size b; outer steps $K$; proximal coefficients $\{\eta^k\}$; maximum outer steps of inner SVRG $T_{\max}$; iterations of the inner loop $m$; inner step sizes $\{\alpha_t^{j}\}$.\\
        \textbf{Output:} $W^{K}$.
        \begin{algorithmic}[1]
        \For{$k = 0,1,\dots,K-1$}
            \State Sample or Compute projection matrix $P^k$.
            \State Initialize $\tilde{B} \gets 0$.
            \For{$t = 0,1,\dots,T_{\max}-1$}
                \State Compute full gradient at current $\tilde{B}$: $\mu \gets \frac{1}{N} (P^k)^\top \sum_{i=1}^N \nabla F_i\!\left(W^k + P^k \tilde{B}\right) + \frac{1}{\eta^k} \cdot \tilde{B}$.
                \State Initialize $B \gets \tilde{B}$.
                \For{$j = 0,1,\dots,m-1$}
                    \State Sample mini-batch $\xi$.
                    \State Compute stochastic gradient at current $B$: $g_{B} \gets \frac{1}{b}(P^{k})^{\!\top}\,\sum_{i \in \xi}\nabla F_{i}\!\left(W^{k}+P^{k}B\right) \;+\; \tfrac{1}{\eta^k} B$.
                    \State Compute stochastic gradient at snapshot point  $\tilde{B}$: 
                    \State $\tilde{g_B} \gets \frac{1}{b}(P^{k})^{\!\top}\,\sum_{i \in \xi}\nabla F_{i}\!\left(W^{k}+P^{k}\tilde{B}\right) \;+\; \tfrac{1}{\eta^k} \tilde{B}$.
                    \State Update $B$ by $g_B \gets g_{B} -\tilde{g_B} + \mu$,\quad $B \gets B - \alpha_t^{j} \, g_B$.
                    \State \textbf{if} {$B$ converges} \textbf{then} $\tilde{B}^{k} \gets B$ and \textbf{goto} end.
                \EndFor
                \State $\tilde{B} \gets B$.
            \EndFor
            \State end: Update $W$ by $W^{k+1} \gets W^{k} + P^{k}\tilde{B}^{k}$.
        \EndFor
        \end{algorithmic}
    \end{algorithm}

\subsubsection{Numerical Example}
For SO-SVRG, we employ the identical numerical experimental configurations as those used in \ref{MLP-setting}.

\begin{figure}[t]
    \center
    \includegraphics{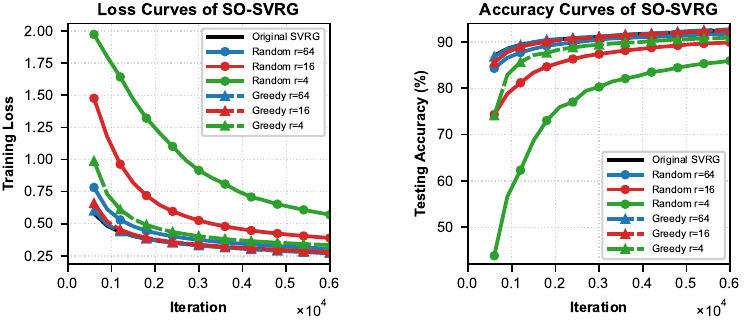}
    \caption{\small Greedy and Random Subspace SVRG with different ranks. Left: Training loss. Right: Test classification accuracy.}
    \label{fig:rso_svrg_ablation}
\end{figure}

\begin{figure}[t]
    \center
    \includegraphics{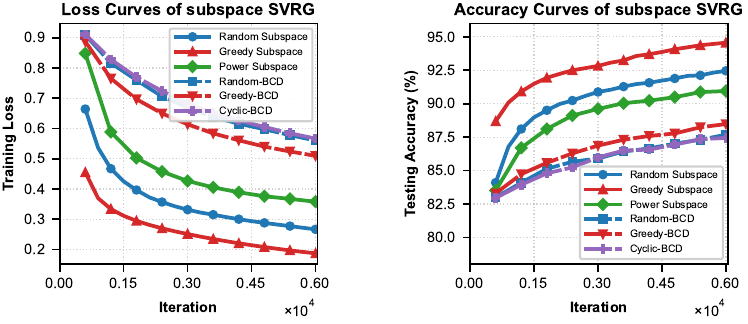}
    \caption{\small Comparison of Subspace Selection Strategy. Left: Training loss. Right: Test classification accuracy. All subspace methods achieve superior performance over BCD variants. In particular, random subspace method even slightly surpasses the performance of the power subspace strategy.}
    \label{fig:svrg_rank16}
\end{figure}

\paragraph{Effect of Subspace Rank}
We first validate the influence of the rank of the subspace on convergence. We compare the greedy subspace SVRG and random subspace SVRG with subspace ranks \( r \in \{64, 16, 4\} \) respectively. Figure \ref{fig:rso_svrg_ablation} shows similar results and insights to the SO-SGD rank ablation.

\paragraph{Ablation on Projection Strategy}
We further investigate the nuanced differences between specific subspace selection mechanisms and traditional coordinate-wise updates. 
We compare six ``subspace'' selection strategies as defined in \ref{subspace-choice}, and the subspace rank is fixed to r = 16 in this experiment. As shown in Figure \ref{fig:svrg_rank16}, the overall trend for SO-SVRG aligns with that of SO-SGD: all subspace-based methods outperform the BCD variants. It is noteworthy that random subspace method even outperforms the power subspace method, which constructs its projection subspace based on current gradient information rather than random sampling. This indicates that the randomness in the selection of the projection matrix provides an advantage for sufficient exploration of the parameter space. By contrast, excessive exploitation of local gradient information may overlook important information contained in other dimensions.

\FloatBarrier
\section{Zeroth-order Optimization}
Zeroth-order (ZO) optimization \cite{nesterov2017random, duchi2015optimal} addresses the problem of minimizing an objective function $f(\bx)$ relying solely on function value evaluations, without access to explicit gradient information. This paradigm is crucial in scenarios such as black-box adversarial attacks, reinforcement learning, and simulation-based engineering design. Recently, ZO methods have been adapted for memory-efficient fine-tuning of LLMs \cite{malladi2023fine} to bypass the storage of activations required for backpropagation.

\subsection{Zeroth-order Stochastic Gradient Descent (ZO-SGD)}
\subsubsection{Preliminary}
A fundamental algorithm in this domain is Zeroth-Order Stochastic Gradient Descent (ZO-SGD) \cite{ghadimi2013stochastic}. ZO-SGD mimics the trajectory of first-order SGD but substitutes the exact gradient with a randomized gradient estimator. The most prevalent estimator is the two-point gradient estimator derived from Gaussian smoothing. Let $\bu$ be a random perturbation vector drawn from a standard normal distribution $\mathcal{N}(0, I)$ and $\mu > 0$ be a smoothing parameter. The gradient is approximated via finite differences:
\begin{equation}
    \hat{\nabla}_{\text{ZO}} f(\bx) = \frac{f(\bx + \mu \bu) - f(\bx - \mu \bu)}{2\mu} \bu.
\end{equation}
At each iteration $t$, ZO-SGD performs: (1) sample $\bu_t \sim \mathcal{N}(0, I)$; (2) query function values $f^+ = f(\bx_t + \mu \bu_t)$ and $f^- = f(\bx_t - \mu \bu_t)$; (3) form the gradient estimate $\hat{\vg}_t = \frac{f^+ - f^-}{2\mu} \bu_t$; and (4) update parameters $\bx_{t+1} = \bx_t - \alpha_t \hat{\vg}_t$, where $\alpha_t$ is the learning rate. Despite the high variance inherent in the gradient estimates, ZO-SGD has been theoretically proven to converge to stationary points for non-convex functions, making it a viable candidate for large-scale model optimization when gradients are inaccessible.

\subsubsection{Update Derivation}
To instantiate the subspace optimization framework in a derivative-free setting, we solve the reduced subproblem with a zeroth-order gradient estimator. At iteration $k$, the framework considers the matrix subspace subproblem:
\begin{equation}\label{eq:rso_sub}
    \min_{\bB} g^k(\bB) := f(\bW^k + \bP^k \bB) + \frac{1}{2\eta^k} \|\bB\|_F^2.
\end{equation}
In the standard first-order setting, one computes $\nabla_{\bB} g^k(\bB) = (\bP^k)^\top \nabla f(\bW) + \frac{1}{\eta^k}\bB$. In the ZO setting, we approximate the gradient of the smooth component using a randomized gradient estimator (e.g., Gaussian smoothing).

Let $\bU$ be a random perturbation matrix drawn from a standard normal distribution $\mathcal{N}(0, I)$ with the same dimensions as $\bB$. The two-point gradient estimator for the parameters within the subspace is given by:
\begin{equation}
    \hat{\nabla}_{\text{ZO}} f(\bW^k + \bP^k \bB) = \frac{f(\bW^k + \bP^k(\bB + \mu \bU)) - f(\bW^k + \bP^k(\bB - \mu \bU))}{2\mu} \bU,
\end{equation}
where $\mu > 0$ is a small smoothing parameter. Consequently, the gradient estimate for the regularized subproblem becomes:
\begin{equation}
    \hat{\bg}_{\bB} = \hat{\nabla}_{\text{ZO}} f(\bW^k + \bP^k \bB) + \frac{1}{\eta^k} \bB.
\end{equation}
This estimator is then used to perform SGD-like updates on $\bB$. The complete procedure is summarized in Algorithm \ref{alg:so-zo}. This approach significantly reduces the dimensionality of the ZO search space from $d$ (dimension of $\bW$) to $r$ (dimension of $\bB$), which is known to improve the convergence rate of ZO methods which typically scale with problem dimension \cite{duchi2015optimal}.

\begin{algorithm}[t]
    \caption{Subspace Optimization with Zeroth-Order SGD (SO-ZO)}\label{alg:so-zo}
    \textbf{Input:} $\bW^0$; outer steps $K$; perturbation $\mu$; proximal coefficients $\{\eta^k\}$; maximum inner steps $T_{\max}$; inner step sizes $\{\alpha_t\}$.\\
    \textbf{Output:} $\bW^K$.
    \begin{algorithmic}[1]
        \For{$k = 0, 1, \dots, K-1$}
            \State Sample or compute a subspace matrix $\bP^k$.
            \State Initialize subspace parameters $\bB \gets 0$.
            \For{$t = 0, 1, \dots, T_{\max}-1$}
                \State Sample mini-batch $\xi$ and perturbation matrix $\bU \sim \mathcal{N}(0, I)$.
                \State Query paired losses $f^{\pm} \leftarrow f(\bW^k + \bP^k(\bB \pm \mu \bU); \xi)$.
                \State Estimate gradient: $\hat{\bg}_{\bB} \leftarrow \frac{f^+ - f^-}{2\mu}\bU + \frac{1}{\eta^k}\bB$.
                \State Update subspace weights: $\bB \leftarrow \bB - \alpha_t \hat{\bg}_{\bB}$.
                \State \textbf{if} {$\bB$ converges} \textbf{then break}
            \EndFor
            \State Update global weights: $\bW^{k+1} \leftarrow \bW^k + \bP^k \bB$.
        \EndFor
    \end{algorithmic}
\end{algorithm}

\subsubsection{Numerical Example}
At a high level, we follow the MNIST classification protocol in Section \ref{MLP-setting}, but the zeroth-order experiment uses a different backbone and data pipeline tailored to black-box optimization. Concretely, we use a bias-enabled MLP-2 with widths $784 \to 256 \to 10$, evaluate on fixed subsets of $10{,}000$ training images and $1{,}000$ test images, and use batch size $512$. The subspace optimization method is applied to all matrix-valued weight blocks; in particular, the dominant matrix variable is the first-layer weight $W_1 \in \mathbb{R}^{256 \times 784}$. As in Section \ref{MLP-setting}, we report training loss and testing accuracy.

\paragraph{Effect of Subspace Rank}
We investigate the trade-off between memory efficiency (determined by rank $r$) and convergence speed. We compare strictly ZO-feasible Random Subspace SO-ZO and Greedy Subspace SO-ZO as a gradient-informed oracle under ranks $r \in \{64, 16, 4\}$, together with the full-rank ZO-SGD baseline. For each matrix block $W_{\ell} \in \mathbb{R}^{m_{\ell} \times n_{\ell}}$, the effective subspace rank is $r_{\ell} = \min\{r, m_{\ell}\}$. Across all rank settings, we keep the same ZO hyperparameters, including the smoothing parameter and step size scaling, so that the comparison isolates the effect of subspace dimension and subspace construction.

\begin{figure}[t]
    \center
    \includegraphics{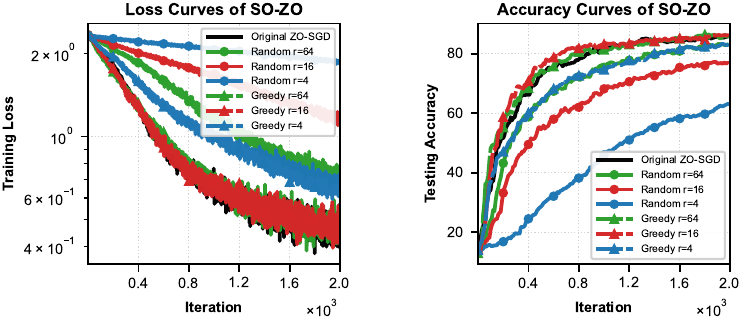}
    \caption{\small Subspace Optimization ZO-SGD with different ranks and subspace construction rules. Left: Training Loss. Right: Test Accuracy. Random Subspace is the strictly ZO-feasible variant, while Greedy Subspace is included as a gradient-informed oracle to diagnose the effect of subspace quality.}\label{fig:so-zo-rank}
\end{figure}

As illustrated in Figure \ref{fig:so-zo-rank}, SO-ZO with larger ranks achieves competitive accuracy compared with the full-rank baseline, validating that optimization can be effectively confined to low-dimensional subspaces even without backpropagating full gradients. The comparison between random subspace and the greedy subspace oracle separates two effects: random subspace tests the strictly memory-efficient black-box setting, whereas the oracle subspace indicates how much improvement is possible when the selected subspace is aligned with first-order information. Reducing the rank to $r=4$ hinders the expressivity of the ZO, leading to higher loss.

\paragraph{Effect of Subspace Selection Strategy}
We next examine how the choice of subspace affects zeroth-order optimization. Following the strategy definitions and same-budget principle in \ref{subspace-choice}, we fix the subspace rank to $r=64$ and compare random subspace with gradient-informed oracle strategies, including greedy subspace, power subspace, and greedy-BCD, against the full-rank ZO-SGD baseline. For each matrix block $W_{\ell} \in \mathbb{R}^{m_{\ell} \times n_{\ell}}$, the effective reduced dimension is $r_{\ell}=\min\{64,m_{\ell}\}$; the coordinate baseline updates $r_{\ell}$ rows, matching the number of free variables in $B_{\ell}\in\mathbb{R}^{r_{\ell}\times n_{\ell}}$.

Compared with the setting in \ref{subspace-choice}, the key difference is that the inner optimization of $\bB$ remains zeroth-order. The greedy subspace, power subspace, and greedy-BCD strategies therefore use the full first-order gradient only to construct $\bP$. This hybrid oracle setup deliberately violates the ZO memory constraints to isolate the quality of the selected subspace from the efficiency of the zeroth-order solver.

\begin{figure}[t]
    \center
    \includegraphics{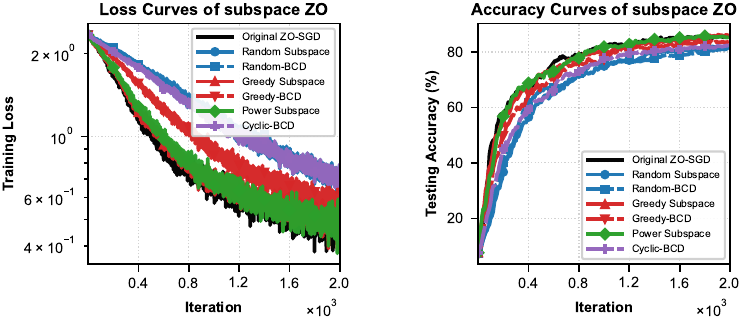}
    \caption{\small Comparison of Subspace Selection Strategy. Left: Training Loss. Right: Test Accuracy. Gradient-informed oracle strategies (Greedy Subspace, Power Subspace) yield faster initial convergence compared to random Gaussian projections, but the Random Subspace eventually catches up.}\label{fig:so-zo-strategy}
\end{figure}

As demonstrated in Figure \ref{fig:so-zo-strategy}, aligning the subspace with the principal components of the current mini-batch gradient (greedy subspace and power subspace) significantly accelerates the initial convergence of SO-ZO, vastly outperforming axis-aligned methods like greedy-BCD. These gradient-informed methods should be interpreted as oracle diagnostics rather than strictly black-box ZO methods, since their subspace construction requires full first-order information. As training progresses, the standard random subspace strategy eventually reaches a comparable loss floor and matching test accuracy. This solidifies a crucial conclusion: while gradient-informed subspaces offer an initial speedup, the pure random subspace approach, which requires zero gradient memory to construct the subspace, remains a highly viable and strictly memory-efficient strategy for black-box zeroth-order optimization.

\section{Second-order Optimization}
Second-order optimization methods \cite{nocedal2006numerical, dennis1996numerical}, such as Newton's method and Quasi-Newton methods, leverage curvature information (Hessian matrix) to achieve faster convergence rates compared to first-order methods. These algorithms are particularly effective for ill-conditioned problems where first-order methods may struggle. However, their application to large-scale deep learning models is limited by the prohibitive cost of computing, storing, and inverting the high-dimensional Hessian matrix \cite{martens2010deep, xu2020second}.

\subsection{Newton's Method}
\subsubsection{Preliminary}
The classical Newton's method updates parameters by minimizing a local quadratic approximation of the objective function $f(\bx)$ at iteration $k$. The Taylor expansion at $\bx^k$ is given by:
\begin{equation}
    f(\bx) \approx f(\bx^k) + \nabla f(\bx^k)^{\top}(\bx-\bx^k) + \frac{1}{2} (\bx-\bx^k)^{\top}\nabla^2 f(\bx^k) (\bx-\bx^k),
\end{equation}
where $\nabla^2 f(\bx^k)$ denotes the Hessian matrix at $\bx^k$. Setting the gradient of this approximation to zero yields the Newton direction $\bd_k$ satisfying
\[
    \nabla^2 f(\bx^k)\bd_k = -\nabla f(\bx^k),
\]
and the update $\bx^{k+1} = \bx^k + \bd_k$. Equivalently, at each iteration $k$, Newton's method: (1) computes the gradient $\bg_k = \nabla f(\bx^k)$ and Hessian $H_k = \nabla^2 f(\bx^k)$; (2) solves the linear system $H_k \bd_k = -\bg_k$ for the search direction, which requires $O(d^3)$ cost in general when using direct solvers, or can be handled by iterative solvers; and (3) updates $\bx^{k+1} = \bx^k + \bd_k$. Despite its quadratic convergence property, the $O(d^3)$ complexity of solving the linear system makes the method impractical when $d$ is large.
\subsubsection{Update Derivation}
For notation clarity, we present the second-order derivation in vector form. The corresponding matrix form follows naturally by replacing the vector variable with a matrix variable and applying the same low-dimensional subspace parametrization. This avoids cumbersome higher-order tensor notation for Hessians while preserving the same subspace-reduction principle.

To enable second-order optimization for large-scale models, we instantiate the subspace optimization framework with Newton updates. Instead of operating on the full parameter space $\mathbb{R}^d$, we perform Newton steps in a significantly lower-dimensional subspace $\mathbb{R}^r$ (where $r \ll d$). At outer iteration $k$, we update the global parameters as $\bx^{k+1} = \bx^k + \bP^k \tilde{\bb}^k$, where $\tilde{\bb}^k$ is obtained by approximately minimizing the regularized subspace objective
\begin{equation}
    \tilde{g}_k(\bb) := f(\bx^k + \bP^k \bb) + \frac{1}{2\eta^k}\|\bb\|^2.
\end{equation}

The subspace Newton method approximates the optimal $\tilde{\bb}^k$ by solving this lower-dimensional optimization problem. The update rule for the inner loop $t$ becomes:
\begin{equation}
    \bb^{k,t+1} \leftarrow \bb^{k,t} -
    \left[
    (\bP^k)^\top \nabla^2 f(\bx^k+\bP^k\bb^{k,t})\bP^k + (\eta^k)^{-1} I
    \right]^{-1}
    \left[
    (\bP^k)^\top \nabla f(\bx^k+\bP^k\bb^{k,t}) + (\eta^k)^{-1}\bb^{k,t}
    \right].
\end{equation}
Crucially, the regularized subspace Hessian is an $r \times r$ matrix, which is much smaller than the full Hessian and reduces the inversion cost from $O(d^3)$ to $O(r^3)$. This makes second-order updates computationally feasible even for large models. The gradients and Hessians in the subspace can be efficiently computed using vector-Jacobian products (VJP) and Hessian-vector products (HVP) without explicitly forming the full Hessian. The complete SO-Newton procedure is outlined in Algorithm \ref{alg:so-newton}.

\begin{algorithm}[t]
    \caption{Subspace Optimization with Newton's Method (SO-Newton)}\label{alg:so-newton}
    \textbf{Input:} $\bx^0$; outer steps $K$; proximal coefficients $\{\eta^k\}$; maximum inner steps $T_{\max}$.\\
    \textbf{Output:} $\bx^K$.
    \begin{algorithmic}[1]
        \For{$k = 0, 1, \dots, K-1$}
            \State Sample or compute a subspace matrix $\bP^k$.
            \State Initialize subspace parameter $\bb \gets 0$.
            \For{$t = 0, 1, \dots, T_{\max}-1$}
                \State Compute regularized subspace gradient: $g_t \gets (\bP^k)^\top \nabla f(\bx^k+\bP^k\bb)+(\eta^k)^{-1}\bb$.
                \State Compute regularized subspace Hessian: $H_t \gets (\bP^k)^\top \nabla^2 f(\bx^k+\bP^k\bb)\bP^k+(\eta^k)^{-1}I$.
                \State Update subspace parameter: $\bb \leftarrow \bb - H_t^{-1} g_t$. \Comment{Requires $O(r^3)$ or iterative solver}
                \State \textbf{if} {$\bb$ converges} \textbf{then break}
            \EndFor
            \State Update global parameters: $\bx^{k+1} \leftarrow \bx^k + \bP^k \bb$.
        \EndFor
    \end{algorithmic}
\end{algorithm}

\subsubsection{Numerical Experiments}
Unlike the MNIST classification setting in Section \ref{MLP-setting}, this section uses a synthetic high-dimensional ridge regression problem with $N=5000$, $D=1024$, and $\lambda=10^{-3}$. This controlled quadratic setting admits a closed-form optimum, allowing us to track both the exact optimality gap and the $R^2$ score. We compare the convergence behavior of SO-Newton against the standard Newton method across different ranks.

\paragraph{Effect of Subspace Dimension}
We investigate how the reduction in subspace dimension $r$ influences the convergence rate towards the optimal solution $f^*$. We evaluate SO-Newton with subspace dimensions $r \in \{256, 64, 8\}$ and compare them with the standard Newton's method. The inner loop performs 2 subspace Newton steps before re-sampling the projection matrix.

\begin{figure}[t]
    \center
    \includegraphics{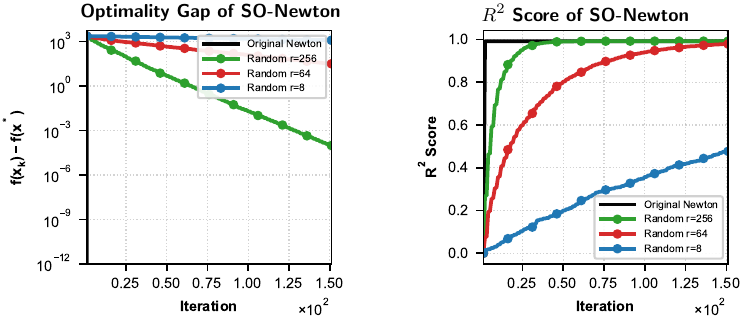} 
    \caption{\small Subspace Optimization Newton (SO-Newton) with different subspace dimensions. Left: Optimality Gap $f(\bx^k)-f^*$. Right: $R^2$ Score. SO-Newton maintains convergence even when operating in significantly lower-dimensional subspaces.}\label{fig:so-newton-rank}
\end{figure}

As depicted in Figure \ref{fig:so-newton-rank}, SO-Newton demonstrates strong convergence properties relative to the full-dimensional baseline. With $r=256$ (25\% of the full dimension), SO-Newton achieves a convergence rate nearly identical to the standard Newton's method, effectively closing the optimality gap within 50 steps. Reducing the subspace dimension to $r=64$ results in a slight slowdown but still maintains a robust convergence profile. Even in the extreme case of $r=8$ ($<$1\% of dimensions), the algorithm consistently descends, albeit at a reduced rate. This confirms that the subspace Newton instantiation can capture useful curvature information within low-dimensional random subspaces, offering a scalable second-order optimization strategy for high-dimensional problems.

\paragraph{Effect of Subspace Selection Strategy}\label{para:so-newton-subspace-selection}
Beyond analyzing the subspace dimension, we explore how different strategies for selecting the low-dimensional subspace affect the convergence of second-order optimization. We use the same comparison logic as in \ref{subspace-choice}, but specialize it to the vector-valued ridge-regression problem. Specifically, we fix the subspace dimension to $r=64$ and compare Random Subspace, Random-CD, Cyclic-CD, and Greedy-CD against the full-dimensional Newton baseline. Since all methods optimize the same vector $w \in \mathbb{R}^{1024}$, this is a same-budget comparison: subspace methods operate in a 64-dimensional subspace, while the coordinate baselines update exactly 64 coordinates per outer step.

In this second-order setting, Random Subspace uses a normalized Gaussian projection, Random-CD and Cyclic-CD use randomly sampled and deterministic coordinate blocks, respectively, and Greedy-CD selects the top-$r$ coordinates with the largest magnitudes in the exact full gradient $\nabla f(w)$ at the current iterate. Greedy-CD is therefore an oracle baseline, since it requires an $O(D)$ exact gradient computation before taking the $O(r^3)$ Newton step.

\begin{figure}[t]
    \center
    \includegraphics{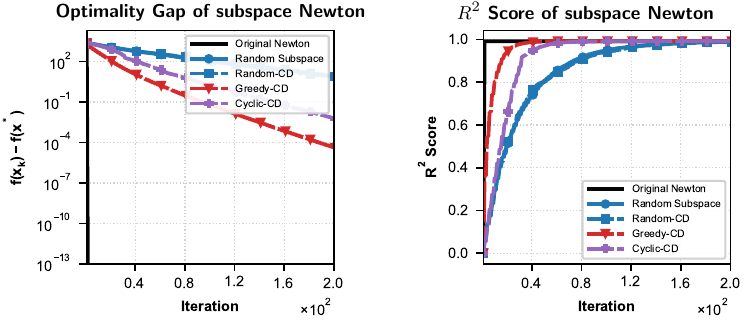}
    \caption{\small Comparison of Subspace Selection Strategy. Left: Optimality Gap. Right: $R^2$ Score. Standard Newton solves the quadratic problem in a single step. Among reduced-dimensional methods, the Greedy-CD oracle is the fastest, followed by the deterministic Cyclic-CD. The Random Subspace approach and Random-CD exhibit similar, slower convergence rates due to the undirected nature of random sampling.}\label{fig:so-newton-strategy}
\end{figure}

As depicted in Figure \ref{fig:so-newton-strategy}, the choice of subspace or coordinate selection significantly impacts the efficiency of the Newton updates. As theoretically expected for a strongly convex quadratic problem, the full-dimensional standard Newton method converges exactly to the optimal solution in a single step. 

Among the dimensionality-reduction methods, the Greedy-CD oracle is the fastest. By explicitly tracking the steepest coordinate directions, namely the coordinates with the largest magnitudes in the exact full gradient $\nabla f(w)$, it accurately targets the most critical dimensions. However, this method breaks the strict computational efficiency of subspace methods, as it mandates the evaluation of the full $D$-dimensional gradient vector to perform the top-$r$ selection at every outer step.

Cyclic-CD emerges as the second-fastest approach. By systematically sweeping through all coordinate blocks, it ensures that every dimension is updated exactly once every $\lceil D/r \rceil$ iterations. This avoids the ``coupon collector's'' delay inherent in random sampling, thereby yielding a steady and reliable linear descent.

Finally, the Random Subspace and Random-CD approaches exhibit similar performance and are the slowest among the tested strategies. Random-CD naturally suffers from the probability of repeatedly updating the same coordinates while neglecting others for extended periods. On the other hand, the random-subspace strategy projects the $D$-dimensional space onto a dense, random $r$-dimensional subspace. For this specific objective, blindly mixing all coordinates via random Gaussian projections does not yield a more effective descent direction than isolating them randomly, causing it to essentially match the convergence rate of uniform random coordinate sampling.

\section{Constrained Optimization}
Constrained optimization represents a ubiquitous challenge within the field of optimization. In recent years, numerous cutting-edge machine learning paradigms have also been formulated within this fundamental framework. Traditional optimization literature provides a wealth of solutions, with representative algorithms including the Augmented Lagrangian Method (ALM), the Alternating Direction Method of Multipliers (ADMM), and the Primal-Dual Hybrid Gradient (PDHG) etc. To enhance the capability of these algorithms in handling high-dimensional problems, this section presents the derivation of update rules integrating \emph{Matrix Subspace Optimization} with these classical methods, alongside comprehensive numerical experiments.

\subsection{Alternating Direction Method of Multipliers (ADMM)}
\subsubsection{Preliminary}
The Alternating Direction Method of Multipliers (ADMM) was originally proposed by Glowinski and Marrocco in \cite{glowinski1975approximation} and independently by Gabay and Mercier in \cite{gabay1976dual}. This algorithm is specifically designed to address optimization problems characterized by a separable objective function subject to linear equality constraints. The standard formulation of such a problem is expressed as:
\begin{equation}\label{eq:admm_formulation1}
\begin{aligned}
& \underset{x, z}{\text{minimize}} & & f(x) + g(z), \\
& \text{subject to} & & A_{1}x + A_{2}z = c.
\end{aligned}
\end{equation}

where variables $x \in \mathbb{R}^n$ and $z \in \mathbb{R}^m$ are coupled by the linear constraint $A_1x + A_2z = c$, and the objective function splits into two convex terms, $f$ and $g$. To solve this, ADMM utilizes the augmented Lagrangian function $\mathcal{L}_{\rho}(x, z, y):=f(x) + g(z) + y^{\top}(A_{1}x + A_{2}z - c) + \frac{\rho}{2}\|A_{1}x + A_{2}z - c\|^{2}$, where $y$ is the dual variable and $\rho > 0$ is a penalty parameter. The algorithm proceeds by iterating through three steps: minimizing $\mathcal{L}_{\rho}$ with respect to $x$, then with respect to $z$, and updating the dual variable $y$. At iteration $k$, the updates are:

\begin{subequations}
\begin{align}
x^{k+1} &:= \underset{x}{\text{argmin}} \ \mathcal{L}_{\rho}(x, z^k, y^k), \\
z^{k+1} &:= \underset{z}{\text{argmin}} \ \mathcal{L}_{\rho}(x^{k+1}, z, y^k), \\
y^{k+1} &:= y^k + \rho (A_{1}x^{k+1} + A_{2}z^{k+1} - c).
\end{align}
\end{subequations}

The theoretical properties and practical utility of ADMM have been studied extensively. In terms of convergence, it has been established that for general convex problems, ADMM converges to an optimal solution with an ergodic convergence rate of $O(\frac{1}{k})$ in the objective function value in \cite{he20121}. Furthermore, under stronger assumptions such as the strong convexity of the objective function or the smoothness of the component functions, the algorithm can achieve a linear convergence rate \cite{deng2016global}.

ADMM's decomposition allows for the separate minimization of $f$ and $g$, making ADMM attractive for distributed optimization and statistical learning where datasets are massive or distributed. The method gained popularity in the modern data science community following the comprehensive survey by Boyd et al. \cite{boyd2011distributed}, which demonstrated its efficacy in Lasso regularization, sparse coding, and consensus optimization.

\subsubsection{Update Derivation}
To apply \emph{Matrix Subspace Optimization} in ADMM's updating rule, we need to transform the problem of full-dimensional variables into the form of subspace variables. The following problem is equivalent to \eqref{eq:admm_formulation1}:

\[
\min_{W_{1}, W_{2}}\quad f(W_{1})+g(W_{2}) + \mathcal{I}\{A_{1}W_{1}+A_{2}W_{2}-C=0\}.
\]

where vector variable $x,z$ are substituted by matrix variable $W_{1},W_{2}$, $\mathcal{I\{\cdot\}}$ is indicator function, i.e.

\[
\mathcal{I}\{AW_{1}+BW_{2}-C=0\} = \begin{cases}
0\qquad &\text{if}\,A_{1}W_{1}+A_{2}W_{2}-C=0,\\
+\infty\qquad &\text{if}\,A_{1}W_{1}+A_{2}W_{2}-C\neq0.
\end{cases}
\]

Now, we substitute the full-dimension variables $W_{1}$, $W_{2}$ by $W_{1}^{k}+P^{k}B_{1}$, $W_{2}^{k}+R^{k}B_{2}$ respectively at outer iteration k from update rule of \emph{Matrix Subspace Optimization} framework, with sampled subspace $P^{k}$ and $R^{k}$, fixed full-dimension variable $W_{1}^{k}$ and $W_{2}^{k}$, and add a regularization term for subspace variable $B_{1}$ and $B_{2}$:

\[
\begin{aligned}
\min_{B_{1},B_{2}}\quad & f(W_{1}^{k}+P^{k}B_{1}) + g(W_{2}^{k}+R^{k}B_{2}) \\
& {}+ \mathcal{I}\{A_{1}(W_{1}^{k}+P^{k}B_{1})+A_{2}(W_{2}^{k}+R^{k}B_{2})-C=0\}
 + \frac{1}{2\eta^{k}}\|B_{1}\|_{F}^{2} + \frac{1}{2\eta^{k}}\|B_{2}\|_{F}^{2}.
\end{aligned}
\]

By eliminating the form of indicator function $\mathcal{I}\{\cdot\}$, we get the subproblem:

\begin{equation}\label{eq:so_admm_formulation1}
\begin{aligned}
& \underset{B_{1}, B_{2}}{\text{minimize}} & & f(W_{1}^{k}+P^{k}B_{1}) + g(W_{2}^{k}+R^{k}B_{2}) + \frac{1}{2\eta^{k}}\|B_{1}\|_{F}^{2} + \frac{1}{2\eta^{k}}\|B_{2}\|_{F}^{2},\\
& \text{subject to} & & A_{1}P^{k}B_{1} + A_{2}R^{k}B_{2} = C-A_{1}W_{1}^{k}-A_{2}W_{2}^{k}.
\end{aligned}
\end{equation}

Now, we can define the augmented Lagrangian function $\mathcal{L}_{\rho}(B_{1}, B_{2}, Y)$ for subproblem \eqref{eq:so_admm_formulation1}:

\begin{equation}
\begin{aligned}
\mathcal{L}_{\rho}^{k}(B_{1}, B_{2}, Y) = &\{f(W_{1}^{k}+P^{k}B_{1}) +\frac{1}{2\eta^{k}}\|B_{1}\|_{F}^{2} \}+\{g(W_{2}^{k}+R^{k}B_{2})+ \frac{1}{2\eta^{k}}\|B_{2}\|_{F}^{2}\},\\
&+<Y, A_{1}P^{k}B_{1} + A_{2}R^{k}B_{2} -\tilde{C}> + \frac{\rho}{2}\|A_{1}P^{k}B_{1} + A_{2}R^{k}B_{2} -\tilde{C}\|_{F}^{2},
\end{aligned}
\end{equation}

where $\tilde{C}=C-A_{1}W_{1}^{k}-A_{2}W_{2}^{k}$. Then follow the update rule of ADMM to solve subproblem \eqref{eq:so_admm_formulation1}:

\begin{subequations}
\begin{align}
B_{1}^{t+1} &:= \underset{B_{1}}{\text{argmin}} \ \mathcal{L}_{\rho}^{k}(B_{1}, B_{2}^{t}, Y^t), \\
B_{2}^{t+1} &:= \underset{B_{2}}{\text{argmin}} \ \mathcal{L}_{\rho}^{k}(B_{1}^{t+1}, B_{2}, Y^t), \\
Y^{t+1} &:= Y^t + \rho (A_{1}P^{k}B_{1}^{t+1} + A_{2}R^{k}B_{2}^{t+1} -\tilde{C}).
\end{align}
\end{subequations}

The complete procedure is summarized in \textbf{Algorithm} \ref{alg:so-admm}. This approach significantly reduces the dimensionality of subproblem of ADMM, from size of $W_{1}$ and $W_{2}$ to much smaller size of $B_{1}$ and $B_{2}$.

\begin{algorithm}[t]
    \caption{Subspace Optimization with ADMM (SO-ADMM)}\label{alg:so-admm}
    \textbf{Input:} $\bW_{1}^0$ and $\bW_{2}^0$; outer steps $K$; proximal coefficients $\{\eta^k\}$; maximum inner steps $T_{\max}$.\\
    \textbf{Output:} $\bW_{1}^{K}$ and $\bW_{2}^{K}$.
    \begin{algorithmic}[1]
        \For{$k = 0, 1, \dots, K-1$}
            \State Sample or Compute projection matrix $\bP^k$ and $\bR^k$.
            \State Compute new constant $\tilde{C}=C-A_{1}W_{1}^k-A_{2}W_{2}^k$
            \State Initialize subspace parameters $\bB_{1}^0$, $\bB_{2}^0$ and $\bY^0$.
            \For{$t = 0, 1, \dots, T_{\max}-1$}
                \State $B_{1} \leftarrow \underset{B}{\text{argmin}} \ \mathcal{L}_{\rho}^{k}(B, B_{2}, Y)$.
                \State $B_{2} \leftarrow \underset{B}{\text{argmin}} \ \mathcal{L}_{\rho}^{k}(B_{1}, B, Y)$.
                \State $Y \leftarrow Y + \rho (A_{1}P^{k}B_{1} + A_{2}R^{k}B_{2}-\tilde{C})$.
            \EndFor
            \State Update global weights: $\bW_{1}^{k+1} \leftarrow \bW_{1}^k + \bP^k \bB_{1}$.
            \State Update global weights: $\bW_{2}^{k+1} \leftarrow \bW_{2}^k + \bR^k \bB_{2}$.
        \EndFor
    \end{algorithmic}
\end{algorithm}

\subsubsection{Numerical Example}
To empirically evaluate the performance of the proposed Subspace Optimization version of Alternating Direction Method of Multipliers, we consider a matrix-variable quadratic programming (QP) problem with equality constraints. The problem is formulated as follows:
\begin{equation}
\begin{aligned}
\min_{X,Z} \quad & \frac{1}{2}\mathrm{tr} \left(X^{\top}Q_{1}X\right)+\left\langle D_{1},X \right\rangle + \frac{1}{2}\mathrm{tr} \left( Z^{\top}Q_{2}Z\right) + \left\langle D_{2}, Z\right\rangle, \\
\mathrm{s.t.} \quad & A_{1}X+A_{2}Z-C=0.
\end{aligned}
\end{equation}
where $X \in \mathbb{R}^{m_{1}\times n}$ and $Z \in \mathbb{R}^{m_{2}\times n}$ are the decision variables, and $A_1$ and $A_2$ are assumed to have full row rank. A fundamental characteristic of this specific experimental setup is the inherent low-rank nature of the optimal solution. We achieve the low-rank property of optimal solution by setting $D_{1}$ and $D_{2}$ as low-rank matrices. In our experiments, the horizontal axis represents the inner iterations, while the vertical axis denotes the relative error (distance) of the variables $X$ and $Z$ with respect to their true low-rank solutions.

The purpose of these numerical experiments is twofold: first, to investigate the impact of the subspace size (i.e., the target rank) on the convergence dynamics; and second, to compare the efficacy of various subspace selection strategies, including a comparative analysis against Block Coordinate Descent (BCD) methodologies.

\begin{figure}[t]
    \center
    \includegraphics{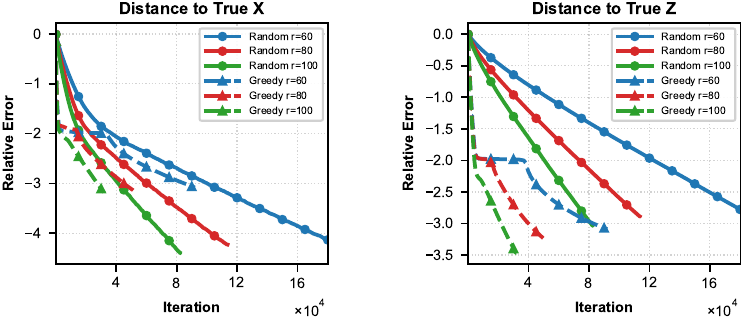}
    \caption{\small Subspace ADMM with different ranks. Left: Relative Error of X. Right: Relative Error of Z. Higher ranks accelerate iteration-wise convergence. Although Greedy outperforms Random, it exhibits a plateau phase—prolonged at lower ranks—due to local over-exploitation.}\label{fig:admm_rank_abla}
\end{figure}

\paragraph{Effect of Subspace Rank}
We first validate the trade-off between per-step computation efficiency (determined by subspace rank r) and convergence speed. We compare Greedy Subspace (SVD) ADMM and Random Subspace ADMM with subspace ranks $r \in \{ 60, 80, 100 \}$ respectively.

Figure \ref{fig:admm_rank_abla} illustrates that whether using the greedy or random subspace selection methods, higher rank choices consistently yield faster iteration-wise convergence. This validates the trade-off between the efficiency of each step's computation and the convergence speed. A higher rank subspace method requires more computational resources per step, but only requires fewer iterations to converge.

Meanwhile, in Figure \ref{fig:admm_rank_abla}, the dashed lines of the same color are always located below the solid lines. This demonstrates that under the same subspace rank, the convergence speed of the greedy subspace method is consistently superior to that of random subspace method.

It is worth noting that all three curves of the greedy subspace with different ranks have shown a period of plateau of relative error in optimization process. This plateau phenomenon stems from over-exploitation of locally dominant residual directions generated by singular value decomposition in greedy subspace method. By prioritizing these directions, the greedy subspace method systematically neglects other structurally critical dimensions. Consequently, a plateau occurs because the algorithm must wait for the dual variables in these unselected dimensions to accumulate sufficient penalty before the greedy criterion finally targets them. A smaller subspace rank inherently restricts the number of dimensions updated per iteration, thereby necessitating a greater number of iterations for sufficient penalty accumulation of dual variables in unexplored directions, which significantly prolongs the plateau phase. In contrast, the random subspace selection strategy uniformly explores the entire parameter space over time. This continuous, unbiased exploration prevents dual accumulation lag, thereby avoiding stagnation and ensuring a steady descent trajectory.

\begin{figure}[t]
    \center
    \includegraphics{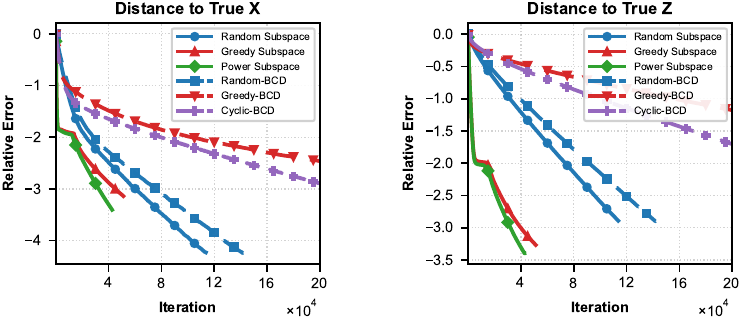}
    \caption{\small Comparison of Subspace Selection Strategy. Left: Relative Error of X. Right: Relative Error of Z. Subspace methods consistently outperform BCD variants. PowerSGD approximates the top-$r$ subspace effectively, even surpassing exact SVD performance slightly. Greedy-BCD underperforms Random/Cyclic-BCD, as row selection mismatches low-rank structure of problem and greedy selection makes it worse.}\label{fig:admm_subspace_abla}
\end{figure}

\paragraph{Effect of Subspace Selection Strategy}
We further investigate the differences between specific subspace selection mechanisms and traditional coordinate-wise updates. 
We compare six ``subspace'' selection strategies defined in \ref{subspace-choice}. Particularly, the gradient of $X$ and $Z$ used in Greedy Subspace, Power Subspace and Greedy-BCD are $\nabla_{X} \mathcal{L}_{\rho}(X, Z, Y)$ and $\nabla_{Z} \mathcal{L}_{\rho}(X, Z, Y)$ respectively.
We set rank $r=80$ as the budget for Subspace methods and BCD methods. Section~\ref{subspace-choice} explains how we match the update budgets. We still use the relative error of $X$ and $Z$ as the comparison metric. Figure \ref{fig:admm_subspace_abla} shows the results.

A prominent conclusion drawn from the empirical results is the consistent superiority of all subspace methods over the BCD variants. The subspace methods achieve significantly lower relative errors at a much faster rate, firmly validating the advantage of exploiting inherent low-rank structures via subspace projections rather than traditional row-wise coordinate updates. Furthermore, the Power Subspace strategy demonstrates a convergence trajectory that is highly comparable to, and even slightly better than, the exact SVD-based Greedy method. This successfully confirms that Power Iteration serves as an excellent and computationally efficient randomized approximation for capturing the top-$r$ dominant eigenspace, matching the descent capability of exact SVD without its prohibitive computational overhead.

Interestingly, the Greedy-BCD strategy—employing the Gauss-Southwell rule for row selection—performs noticeably worse than both Random-BCD and Cyclic-BCD. This behavior most likely stems from a critical combination of structural mismatch and algorithmic over-exploitation. Fundamentally, strict row-wise updates fail to capture the underlying low-rank structure of the solution, which inherently spans across multiple coordinates rather than isolated rows. While Random-BCD and Cyclic-BCD partially mitigate this limitation by enforcing uniform, periodic exploration of all dimensions, Greedy-BCD severely over-exploits a small subset of ``active'' rows with temporarily large gradients. Under the coupled equality constraints of ADMM, this greedy myopia traps the algorithm in local sub-manifolds.

\subsection{Primal-Dual Hybrid Gradient Algorithm (PDHG)}
\subsubsection{Preliminary}
The Primal-Dual Hybrid Gradient (PDHG) method is a highly efficient first-order algorithmic framework tailored for solving large-scale, non-smooth convex optimization problems. Chambolle and Pock established an $O(1/k)$ ergodic primal-dual gap bound for their extrapolated primal-dual algorithm under suitable step-size conditions \cite{chambolle2011first}.

Fundamentally, PDHG is designed to tackle saddle-point (min-max) problems, which have a direct connection to constrained optimization. Consider a general constrained optimization problem of the form:
\begin{equation}
\begin{aligned}
\underset{x, y}{\text{minimize}} \quad & g(x) + f(y), \\
\text{subject to} \quad & Ax = y.
\end{aligned}
\end{equation}
By introducing a dual variable  $z$ to handle the linear constraint and utilizing the Fenchel convex conjugate $f^*$, this constrained problem can be transformed into the unconstrained primal-dual saddle-point formulation:
\begin{equation} \label{eq:pdhg_saddle}
\min_{x} \max_{z} \ \langle Ax, z \rangle + g(x) - f^*(z).
\end{equation}
where $g$ and $f^*$ are proper, lower semi-continuous convex functions, and $A$ is a continuous linear operator.

To solve this minimax problem, the standard PDHG algorithm alternates in the dual and primal spaces. The iterative updates at step k are formulated as follows:
\begin{subequations}
\begin{align}
z^{k+1} &= \underset{z}{\text{argmax}}\{-f^*(z) + \langle Ax^{k}, z \rangle - \frac{1}{2\delta_{k}}\|z-z^{k}\|_{2}^{2}\}, \\
x^{k+1} &= \underset{x}{\text{argmin}}\{g(x) + \langle A^{\top}z^{k+1}, x\rangle + \frac{1}{2\alpha_{k}}\| x-x^{k} \|_{2}^{2}\}.
\end{align}
\end{subequations}

\subsubsection{Update Derivation}
To apply \emph{Matrix Subspace Optimization} in PDHG's updating rule, we need to transform the full-dimensional saddle-point problem into the form of subspace variables. The following problem is matrix version of \eqref{eq:pdhg_saddle}:

\begin{equation} \label{eq:pdhg_saddle_problem}
\min_{W_{1}} \max_{W_{2}} \ \langle AW_{1}, W_{2} \rangle + g(W_{1}) - f^*(W_{2}).
\end{equation}

where vector variable $x$ and $z$ are substituted by matrix variable $W_{1}$ and $W_{2}$ respectively.

Now, we substitute the full-dimension variable $W_{1}$ and $W_{2}$ by $W_{1}^{k}+P^{k}B_{1}$ and $W_{2}^{k}+R^{k}B_{2}$ respectively at outer iteration k, where $P^k$ and $R^k$ are sampled subspace and $W_{1}^{k}$ and $W_{2}^{k}$ are fixed full-dimension variable. Besides, we will add a regularization term for subspace variable $B_1$ and $B_2$. Then we get:

\begin{equation} \label{eq:so_pdhg_saddle_problem}
\min_{B_{1}} \max_{B_{2}} \ \langle A(W_{1}^{k}+P^{k}B_{1}), W_{2}^{k}+R^{k}B_{2} \rangle + \left\{g(W_{1}^{k}+P^{k}B_{1}) + \frac{1}{2\eta^{k}} \|B_1\|_{F}^{2} \right\} - \left\{f^*(W_{2}^{k}+R^{k}B_{2}) + \frac{1}{2\eta^{k}}\|B_2\|_{F}^{2}\right\}.
\end{equation}

Now, we can follow the classical update rule of PDHG to solve subproblem \eqref{eq:so_pdhg_saddle_problem}:

\begin{subequations}
\begin{align}
B_{2}^{t+1} &= \underset{B_{2}}{\text{argmax}}\left\{- \left[f^*(W_{2}^{k}+R^{k}B_{2}) + \frac{1}{2\eta^{k}}\|B_2\|_{F}^{2}\right] + \left< A(W_{1}^{k}+P^{k}B_{1}^{t}), W_{2}^{k}+R^{k}B_{2} \right> - \frac{1}{2\delta_{k}}\|B_{2}-B_{2}^{t}\|_{F}^{2}\right\}, \\
B_{1}^{t+1} &= \underset{B_{1}}{\text{argmin}}\left\{\left[g(W_{1}^{k}+P^{k}B_{1}) + \frac{1}{2\eta^{k}} \|B_1\|_{F}^{2} \right]  + \left< A^{\top}(W_{2}^{k}+R^{k}B_{2}^{t+1}), W_{1}^{k}+P^{k}B\right> + \frac{1}{2\alpha_{k}}\| B_{1}-B_{1}^{t} \|_{F}^{2}\right\}.
\end{align}
\end{subequations}

To simplify the form of updating rule, we introduce the following notations:

\begin{align}
\Phi_{k}^{(t)}(B;B_{1}^{t},B_{2}^{t})&:=- \left[f^*(W_{2}^{k}+R^{k}B) + \frac{1}{2\eta^{k}}\|B\|_{F}^{2}\right] + \left< A(W_{1}^{k}+P^{k}B_{1}^{t}), W_{2}^{k}+R^{k}B \right> - \frac{1}{2\delta_{k}}\|B-B_{2}^{t}\|_{F}^{2},\\
\Psi_{k}^{(t)}(B;B_{1}^{t},B_{2}^{t+1})&:=\left[g(W_{1}^{k}+P^{k}B) + \frac{1}{2\eta^{k}} \|B\|_{F}^{2} \right]  + \left< A^{\top}(W_{2}^{k}+R^{k}B_{2}^{t+1}), W_{1}^{k}+P^{k}B\right> + \frac{1}{2\alpha_{k}}\| B-B_{1}^{t} \|_{F}^{2}.
\end{align}

Then the complete procedure is summarized in \textbf{Algorithm} \ref{alg:so-pdhg}. This approach significantly reduces the dimensionality of subproblem of PDHG, from size of $W_{1}$ and $W_{2}$ to much smaller size of $B_{1}$ and $B_{2}$.

\begin{algorithm}[t]
    \caption{Subspace Optimization with PDHG (SO-PDHG)}\label{alg:so-pdhg}
    \textbf{Input:} $\bW_{1}^0$ and $\bW_{2}^0$; outer steps $K$; proximal coefficients $\{\eta^k\}$; maximum inner steps $T_{\max}$; regularization coefficients $\{\delta_{k}\}$ and $\{\alpha_{k}\}$.\\
    \textbf{Output:} $\bW_{1}^{K}$ and $\bW_{2}^{K}$.
    \begin{algorithmic}[1]
        \For{$k = 0, 1, \dots, K-1$}
            \State Sample or Compute projection matrix $\bP^k$ and $\bR^k$.
            \State Initialize subspace parameters $\bB_{1}^0$ and $\bB_{2}^0$.
            \For{$t = 0, 1, \dots, T_{\max}-1$}
                \State $B_{2} \leftarrow \underset{B}{\text{argmax}}  
                \Phi_{k}^{(t)}(B;B_{1},B_{2})$.
                \State $B_{1} \leftarrow \underset{B}{\text{argmin}} \ \Psi_{k}^{(t)}(B;B_{1},B_{2})$.
            \EndFor
            \State Update global weights: $\bW_{1}^{k+1} \leftarrow \bW_{1}^k + \bP^k \bB_{1}$.
            \State Update global weights: $\bW_{2}^{k+1} \leftarrow \bW_{2}^k + \bR^k \bB_{2}$.
        \EndFor
    \end{algorithmic}
\end{algorithm}

\subsubsection{Numerical Example}

To evaluate the proposed Subspace Primal-Dual Hybrid Gradient (Subspace PDHG) method, we consider a class of saddle-point problems featuring an inherent low-rank structure. The problem is derived from a quadratic programming formulation with a proximal term, expressed as the unconstrained minimax problem:
\begin{equation}
\min_{X} \max_{Z} \ f(X) - h^*(Z) + \mathrm{tr}(Z^\top A X).
\end{equation}
where $f(X) = \frac{1}{2}\mathrm{tr}(X^\top Q X) + \langle D, X \rangle$ and $h(Z) = \frac{1}{2}\|Z-B\|_F^2$. After using the definition of conjugate function, we can get $h^*(Z) = \frac{1}{2}\|Z\|_F^2 + \langle Z, B \rangle$. For simplicity, we denote $\mathcal{L}(X,Z) = f(X) - h^*(Z) + \mathrm{tr}(Z^\top A X)$.

A fundamental characteristic of this experimental setup is that the optimal solutions $X$ and $Z$ inherently reside on a low-rank manifold. We achieve this property of low-rank by setting matrix $B$ and $D$ as low-rank matrices. In our experiments, the horizontal axis represents the inner iterations, while the vertical axis denotes the relative error (distance) of the variables $X$ and $Z$ with respect to their true underlying low-rank solutions. The primary objectives of these numerical experiments are to investigate the ablation of subspace rank sizes and to systematically compare subspace selection strategies against BCD methods.

\begin{figure}[t]
    \center
    \includegraphics{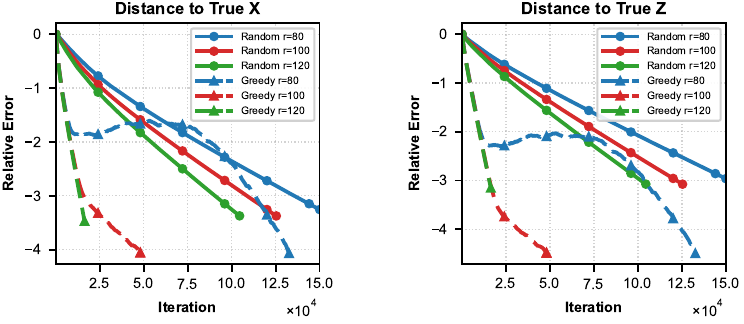}
    \caption{\small Subspace PDHG with different ranks. Left: Relative Error of X. Right: Relative Error of Z. Higher ranks accelerate convergence. Greedy outperforms Random but exhibits a transient relative error increase at smaller ranks due to subspace drift.}\label{fig:pdhg_rank_abla}
\end{figure}

\paragraph{Effect of Subspace Rank}
We first validate the trade-off between per-step computation efficiency (determined by the subspace rank $r$) and convergence speed. We compare the Greedy Subspace (SVD) PDHG and the Random Subspace PDHG with subspace ranks $r \in \{80, 100, 120\}$ respectively.

Figure \ref{fig:pdhg_rank_abla} illustrates that, whether utilizing greedy or random subspace selection methods, higher rank choices consistently yield faster iteration-wise convergence. A higher rank subspace method incorporates a larger volume of gradient information per step, thus requiring fewer iterations to converge. Concurrently, under the same subspace rank, the convergence iterations of the greedy subspace method are consistently fewer than those of the random subspace method, highlighting the efficiency of principal gradient directions.

It is particularly noteworthy that for the greedy subspace method with a smaller rank (e.g., $r=80$), the relative error exhibits a transient period of increase during the optimization process. During this phase, the algorithmic primal and dual residuals are actually still monotonically decreasing. This counterintuitive relative error inflation occurs because a small-rank greedy strategy intensely over-exploits a narrowly selected subspace based on the current local gradient. This causes the variables to temporarily drift away from the global true low-rank manifold. Once the dual variables accumulate sufficient information to correct this trajectory, the relative error resumes its rapid descent.

\begin{figure}[t]
    \center
    \includegraphics{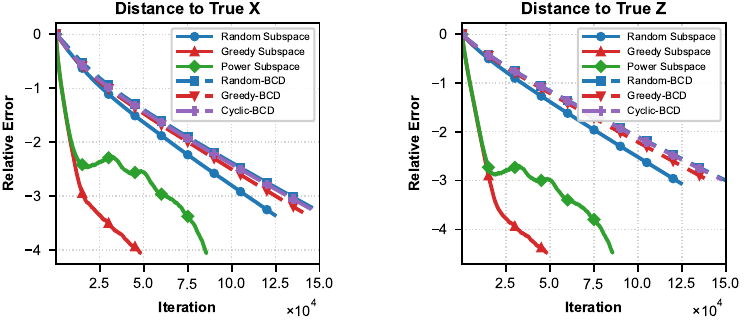}
    \caption{\small Comparison of Subspace Selection Strategy. Left: Relative Error of X. Right: Relative Error of Z. Subspace methods consistently outperform BCD variants, capturing the low-rank structure effectively. PowerSGD perfectly approximates the SVD trajectory initially but jitters later. Greedy-BCD resembles other BCDs, as strict row-wise updates cannot capture global low-rank geometry.}\label{fig:pdhg_subspace_abla}
\end{figure}

\paragraph{Effect of Subspace Selection Strategy}
We further investigate the differences between specific subspace selection mechanisms and traditional coordinate-wise updates. 
We compare six ``subspace'' selection strategies as defined in \ref{subspace-choice}. Particularly, the gradient of $X$ and $Z$ used in Greedy Subspace, Power Subspace and Greedy-BCD are $\nabla_{X} \mathcal{L}_{\rho}(X, Z)$ and $\nabla_{Z} \mathcal{L}_{\rho}(X, Z)$ respectively.
We set rank $r=100$ as the budget for Subspace methods and BCD methods. Section~\ref{subspace-choice} explains how we match the update budgets. We still use the relative error of variables $X$ and $Z$ as the comparison metric. Figure \ref{fig:pdhg_subspace_abla} shows the results.

A prominent conclusion drawn from the empirical results in Figure \ref{fig:pdhg_subspace_abla} is the consistent superiority of all subspace methods over the BCD variants. The subspace methods achieve significantly lower relative errors at a much faster rate, firmly validating the advantage of exploiting inherent low-rank structures via subspace projections. Furthermore, the convergence curve of the PowerSGD strategy nearly perfectly overlaps with the exact SVD (Greedy Subspace) curve in the early stages. This confirms that PowerSGD serves as an excellent randomized approximation for capturing the top-$r$ dominant eigenspace.

However, two anomalous behaviors warrant detailed analysis. First, the PowerSGD curve exhibits noticeable jittering in the later stages of optimization. This is probably because, as the algorithm approaches convergence, the gradients become vanishingly small, and the iterative power method struggles with numerical precision, leading to noisy approximations of the singular vectors. Second, the Greedy-BCD strategy performs almost identically poorly to the Random-BCD and Cyclic-BCD methods. Just like the relevant phenomenon in the ADMM case, this most likely stems from a critical structural mismatch: strict row-wise blocking fundamentally fails to capture the underlying low-rank structure of the solution, which inherently spans across multiple coordinates. Consequently, employing the Gauss-Southwell rule to pick the ``best'' rows offers no meaningful advantage when the very premise of isolated row-updates ignores the global low-rank geometry.

\section{Federated Learning}

Federated Learning (FL) \cite{mcmahan2017communication} has emerged as a transformative paradigm for privacy-preserving distributed optimization, enabling multiple clients to collaboratively train a shared model without exposing their raw data. Consider a federated system with $N$ clients. Let $\bW \in \mathbb{R}^{m \times n}$ denote the global model parameters, and let $F_i(\bW) = \mathbb{E}_{\xi \sim \mathcal{D}_i}[\ell(\bW;\xi)]$ represent the local empirical loss of client $i$ over its private data distribution $\mathcal{D}_i$. The canonical global objective of FL is to minimize the weighted average of these local losses:
\begin{equation}
\label{eq:fl_global_objective}
\min_{\bW \in \mathbb{R}^{m \times n}} f(\bW) := \sum_{i=1}^{N} p_i F_i(\bW),
\end{equation}
where $p_i \ge 0$, $\sum_i^N p_i = 1$ (typically $p_i = |\mathcal{D}_i| / \sum_j |\mathcal{D}_j|$). This formulation addresses the fundamental tension between collaborative learning and data sovereignty.

Despite its conceptual appeal, the practical deployment of FL faces two interconnected challenges that have shaped the evolution of FL algorithms. First, the Communication Bottleneck: synchronizing millions of parameters across bandwidth-limited networks remains prohibitively expensive. Second, and more fundamentally, Statistical Heterogeneity: the ubiquitous non-Independent and Identically Distributed (non-IID) nature of client data causes local updates to drift from the global optimum—a phenomenon known as client drift—leading to slower convergence and degraded global model performance.

These challenges have spurred a rich line of algorithmic innovations. Early efforts focused on stabilizing local training through proximal regularization (e.g., FedProx \cite{li2020federated}) or variance reduction via control variates (e.g., SCAFFOLD \cite{karimireddy2020scaffold}). More recently, the field has witnessed a diversification of approaches tailored to specific bottlenecks. For instance, gradient compression techniques—such as quantization (FedPAQ \cite{reisizadeh2020fedpaq}) and sparsification (DGC \cite{lin2017deep})—have been extensively explored to reduce communication overhead, often combined with error feedback to maintain convergence. Another major direction tackles statistical heterogeneity through personalized federated learning, where methods like pFedMe \cite{t2020personalized} and FedPer \cite{arivazhagan2019federated} learn locally adapted models or via multi-task learning frameworks. Adaptive optimization strategies (e.g., FedAdam, FedYogi \cite{reddi2020adaptive}) have also been proposed to improve convergence under non-IID data, while decentralized architectures (e.g., D-PSGD \cite{lian2017can}) remove the central server to enhance scalability and robustness. 

Notably, a common thread across these advances is the recognition that communication efficiency and robustness to heterogeneity are deeply intertwined. Gradient compression techniques, while reducing communication overhead, often amplify the bias introduced by non-IID data; conversely, methods that excel at handling heterogeneity (e.g., SCAFFOLD) incur additional communication or memory costs. This observation motivates our work: we posit that a parameter-efficient optimization framework, which inherently constrains the optimization trajectory to a low-dimensional subspace, can simultaneously alleviate the communication burden and provide an implicit regularization effect that stabilizes training under heterogeneity.

To this end, this section extends the \emph{Matrix Subspace Optimization} framework to the FL setting. By restricting both local updates and global aggregation to a randomly generated, low-dimensional subspace, Matrix Subspace Optimization offers a principled approach to compressing the communication volume while implicitly regularizing the solution space. This dual benefit positions it as a promising candidate to address the twin challenges of communication efficiency and client drift in a unified manner.

\subsection{Federated Averaging (FedAvg)}

\subsubsection{Preliminary}

\textbf{Federated Averaging (FedAvg)} \cite{mcmahan2017communication} is one of the most fundamental methods in the field of federated learning. In the FedAvg algorithm, each client maintains a local model, while the server maintains a centralized global model. In each iteration, the clients upload their model parameters to the server after local training, and the server then performs weighted average aggregation to complete the update of the global model.

We denote the global model as $W$, and the model of the $i$-th client as $W_i$. Hence, the global model and client model at the iteration $k$ are denoted as $W^k$ and $W_i^k$, respectively. $W_i$ is derived by optimizing the objective function $F_i(W)$ of the local loss. At the $k$-th iteration, after receiving $W^{k}$ from the server, each client performs $T$ steps of SGD locally:

\begin{equation}
    \bW_i^{k,t+1} = \bW_i^{k,t} - \eta \bG_i^{k,t}, \quad t = 0, 1, \dots, T-1.
\end{equation}

Here $\bW_i^{k,0} = \bW^k$, $\eta$ is the learning rate and $\bG_i^{k,t}$ denotes the gradient $\nabla \bF_i(\bW_i^{k,t})$ calculated using a random batch from $\mathcal{D}_i$. Then the resulting
models $\bW_i^k = \bW_i^{k,T}$ are updated and averaged to obtain $W^{k+1}$:

\begin{equation}
    \bW^{k+1} = \sum_{i} p_i \bW_i^k.
\end{equation}

By repeating the process of local client training and server aggregation, we obtain the final global model $W$ of the FedAvg algorithm. In each iteration, the server distributes the global model to the \(N\) clients, and clients upload their distinct models after local training. The total number of transmitted elements is \(2Nmn\).

\subsubsection{Update Derivation}

To mitigate the communication and computation cost, we introduce Matrix Subspace Optimization and constrain the updates to a low-dimensional subspace defined by projection matrix $\bP^k \in \mathbb{R}^{m \times r}$. Instead of updating the model $\bW_i$, client $i$ optimizes a small coefficient matrix $\bB_i \in \mathbb{R}^{r \times n}$. The local objective for SO-FedAvg is formulated as:
\begin{equation}
    \bB_i^k \approx \argmin_{\bB} \biggl\{ F_i(\bW^k + \bP^k \bB) + \frac{1}{2\eta^k} \|\bB\|_F^2 \biggl\},
\end{equation}

We further rewrite the global function of $k$-th iteration as $f(\bW^k+\bP^k \bB)$:
\begin{equation}
    \label{eq:fl_so_objective}
    f(\bW^k+\bP^k \bB) = \sum_{i} p_i \biggl\{ F_i(\bW^k + \bP^k \bB) + \frac{1}{2\eta^k} \|\bB\|_F^2 \biggl\},
\end{equation}
where the term $\frac{1}{2 \eta^k} \|B\|^2_F$ is a proximal term that regularizes the subproblem. Thus, we can project the FL problem onto the low-dimensional subspace, with the specific procedure as follows.

At outer round $k$, server broadcasts subspace projection matrix $\bP^k \in \mathbb{R}^{m \times r}$ (with $r \ll m$) and the current global model $\bW^k$ to all participating clients. Each client $i$ then initializes its subspace parameter $\bB_i^{k,0} = \mathbf{0}$ and performs a $T$ step local training. For $t = 0,1, \dots, T-1$:
\begin{align}
        \bG_i^{k,t} &= (\bP^k)^\top \nabla F_i(\bW^k + \bP^k \bB_i^{k,t}) + \frac{1}{\eta^k} \bB_i^{k,t},\\
        \bB_i^{k,t+1} &= \bB_i^{k,t} - \alpha_t \bG_i^{k,t}.
\end{align}

Client $i$ sends $\bB_i^k = \bB_i^{k,T}$ back to the server, and server aggregates the received subspace updates via weighted averaging, then updates the global model. The detailed procedure is summarized in \ref{alg:so-fedavg}.

\begin{algorithm}[t]
    \caption{Subspace Optimization for Federated Averaging (SO-FedAvg)}
    \label{alg:so-fedavg}
        \textbf{Input:} $\bW^0$, outer steps $K$, proximal coefficients $\{\eta^k\}$, maximum inner steps $T_{\max}$, inner step sizes $\{\alpha_t\}$, client sampling probabilities $\{p_i\}$.\\
        \textbf{Output:} $\bW^K$
        \begin{algorithmic}[1]
        \For{$k = 0, 1, \ldots, K-1$}
            \State Server samples or computes a subspace matrix $\bP^k$.
            \State Server broadcasts $\bW^k$ and $\bP^k$ to selected clients.
            \For{each selected client $i$ \textbf{in parallel}}
                \State Initialize $\bB_i \gets \mathbf{0}$.
                \For{$t = 0, 1, \ldots, T_{\max}-1$}
                    \State Compute low-dimensional gradient:
                        $\bG_i \gets (\bP^k)^\top \nabla F_i\!\left( \bW^k + \bP^k \bB_i \right)+ \frac{1}{\eta^k} \bB_i$.
                    \State Update low-rank parameters: $\bB_i \gets \bB_i - \alpha_t \bG_i$.
                    \State \textbf{if} {$B_i$ converges} \textbf{then} $B_i^{k} \gets B_i$ and \textbf{break}
                \EndFor
                \State Client sends $\bB_i^k$ to server.
            \EndFor
            \State Server aggregates $\bB^k \gets \sum_{i} p_i \bB_i^k$.
            \State Server updates global model: $\bW^{k+1} \gets \bW^k + \bP^k \bB^k$.
        \EndFor
    \end{algorithmic}
\end{algorithm}

For the SO-FedAvg algorithm, in each iteration, the server distributes the initial global model and the projection matrix to the \(N\) clients. The \(N\) clients upload the low-dimensional parameter \(B\) to the server, and the total number of elements transmitted is \(N(mn+mr+rn)\). When \(r \ll m\), we may assume that \(m\) and \(n\) are of similar magnitude, then \(N(mn+mr+rn) \ll 2Nmn\). It can be seen that SO-FedAvg significantly reduces the communication overhead compared with FedAvg.

\subsubsection{Numerical Example}

To verify the effectiveness of algorithms of SO-FedAvg, we conduct experiments on classification tasks using the MNIST dataset. We adopt the classification model as defined in \ref{MLP-setting}. The experimental objective is chosen as follows: for a FL system of $N=10$ clients, SO-FedAvg is used to handle the scenario where client data is independently and identically distributed (i.i.d.). We compare the training process under different subspace ranks and different subspace selection strategies, and plot the curves of the average loss of client models versus training steps, as well as the test accuracy of the server model versus training steps.

\paragraph{Effect of Subspace Rank}
First, we investigate the impact of the subspace rank $r$ on the FedAvg training process. Specifically, we test the cases where the rank takes values $r \in \{16, 32, 64\}$, and the projection subspaces are selected via Random Subspace and Greedy Subspace (implemented by the SVD method) strategies, while the full-rank case for this problem is $r=128$. The outer loop of FedAvg consists of 15 rounds, with each client performing 100 local training steps per round. Since the model of the server does not inherently have a gradient, the greedy subspace method requires each client to additionally upload their gradients to the server. Server then computes the average gradient and extracts its top $r$ singular vectors to construct the subspace.

\begin{figure}[t]
    \center
    \includegraphics{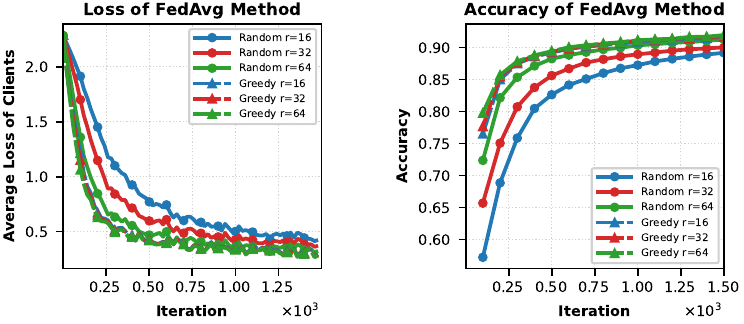} 
    \caption{\small Comparison plots of the training performance of the Subspace+FedAvg algorithm with randomly selected subspaces and greedily selected subspaces under different subspace ranks. The plots show the trends of the average client loss and server test accuracy as functions of the total number of iterations.}\label{fig:so-fedavg-rank}
\end{figure}

As shown in figure~\ref{fig:so-fedavg-rank}, SO-FedAvg exhibits strong robustness to the reduction of rank. In the Loss curve subfigure, the dashed lines of the same color consistently lie below the solid lines, whereas the opposite holds in the Accuracy curve subfigure. This indicates that given the same subspace rank, the convergence speed of the greedy subspace method is consistently superior to that of the random subspace method.

Meanwhile, both subfigures show that, whether using the greedy subspace method or the random subspace method, higher rank choices consistently yield faster iteration-wise convergence. This validates the trade-off between per-step computational efficiency and convergence speed: a higher-rank subspace method requires more computational resources per step but needs fewer iterations to converge.

\paragraph{Effect of Subspace Selection Strategy}
Next, we fix the subspace rank at $r=16$ and conduct a horizontal comparison of the impact of different projection subspace selection strategies on training. Specifically, we adopt 6 projection subspace selection strategies, with the same 15 outer FedAvg loops and 100 local training steps per client per loop. The 6 strategies of selecting the subspace are roughly the same as those in \ref{subspace-choice}. The difference is that, since the model at the server does not have a gradient, we need to use the mean of the gradients of the local models from all clients as a replacement to carry out the computation of Greedy Subspace, Power Subspace and Greedy-BCD.

\begin{figure}[t]
    \center
    \includegraphics{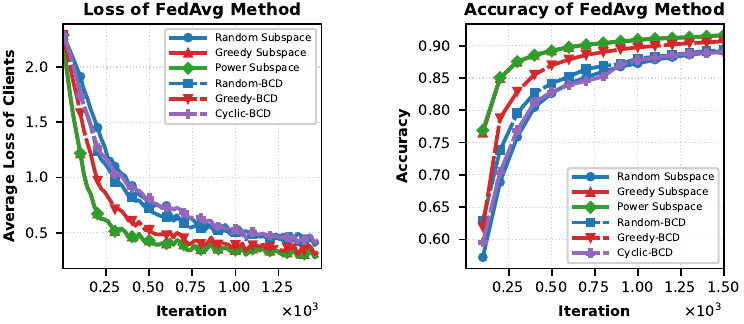} 
    \caption{\small Comparison plots of the training performance of the Subspace+FedAvg algorithm with different subspace selection strategies. The plots show the trends of the average client loss and server test accuracy as functions of the total number of iterations.}\label{fig:so-fedavg-method}
\end{figure}

As can be seen from figure~\ref{fig:so-fedavg-method}, the Greedy Subspace and Power Subspace methods achieve the fastest iteration-wise convergence, while the BCD-type methods and the Random Subspace method converge more slowly. This indicates that the power iteration strategy in this experiment effectively approximates the subspace spanned by the top-$r$ singular vectors of the average gradient, significantly reducing the computational overhead at the server. Greedy-BCD notably outperforms the other two BCD methods and random subspace method; however, Greedy-BCD also requires each client to additionally upload gradients, which implies the increase in communication overhead.

Based on the above analysis, it is evident that SO-FedAvg can leverage the inherent low-rank structure in FL, demonstrating clear advantages when comprehensively considering communication overhead, computational cost, and convergence speed, thus offering significant potential for research and application.

\subsection{Federated Learning with Gradient Tracking (FLGT)}

\subsubsection{Preliminary}

In real-world scenarios, the data distribution among clients is usually non-i.i.d., and the properties of clients vary across different participants, which results in the phenomenon of client drift. Studies have shown that FedAvg struggles with heterogeneity. To alleviate this problem, we introduce the \textbf{gradient tracking (GT)} idea from distributed optimization \cite{nedic2017achieving} and propose the algorithm named \textbf{Federated Learning Gradient Tracking (FLGT)}.

We first introduce the original \textbf{Gradient Tracking (GT)} algorithm. In the following distributed optimization problem, \( \mathbf{x} \) is an \( m \times n \) matrix, whose \( i \)-th row represents the optimization variable locally maintained by agent \( i \). The update rule for \( \mathbf{x} \) is given by:
\begin{equation}
\mathbf{x}(k+1)=\vW\mathbf{x}(k) - \alpha \mathbf{y}(k),
\end{equation}
where \( \mathbf{W} \) is the mixing matrix for communication and averaging, and \( \mathbf{y} \) denotes the gradient tracking variable. The update formula for \( \mathbf{y} \) is as follows.
\begin{equation}
\mathbf{y}(k+1)=\vW\mathbf{y}(k)+ \nabla \vF (\mathbf{x}(k+1))-\nabla \vF (\mathbf{x}(k)).
\end{equation}

By maintaining the auxiliary variable $y$ that estimates the global gradient, GT allows each agent to align its local updates with the overall descent direction, effectively reducing client drift. This leads to faster convergence and enhanced robustness to non-IID data distributions and time-varying network topologies.

We now extend the concept of GT to the FL setting, resulting in the algorithm of FLGT. FLGT differs from FedAvg in both local training and communication procedures. In each iteration, FLGT consists of the following steps: clients maintain gradient tracking variables locally and perform local training; decentralized communication is performed periodically among clients to average model parameters and GT variables; the server aggregates the uploaded model parameters from clients and update the global model. At the beginning of round $k$, client $i$ initializes the local tracking variable $Y_i \in \mathbb{R}^{r \times n}$ to $\mathbf{0}$ and receives the current global model $W^k$. After that the local model $W_i$ is updated according to the following rules:
\begin{align}
    \bW_i^{k,t+1} &= \bW_i^{k,t} - \alpha_t \bY_i^{k,t},  \\
    \bY_i^{k,t+1} &= \bY_i^{k,t} + \left( \nabla F_i\!\left( \bW_i^{k,t+1} \right) - \nabla F_i\!\left( \bW_i^{k,t} \right) \right).
\end{align}

Here $t = 0,1, \dots, T-1$. In the decentralized optimization setting, each update of GT method performs communication and aggregation via the matrix $\mathbf{W}$. However, in FL we conduct such aggregation only once every $\tau$ steps. The model parameter and GT variables are averaged: $\bar{\bW} = \sum_{j} p_j \bW_j$, $\bar{\bY} = \sum_{j} p_j \bY_j$. We then set $\bW_i = \bar{\bW}$, $\bY_i = \bar{\bY}$ for all participating clients. After $T$ steps of training, the server aggregates client models and updates the global model. The algorithm then proceeds to the next iteration.

\subsubsection{Update Derivation}

Combining Matrix Subspace Optimization with FLGT yields SO-FLGT. This algorithm is particularly powerful as it reduces communication while explicitly addressing data heterogeneity through the gradient tracking variable in the subspace. While the notations and problem setting are similar to subsection of SO-FedAvg above, here we present the procedure as follows. In order to solve the local objective $\min_{\bB_i} F_i(\bW^k + \bP^k \bB_i) + \frac{1}{2\eta^k} \|\bB_i\|_F^2$, each client $i$ maintains a low-dimensional model parameter $\bB_i \in \mathbb{R}^{r \times n}$ and a gradient tracking variable $Y_i \in \mathbb{R}^{r \times n}$.

At outer round $k$, server broadcasts subspace projection matrix $\bP^k \in \mathbb{R}^{m \times r}$ (with $r \ll m$) and the current global model $\bW^k$ to all participating clients. Each client $i$ then initializes $\bB_i^0 = \mathbf{0}$, $\bY_i^0 = \mathbf{0}$ and performs $T$ steps of local training. For $t = 0, \dots, T-1$:
\begin{align}
    \bB_i^{k,t+1} &= \bB_i^{k,t} - \alpha_t \bY_i^{k,t},  \\
    \bY_i^{k,t+1} &= \bY_i^{k,t} + \left( \bG_i(\bB_i^{k,t+1}) -  \bG_i(\bB_i^{k,t})\right).
\end{align}
where $G_i$ denotes the gradient of the objective function for client $i$: $\bG_i(\bB) = (\bP^k)^\top \nabla F_i\!\left( \bW^k + \bP^k \bB \right) + \frac{1}{\eta^k} \bB$. When $(t+1) \bmod \tau = 0$, clients communicate and average the models and gradient tracking variables.
\begin{align}
    \bar{\bB} = \sum_{j} p_j \bB_j,  \\
    \bar{\bY} = \sum_{j} p_j \bY_j.
\end{align}
Then for each client $i$, $\bB_i \gets \bar{\bB}$, $\bY_i \gets \bar{\bY}$. When the training process is completed, client $i$ sends the accumulated subspace update $\bB_i^k \gets \bB_i^{k,T-1}$ back to the server. Server aggregates the received subspace updates via weighted averaging, then updates the global model. The complete procedure is presented in algorithm \ref{alg:so-flgt}. Similar to the analysis in the SO-FedAvg section, the SO-FLGT algorithm can also significantly reduce the communication overhead compared with the original FLGT algorithm.

\begin{algorithm}[t]
    \caption{Subspace Optimization with Gradient Tracking for Federated Learning (SO-FLGT)}
    \label{alg:so-flgt}
    \textbf{Input:} $\bW^0$, outer steps $K$, proximal coefficients $\{\eta^k\}$, maximum inner steps $T_{\max}$, inner step sizes $\{\alpha_t\}$, client sampling probabilities $\{p_i\}$, GT aggregation cycle $\tau$.\\
    \textbf{Output:} $\bW^K$
    \begin{algorithmic}[1]
        \For{$k = 0, 1, \ldots, K-1$}
            \State Server samples or computes a subspace matrix $\bP^k$.
            \State Server broadcasts $\bW^k$ and $\bP^k$ to clients.
            \For{each selected client $i$ \textbf{in parallel}}
                \State Initialize $\bB_i \gets \mathbf{0}$, $\bY_i \gets \mathbf{0}$, $\bG_i^{\text{old}} \gets \mathbf{0}$.
                \For{$t = 0, 1, \ldots, T_{\max}-1$}
                    \State Update low-rank parameters: $\bB_i \gets \bB_i - \alpha_t \bY_i$.
                    \State Compute low-dimensional gradient:
                        $\bG_i \gets (\bP^k)^\top \nabla F_i\!\left( \bW^k + \bP^k \bB_i \right) + \frac{1}{\eta^k} \bB_i$.
                    \State Update gradient tracker: $\bY_i \gets \bY_i + \left( \bG_i - \bG_i^{\text{old}} \right)$.
                    \State Store old gradient: $\bG_i^{\text{old}} \gets \bG_i$.
                    \If{$(t+1) \bmod \tau = 0$}
                        \State Aggregates $\bB \gets \sum_{i} p_i \bB_i$, $\bY \gets \sum_{i} p_i \bY_i$.
                        \State Broadcasts $\bB$ and $\bY$ to all clients: $\bB_i \gets \bB$ , $\bY_i \gets \bY$.
                    \EndIf
                    \State \textbf{if} {$B_i$ converges} \textbf{then} $B_i^{k} \gets B_i$ and \textbf{break}
                \EndFor
                \State Client sends $\bB_i^k$ to server.
            \EndFor
            \State Server aggregates $\bB^k \gets \sum_{i} p_i  \bB_i^k$.
            \State Server updates global model: $\bW^{k+1} \gets \bW^k + \bP^k \bB^k$.
        \EndFor
    \end{algorithmic}
\end{algorithm}

\subsubsection{Numerical Example}

To verify the effectiveness of algorithms of SO-FLGT, we conduct experiments approximately the same as experiments for SO-FedAvg. Nevertheless, the client data distribution is non-i.i.d. during the training. We compare the training process under different subspace ranks and different subspace selection strategies.

\paragraph{Effect of Subspace Rank}
First, we investigate the impact of the subspace rank $r$ on the FLGT training process. The main experimental parameters and subspace selection strategies are consistent with those in the corresponding experiments of SO-FedAvg. However, due to the presence of gradient tracking variables in the GT algorithm, each client performs mixing of the optimization variable B and the gradient tracking variable Y every 20 local training steps. That is, each client conducts 100 local training steps between two updates of the projection subspace by the server, but performs 20 local training steps between two averaging of parameters $\bB$ and $\bY$. We test the accuracy of the server-side model after each averaging.

\begin{figure}[t]
    \center
    \includegraphics{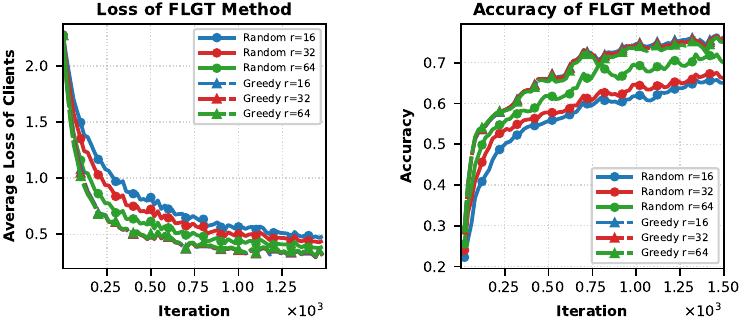} 
    \caption{\small Comparison plots of the training performance of the Subspace+FLGT algorithm with randomly selected subspaces and greedily selected subspaces under different subspace ranks. The plots show the trends of the average client loss and server test accuracy as functions of the total number of iterations.}\label{fig:so-flgt-rank}
\end{figure}

As shown in figure~\ref{fig:so-flgt-rank}, SO-FLGT also exhibits robustness to the reduction of rank. Given the same subspace rank, the convergence speed of the greedy subspace method is consistently superior to that of the random subspace method, but compared to random subspace method, applying the greedy subspace approach in FLGT leads to increase in communication and computation overhead.

Both subfigures show that for the random subspace method, choosing a larger rank yields faster iteration-wise convergence. Nevertheless, for the greedy subspace method, the convergence speeds under different ranks are almost identical. This not only validates that higher-rank subspace methods, despite requiring more computational resources per step, need fewer iterations to converge, but also indicates that there exists an upper bound on the iteration-wise convergence speed of subspace methods that is independent of the rank.

\paragraph{Effect of Subspace Selection Strategy}
Next, we fix the subspace rank at $r=16$ and conduct a comparison of the impact of different selection strategies of projection subspace on FLGT training the same way as we have experimented on FedAvg. We adopt 6 strategies, with the same 15 outer FLGT loops, where each client performs 100 local training steps per outer loop and aggregation is conducted every 20 steps. The 6 projection strategies are defined the same as those for SO-FedAvg experiments.

\begin{figure}[t]
    \center
    \includegraphics{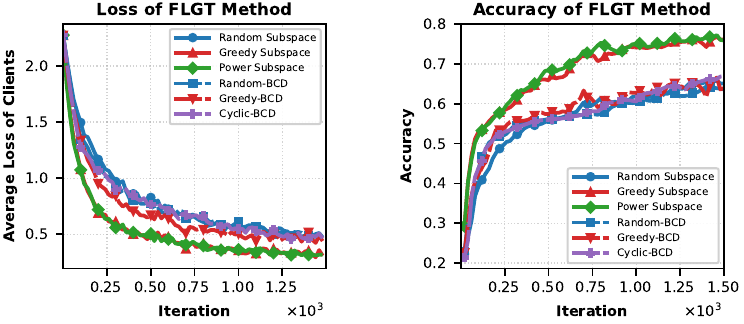} 
    \caption{\small Comparison plots of the training performance of the Subspace+FLGT algorithm with different subspace selection strategies. The plots show the trends of the average client loss and server test accuracy as functions of the total number of iterations.}\label{fig:so-flgt-method}
\end{figure}

As can be seen from figure~\ref{fig:so-flgt-method}, the Greedy Subspace and Power Subspace methods achieve the fastest iteration-wise convergence, while the BCD-type methods and the Random Subspace method converge more slowly. This indicates that the power iteration strategy in this experiment effectively approximates the subspace spanned by the top-$r$ singular vectors of the average gradient, significantly reducing the computational overhead at the server. The strategies of Greedy Subspace, Power Subspace and Greedy-BCD require each client to additionally upload gradients, which implies the increase in communication overhead, while strategies including Random Subspace do not need the gradient communication process. We conclude from the experiment results that SO-FLGT can also leverage the inherent low-rank structure in federated learning, demonstrating clear advantages when comprehensively considering communication overhead, computational cost, and convergence speed, thus offering significant potential for research and application.

Furthermore, compared to the experimental results in previous section using FedAvg under the i.i.d. data distribution, the convergence curves in this experiment are less smooth, which is attributed to the client drift phenomenon during training. Nevertheless, SO-FLGT still achieves convergence with a satisfactory speed, indicating that gradient tracking is essential for maintaining subspace update quality under non-i.i.d. data.

In summary, the results of numerical experiments demonstrate that the model parameter space in FL tasks has a certain degree of redundancy and can be well solved by low-dimensional optimization methods. Meanwhile, instead of transmitting all model parameters between the server and clients, the server issues projection matrices and model parameters while clients upload low-dimensional parameters, which means introducing Matrix Subspace Optimization into FL significantly reduces the communication overhead.

\section{Decentralized Optimization}
In recent years, decentralized optimization has become an important approach to problems concerning multi-node systems. 
Unlike traditional centralized methods that rely on a central server for data collection and global computation, decentralized methods require a connected network of nodes to solve a global optimization problem collaboratively through local computations and communication with immediate neighbors only, effectively mitigating the impact of node failure and alleviating communication bottlenecks at the single central node. 
Such methods facilitate tasks including collaborative AI training, which rely on large distributed datasets, enabling efficient learning while preserving data privacy at the same time.

We consider a decentralized consensus optimization problem over a network of $N$ nodes. We assume that the nodes communicate over a connected network represented by a sparse doubly stochastic matrix $A \in \mathbb{R}^{N \times N}$. Each node $i$ maintains a local copy $W_i\in\mathbb{R}^{m\times n}$ of the global matrix variable $W$ and holds a local private objective function $f_i:\mathbb{R}^{m\times n} \to \mathbb{R}$. The general goal is to solve the following problem:
\begin{equation}
\min_{W \in \mathbb{R}^{m\times n}} F(W) := \frac{1}{N}\sum_{i=1}^N f_i(W).
\end{equation}

All local parameters, gradients, and tracking variables are matrices. Communication is performed through weighted averaging of the corresponding local matrices.

Using local copies, the problem can be equivalently written as:
\begin{equation}\label{eq:do-prob}
\min_{W_1,\ldots,W_N} \frac{1}{N}\sum_{i=1}^N f_i(W_i)
\quad \text{s.t. } W_1=\cdots=W_N.
\end{equation}

Several algorithms has been developed for decentralized optimization problems that emphasize consensus among nodes. Decentralized gradient descent (DGD)~\cite{nedic2009distributed} is one of the most intuitive one, yet it fails to achieve exact convergence under fixed learning rates~\cite{yuan2016convergence}. EXTRA~\cite{shi2015extra} overcomes this issue by leveraging the state and gradient differences between consecutive iterations, achieving exact convergence with a constant learning rate.
Gradient tracking (GT) method as proposed in \cite{nedic2017achieving} further introduces an auxiliary variable to dynamically track and estimate the global average gradient, achieving rapid and accurate convergence.

This chapter presents the integration of \textit{Matrix Subspace Optimization} framework with DGD and GT with update derivation and numerical experiments. Although the vanilla versions of the algorithms are proposed in a Combine-then-Adapt (CTA) paradigm, this paper adopts the Adapt-then-Combine (ATC) structure. According to Sayed~\cite{sayed2014diffusion}, ATC and CTA share a theoretical duality, yet ATC incorporates new local gradient information before diffusion, exhibiting  stability and faster convergence under the same learning rate.

\subsection{Decentralized Gradient Descent (DGD)}
\subsubsection{Preliminary}
Decentralized Gradient Descent (DGD) serves as a fundamental deterministic algorithm for distributed optimization over sparse networks. In each iteration, a node updates its state by performing a local gradient descent step (adapt phase) and subsequently computing a weighted average of its immediate neighbors' states to promote network consensus (combine phase). The algorithm is a combination of the two phases.

We denote the local parameter matrix at node $i$ in iteration $k$ by $W_i^k$. The corresponding local gradient is $\nabla f_i(W_i^k)$. The DGD update rule can be expressed as:
  \begin{equation}
    W_i^{k+1} = \sum_{j=1}^N a_{ij}\left(W_j^k-\alpha\nabla f_j(W_j^k)\right).
  \end{equation}
Despite its computational simplicity, DGD exhibits a well-documented limitation: it fails to achieve exact convergence to the global optimum when operating with a fixed learning rate $\alpha$. The inherent variance among the local objective functions creates a non-vanishing steady-state error, forcing a trade-off between the convergence rate and the exactness of the final solution.\cite{yuan2016convergence}
\subsubsection{Update Derivation}
To apply the proposed framework to DGD, we first reintroduce the problem with subspace variable matrices.

In outer iteration $k$, let $B_i^k\in\mathbb{R}^{r\times n}(r\ll m)$ denote the subspace coefficient matrix held by node $i$, and let $P^k\in\mathbb{R}^{m\times r}$ denote the projection matrix.
We wish to obtain a subspace formulation of \ref{eq:do-prob}. Let $G(B):=\sum_{i=1}^N g_i(B_i^k)$ be the objective to be minimized with $g_i:\mathbb{R}^{r\times n}\to\mathbb{R}$ being the local function at node $i$ with respect to the local subspace coefficient matrix.

We first substitute $W_i^k+P^kB_i^k$ for the local full-dimension matrix $W_i$ in \ref{eq:do-prob} and then add a regularization term for the subspace variables. The local subspace objective is $g_i(B_i^k):=f_i(W_i^k+P^kB_i^k)+\frac{1}{2\eta^k}\|B_i^k\|_F^2$. This leads us to the resulting subproblem as follows:
\begin{equation}\label{eq:do-subprob}
    \min_{B} \quad G(B):=\sum_{i=1}^N g_i(B_i^k).
\end{equation}

Taking the gradient of $g_i(B_i^k)$ with respect to $B_i^k$, we obtain the local subspace gradient matrix:
\begin{equation}\label{eq:do-sub-derivative}
    \nabla g_i(B_i^k)=(P^k)^\top\nabla f_i(W_i^k+P^kB_i^k)+\frac{1}{\eta^k}B_i^k.
\end{equation}

We implement the \textit{Matrix Subspace Optimization} framework for DGD by solving the above subproblem according to the DGD update rules. Within the inner loops of the $k$-th outer iteration, nodes are first updated by gradient descent with a fixed learning rate $\alpha$. Subsequently, the variables are communicated and combined across the network. Denote $B_i^k$ at the $t$-th inner loop by $B_i^{k,t}$, the matrix update gives:
\begin{equation}\label{eq:do-dgd-1}
    B_i^{k,t+1}=\sum_{j=1}^N a_{ij}\left(B_j^{k,t}-\alpha\nabla g_j(B_j^{k,t})\right).
\end{equation}
At the end of each outer iteration, the local full-dimensional variable matrix $W_i^k$ is updated:
\begin{equation}\label{eq:do-dgd-2}
    W_i^{k+1}=W_i^k+P^kB_i^{k,T_{\max}}.
\end{equation}
A recursion of \eqref{eq:do-sub-derivative}, \eqref{eq:do-dgd-1}, and \eqref{eq:do-dgd-2} gives the SO-DGD update rule, which is summarized in Algorithm \ref{alg:SO-DGD}. 
\begin{algorithm}[t]
    \caption{Subspace Optimization with Decentralized Gradient Descent (SO-DGD)}\label{alg:SO-DGD}
    \textbf{Input:} Number of nodes $N$; $\{W_i^0\}_{i=1}^N$; connected network $A=(a_{ij}) \in \mathbb{R}^{N \times N}$; outer steps $K$; proximal coefficients $\{\eta^k\}$; maximum inner steps $T_{\max}$.\\
    \textbf{Output:} $\{W_i^K\}_{i=1}^N$.
        \begin{algorithmic}[1]
        \For{$k = 0, 1, \dots, K-1$}
            \State Sample or compute a subspace matrix $P^k$.
            \State Initialize subspace parameter matrix $B_i^{k,0}$ for each node $i=1,\cdots,N$.
            \For{$t = 0, 1, \dots, T_{\max}-1$}
                \State Compute subspace gradient for each node: $\nabla g_i(B_i^{k,t}) \leftarrow (P^k)^\top \nabla f_i(W_i^k + P^k B_i^{k,t}) + \frac{1}{\eta^k} B_i^{k,t}$.
                \State Update each node: $B_i^{k,t+1/2} \leftarrow B_i^{k,t} - \alpha \nabla g_i(B_i^{k,t})$.
                \State Communicate and combine at each node: $B_i^{k,t+1} \leftarrow \sum_{j=1}^N a_{ij} B_j^{k,t+1/2} $.
                \State \textbf{if} {$B_i$ converges for all nodes} \textbf{then} $\tilde{B}_i^{k} \gets B_i$ and \textbf{break}
            \EndFor
            \State Update full-dimension variable $W_i$ for each node:  $W_i^{k+1} \leftarrow W_i^k + P^k \tilde{B}_i^{k}$.
        \EndFor
    \end{algorithmic}
\end{algorithm}

The primary advantage of SO-DGD is communication and memory efficiency. In standard DGD, nodes exchange and store full-dimensional parameters in each iteration. For modern high-dimensional non-convex problems where the parameter matrix has $mn$ entries, this introduces an $\mathcal{O}(mn)$ per-iteration communication and memory cost. In contrast, our \textit{Matrix Subspace Optimization} framework confines high-frequency inner-loop consensus and tracking updates to coefficient matrices with $rn$ entries, with $r\ll m$. Therefore, nodes only transmit the coefficient matrices, reducing the communication burden from $\mathcal{O}(mn)$ to $\mathcal{O}(rn)$ per inner iteration. Furthermore, memory demand is reduced since nodes only store one full dimensional variable.

\subsubsection{Numerical Experiments}
\label{sec:ex}
We test SO-DGD on MNIST image classification tasks, the settings of which are introduced in \ref{MLP-setting}.

\begin{figure}[t]
    \centering
    \includegraphics{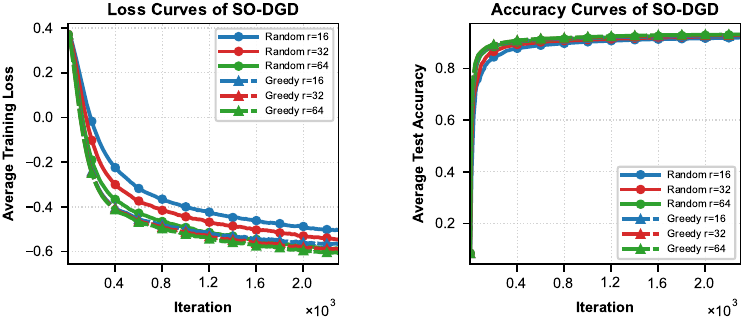}
    \caption{\small Convergence of Random and Greedy Subspace DGD at different ranks. Left: Average training loss across nodes. Right: Average test accuracy across nodes. 
   Both Random and Greedy Subspace DGD converge robustly even with significantly reduced ranks. All methods demonstrate high accuracy on the test set.}
    \label{fig:dgd-rank}
\end{figure}

We perform the tasks on an exponential topology (or a power-of-two circulant graph) of 20 nodes. The standard MNIST dataset is distributed among nodes, each node independently accessing the global training pool. At each iteration, every node samples a local mini-batch of 128 images. To measure global performance of the algorithm, we calculate the average training loss across all nodes using their current local batches. Meanwhile, we evaluate the global test accuracy across nodes by averaging the accuracy of all 20 local models across the entire MNIST test set. We then analyze the impact of subspace rank and selection. 
\paragraph{Effect of Subspace Rank} We first investigate the rank selection for the subspace. We test scenarios with subspace rank $r\in\{16,32,64\}$. We compare Random Subspace DGD where subspaces are randomly selected against Greedy Subspace DGD where projection matrices are selected via SVD decomposition of local gradients. To ensure a fair comparison, all experiments are conducted over a total of 2340 iterations, equivalent to 5 epochs on the MNIST dataset with a batch size of 128. For the Random Subspace algorithms, the training process is structured into 117 outer loops. Within each outer loop,  the projection matrix is updated and fixed, followed by 20 inner loop iterations where the parameters are optimized within the subspace. The learning rates are fixed at 0.01. The results are shown in Figure \ref{fig:dgd-rank}.

As figure \ref{fig:dgd-rank} shows, both Random and Greedy Subspace DGD show stable performance even when ranks are greatly reduced. As the subspace rank grows, the initial speed of convergence grows, but in smaller rank cases convergence can be eventually reached. This implies that both Random and Greedy Subspace DGD perform well on the task. As rank becomes larger, the initial convergence speed of the algorithm gets faster, which is intuitively clear since a larger subspace contains more information about the update direction.

\begin{figure}[t]
    \centering
    \includegraphics{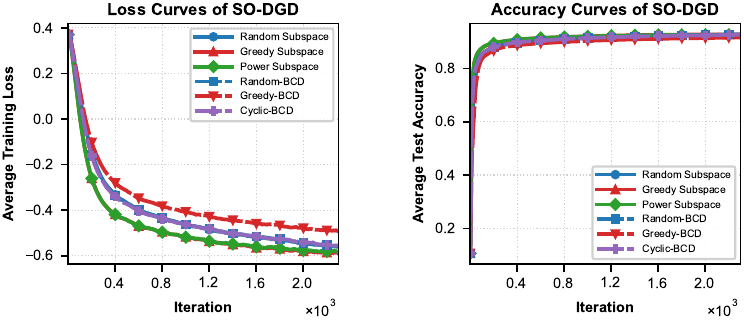}
    \caption{\small Comparison of Subspace Selection Strategy. All strategies select subspaces of rank $r = 32$. Left: Average training loss across nodes. Right: Average test accuracy across nodes. 
    All methods show considerable accuracy on test set. Greedy and Power Subspace DGD converge fast initially, while Random Subspace DGD converges in similar pace with Cyclic-BCD and Random-BCD methods. Greedy-BCD Subspace has the slowest rate of convergence. }
    \label{fig:dgd-p}
\end{figure}

Also, the performance of Greedy Subspace DGD seems to be always better than that of Random Subspace DGD. Moreover, Greedy Subspace DGD with different ranks shows similar performance. A possible explanation concerns heavy-tailed spectral distribution of gradients in the neural network, where the majority of update information is concentrated in a few top principal components. 
By targeting these highest-magnitude components, the Greedy Subspace DGD successfully captures almost all essential signals at a low rank of $r = 16$ . Expanding the rank to $r=32$ or $r=64$ yields only negligible additional gradient information.
However, Random Subspace DGD samples the subspace blindly, so at a highly compressed rank it is very likely to miss the most critical descent directions, leading to slower convergence.
\paragraph{Effect of Subspace Selection Strategy}
As the selection of subspace might impact the performance of the \textit{Matrix Subspace Optimization} methods, we further implement a subspace selection ablation test. Fixing the subspace rank at $r = 32$, we use different projection matrices $P$ indicating various subspace selection strategies. We test the six projection matrix selections that are detailed in \ref{subspace-choice}.
The training spans 5 epochs and the learning rates are fixed at $0.01$. The results are shown in Figure \ref{fig:dgd-p}.


As figure \ref{fig:dgd-p} shows, Greedy and Power Subspace DGD outperform other subspace selection strategies and the resulting loss curves nearly overlap. Random Subspace DGD gives performance comparable to that of Cyclic-BCD and Random-BCD methods, reaching convergence and achieving high accuracy eventually. Greedy-BCD underperforms all other methods, but still manages to converge.

Though Power and Greedy Subspace methods perform best, they require communication of full gradient information, which is not memory-efficient.
The poor performance of Greedy-BCD might be attributed to the mismatch in coordinate updates among nodes. Each node independently selects the best coordinate blocks based only on the local gradient. Due to differences in local datasets and mini-batch gradients, different nodes naturally end up selecting and transmitting different sets of coordinates, leading to a mismatch.

\subsection{Gradient Tracking (GT)}
\subsubsection{Preliminary}
To rapidly gain exact convergence with a constant learning rate, Gradient Tracking (GT) introduces a dynamic average consensus mechanism. By adopting a tracked global gradient estimator for the estimation of the local gradient during the adaptation phase, the algorithm ensures that the local update directions asymptotically align with the true global descent direction, enforcing consensus among nodes.

We denote the local gradient tracking matrix of node $i$ in iteration $k$ by $Y_i^k$. The tracker is initialized by the gradient evaluated at the starting point such that $Y_i^0=\nabla f_i(W_i^0)$. During each iteration $k$, the network synchronously updates both the local parameter matrices $W_i^k$ and the local tracking matrices $Y_i^k$. The update of $W_i^k$ follows a straightforward gradient descent rule:
 \begin{equation}
 W_i^{k+1}=\sum_{j=1}^N a_{ij}\left(W_j^k-\alpha Y_j^k\right).
 \end{equation}
$Y_i^k$ is estimated through network diffusion after an update with the difference between consecutive gradients:
\begin{equation}
Y_i^{k+1}=\sum_{j=1}^N a_{ij}\left(Y_j^k+\nabla f_j(W_j^{k+1})-\nabla f_j(W_j^k)\right).
\end{equation}
In this way, $Y_i^k$ asymptotically tracks the average of local gradients, which coincides with the global gradient under consensus. This mechanism eliminates the bias in standard decentralized gradient descent and guarantees exact convergence to the optimal solution in deterministic settings.
\subsubsection{Update Derivation}
Our aim here is to solve a decentralized consensus problem, so the goal is the subproblem described in \ref{eq:do-prob}.

We apply the \textit{Matrix Subspace Optimization} framework to gradient tracking. Similarly to the vanilla gradient tracking algorithm, we first introduce a local auxiliary tracking matrix $Y_i^{k,t}$ with the same shape as $B_i^{k,t}$ to track the global average of projected gradients in the subspace. In the inner loops of the $k$-th outer iteration, we first perform a local update along the tracked direction matrix $Y_i^{k,t}$:
\begin{equation}\label{eq:gt-update-1}
    B_i^{k,t+1/2}=B_i^{k,t}-\alpha Y_i^{k,t},\qquad
    B_i^{k,t+1}=\sum_{j=1}^N a_{ij}B_j^{k,t+1/2}.
\end{equation}
The tracking matrix is initialized as the local subspace gradient $Y_i^{k,0}=\nabla g_i(B_i^{k,0})$. In each inner step $t$, it is updated by aggregating the neighborhood trackers and injecting the gradient difference after the local update:
\begin{equation}\label{eq:gt-update-2}
    Y_i^{k,t+1}=\sum_{j=1}^N a_{ij}\left(Y_j^{k,t}+\nabla g_j(B_j^{k,t+1})-\nabla g_j(B_j^{k,t})\right).
\end{equation}
A recursion of \ref{eq:gt-update-1} and \ref{eq:gt-update-2} , along with the global update as in \ref{eq:do-dgd-2} with a local gradient calculated by \ref{eq:do-sub-derivative} outlines the update of Subspace Optimization with Gradient Tracking, which is detailed in Algorithm \ref{alg:SO-GT}.

\begin{algorithm}[t]
    \caption{Subspace Optimization with Gradient Tracking (SO-GT)}\label{alg:SO-GT}
    \textbf{Input:} Number of nodes $N$; $\{W_i^0\}_{i=1}^N$; connected network $A=(a_{ij}) \in \mathbb{R}^{N \times N}$; outer steps $K$; proximal coefficients $\{\eta^k\}$; maximum inner steps $T_{\max}$.\\
    \textbf{Output:} $\{W_i^K\}_{i=1}^N$.
    \begin{algorithmic}[1]
        \For{$k = 0, 1, \dots, K-1$}
            \State Sample or compute a subspace matrix $P^k$.
            \State Initialize subspace parameter matrix $B_i^{k,0}$ for each node $i = 1,\cdots, N$.
            \State Initialize auxiliary tracking matrix for each node: $Y_i^{k,0} \leftarrow \nabla g_i (B_i^{k,0})$.
            \For{$t = 0, 1, \dots T_{\max}-1$}
                \State Update each node: $B_i^{k,t+1/2} \leftarrow B_i^{k,t} - \alpha Y_i^{k,t}$.
                \State Communicate and combine at each node: $B_i^{k,t+1} \leftarrow \sum_{j=1}^N a_{ij} B_j^{k,t+1/2}$.
                \State Compute subspace gradient $\nabla g_i (B_i^{k,t+1})$ for each node.
                \State Gradient tracking for each node: $Y_i^{k,t+1} \leftarrow \sum_{j=1}^N a_{ij} \left(Y_j^{k,t} + \nabla g_j(B_j^{k,t+1}) - \nabla g_j(B_j^{k,t})\right)$.
                \State \textbf{if} {$B_i$ converges for all nodes} \textbf{then} $\tilde{B}_i^{k} \gets B_i$ and \textbf{break}
            \EndFor
            \State Update full-dimension variable for each node: $W_i^{k+1} \leftarrow W_i^k + P^k \tilde{B}_i^{k}$.
        \EndFor
    \end{algorithmic}
\end{algorithm}

\subsubsection{Numerical Experiments}
We perform MNIST image classification tasks as described in \ref{MLP-setting} on an exponential graph with 20 nodes to evaluate the efficacy of SO-GT. Each node independently samples a 128-image mini-batch from the MNIST training set in each iteration step. The overall performance is measured by the average local training loss across the whole topology and the average accuracy of the nodes on the standard test set.

\begin{figure}[t]
    \centering
    \includegraphics{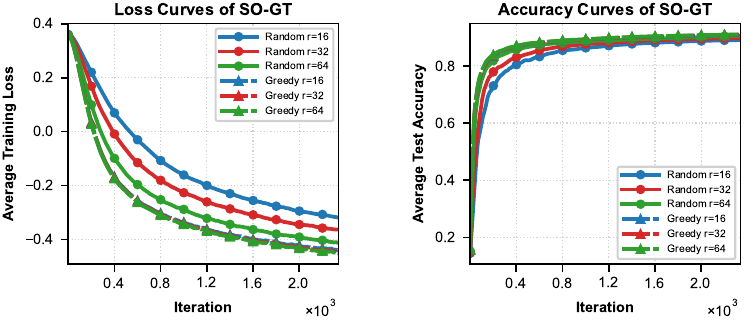}
    \caption{\small Greedy and Random Subspace GT with different subspace ranks. Left: Average training loss across nodes. Right: Average test accuracy across nodes. 
 Random Subspace GT methods converge robustly with low rank subspaces, although they converge slower than Greedy counterparts initially. All methods demonstrate high accuracy.}
    \label{fig:gt-rank}

    \par\medskip
    \includegraphics{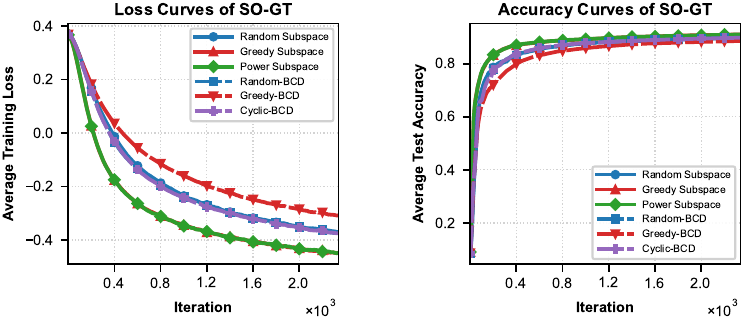}
    \caption{\small Comparison of Subspace Selection Strategy. All subspaces have rank $r = 32$. Left: Average training loss across nodes. Right: Average test accuracy across nodes. 
    All methods demonstrate high accuracy on test set. Greedy and Power Subspace GT converge the fastest, while Random Subspace GT later catches up, converging in similar pace with Cyclic-BCD and Random-BCD methods. The Greedy-BCD method converges the slowest. }
    \label{fig:gt-p}
\end{figure}

\paragraph{Effect of Subspace Rank} The same task for SO-DGD is performed here to compare the performance of SO-GT with different subspace ranks. The ranks are $r\in\{16,32,64\}$ out of a total rank of $128$ and the learning rates are fixed at $0.01$. The results are shown in Figure \ref{fig:gt-rank}.

As figure \ref{fig:gt-rank} depicts, both Greedy and Random Subspace GT methods perform well. As the rank of the selected subspace grows, the initial convergence speed increases and convergence is achieved in all cases. The Greedy Subspace GT converges faster than the Random Subspace GT with the same subspace dimension. Meanwhile, greedy subspace methods with various ranks achieve close convergence rates.

\paragraph{Effect of Subspace Selection Strategy} To study the effect of the subspace selection strategy, we try the same six ways to sample $P$ as described in \ref{subspace-choice}. The subspace rank is fixed at $r = 32$ and the learning rates are fixed at $0.01$. The overall training process is partitioned into 117 outer loops. In each outer loop, a subspace is selected and fixed in the subsequent 20 inner iterations. The results are shown in Figure \ref{fig:gt-p}.

As figure \ref{fig:gt-p} shows, the Greedy and Power Subspace  methods have the fastest convergence rates among all methods. The Random Subspace, Cyclic-BCD and Random-BCD Subspace methods later catch up, and they show almost overlapping loss curves. The Greedy-BCD method performs the worst, but still converges in the end, achieving high accuracy. A similar phenomenon has been observed in SO-DGD tests.

In summary, SO-DGD and SO-GT performs well in reduced subspaces. This implies that the \textit{Matrix Subspace Optimization} framework might be promising for optimization in decentralized settings.

\FloatBarrier
\section{Bilevel Optimization}
Bilevel optimization \cite{bracken1973mathematical, sinha2017review} involves problems where one optimization task (the outer problem) is constrained by the solution of another (the inner problem). This hierarchical structure naturally models applications such as hyperparameter optimization \cite{franceschi2018bilevel}, meta-learning \cite{finn2017model}, and neural architecture search \cite{liu2018darts}. Formally, a bilevel problem can be written as:
\begin{equation}\label{eq:bilevel_general}
    \begin{array}{ll}
        \min_{\bx \in \mathbb{R}^{d}} & h(\bx) = F(\bz^*(\bx),\bx), \\
        \text{s.t.} & \bz^*(\bx) \in \argmin_{\bz \in \mathbb{R}^{p}} G(\bz,\bx),
    \end{array}
\end{equation}
where $F$ and $G$ are the outer and inner objective functions, respectively. Solving such problems is challenging due to the implicit dependence of the optimal inner variable $\bz^*(\bx)$ on the outer variable $\bx$, which requires computing the hypergradient $\nabla h(\bx)$.

\subsection{Stochastic One-stage Bilevel Algorithm (SOBA)}
\subsubsection{Preliminary}
To mitigate the computational cost of nested loops, single-loop algorithms have gained popularity. Among them, the Stochastic One-stage Bilevel Algorithm (SOBA) \cite{dagreou2022framework} stands out for its efficiency. The core idea of SOBA is to update the outer variable $\bx$, the inner variable $\bz$, and an auxiliary Lagrange multiplier (or adjoint variable) $\bv$ simultaneously in a single loop. 

Specifically, the hypergradient is given by $\nabla h(\bx) = \nabla_2 F(\bz^*, \bx) - \nabla_{21}^2 G(\bz^*, \bx) [\nabla_{11}^2 G(\bz^*, \bx)]^{-1} \nabla_1 F(\bz^*, \bx)$. SOBA introduces the variable $\bv$ to estimate the term $[\nabla_{11}^2 G(\bz^*, \bx)]^{-1} \nabla_1 F(\bz^*, \bx)$ by viewing it as the solution to a linear system, or equivalently, the minimizer of a quadratic approximation. At each iteration $k$, SOBA performs three steps: 
(1) a descent step on $\bz$ to minimize the inner objective $G$; 
(2) a descent step on $\bv$ to solve the linear system involving the Hessian $\nabla_{11}^2 G$; 
and (3) a descent step on $\bx$ using the estimated hypergradient constructed from current $\bz$ and $\bv$. This approach relies on the distinct timescales of the step sizes ($\rho_k \gg \gamma_k$) to ensure that $\bz$ and $\bv$ converge faster than $\bx$, thereby tracking the optimal trajectory.

\subsubsection{Update Derivation}
For notation clarity, we present the bilevel derivation in vector form. The corresponding matrix form follows naturally by replacing the outer vector variable with a matrix variable and applying the same subspace parametrization. This keeps the hypergradient, Hessian-vector, and mixed-Jacobian notation readable without changing the subspace-reduction principle.

To alleviate the computational and memory burden in high-dimensional bilevel problems, we instantiate the subspace optimization framework with the SOBA algorithm. Assuming the inner problem has a unique solution, we reformulate the bilevel update within a selected low-dimensional subspace. At outer iteration $k$, we seek an update direction $\tilde{\bb}^k$ in a subspace defined by $\bP^k \in \mathbb{R}^{d \times r}$.

We define the subspace-restricted inner and outer functions as $\widetilde{G}_k(\bz,\bb) := G(\bz, \bx^k + \bP^k \bb)$ and $\widetilde{F}_k(\bz,\bb) := F(\bz, \bx^k + \bP^k \bb)$. Let $\tilde{h}_k(\bb) := \widetilde{F}_k(\tilde{\bz}_k^*(\bb),\bb)$. The subproblem for finding the update $\tilde{\bb}^k$ becomes a lower-dimensional regularized bilevel problem:
\begin{equation}
    \begin{array}{ll} 
        \min_{\bb \in \mathbb{R}^{r}} & \tilde{h}_k(\bb) + \frac{1}{2\eta^k}\|\bb\|^2, \\ 
        \text{s.t.} & \tilde{\bz}_k^*(\bb) \in \argmin_{\bz \in \mathbb{R}^{p}} \widetilde{G}_k(\bz,\bb).
    \end{array}
\end{equation}
By applying the SOBA strategy to this reduced problem, we update the inner variable $\bz$, the auxiliary variable $\bv$ (for the subspace Hessian inverse), and the subspace parameter $\bb$ concurrently. The update rules for the inner loop $t$ are given by:
\begin{subequations}
    \begin{align}
        \bz^{k,t+1} &= \bz^{k,t} - \rho_t \nabla_1 \widetilde{G}_k(\bz^{k,t},\bb^{k,t}), \\
        \bv^{k,t+1} &= \bv^{k,t} - \rho_t (\nabla_{11}^2 \widetilde{G}_k(\bz^{k,t},\bb^{k,t}) \bv^{k,t} + \nabla_1 \widetilde{F}_k(\bz^{k,t},\bb^{k,t})), \\
        \bb^{k,t+1} &= \bb^{k,t} - \gamma_t \left(\nabla_{21}^2 \widetilde{G}_k(\bz^{k,t},\bb^{k,t}) \bv^{k,t} + \nabla_2 \widetilde{F}_k(\bz^{k,t},\bb^{k,t}) + (\eta^k)^{-1}\bb^{k,t}\right).
    \end{align}
\end{subequations}
Crucially, the terms $\nabla_{21}^2 \widetilde{G}_k$ and $\nabla_2 \widetilde{F}_k$ involve differentiation with respect to the low-dimensional $\bb$ rather than the full-dimensional $\bx$. This significantly reduces the size of the Jacobian-vector products from $d \times p$ (or $d$) to $r \times p$ (or $r$), leading to memory savings. The SO-SOBA procedure is detailed in Algorithm \ref{alg:so-soba}.

\begin{algorithm}[t]
    \caption{Subspace Optimization with SOBA (SO-SOBA)}\label{alg:so-soba}
    \textbf{Input:} Initial outer variable $\bx^0$; inner variable $\bz^{0,0}$; auxiliary variable $\bv^{0,0}$; outer steps $K$; proximal coefficients $\{\eta^k\}$; maximum inner steps $T_{\max}$; step sizes $\{\rho_t,\gamma_t\}$.\\
    \textbf{Output:} $\bx^K$.
    \begin{algorithmic}[1]
        \For{$k = 0, 1, \dots, K-1$}
            \State Sample or compute a subspace matrix $\bP^k \in \mathbb{R}^{d \times r}$.
            \State Initialize subspace parameter $\bb^{k,0} \gets \mathbf{0}$.
            \For{$t = 0, 1, \dots, T_{\max}-1$}
                \State Update inner variable $\bz^{k,t+1} \leftarrow \bz^{k,t} - \rho_t \nabla_1 \widetilde{G}_k(\bz^{k,t},\bb^{k,t})$.
                \State Update auxiliary variable $\bv^{k,t+1} \leftarrow \bv^{k,t} - \rho_t (\nabla_{11}^2 \widetilde{G}_k(\bz^{k,t},\bb^{k,t}) \bv^{k,t} + \nabla_1 \widetilde{F}_k(\bz^{k,t},\bb^{k,t}))$.
                \State Update subspace variable $\bb^{k,t+1} \leftarrow \bb^{k,t} - \gamma_t (\nabla_{21}^2 \widetilde{G}_k(\bz^{k,t},\bb^{k,t}) \bv^{k,t} + \nabla_2 \widetilde{F}_k(\bz^{k,t},\bb^{k,t}) + (\eta^k)^{-1}\bb^{k,t})$.
                \State \textbf{if} {$\bb$ converges} \textbf{then break}
            \EndFor
            \State Update outer variable: $\bx^{k+1} \leftarrow \bx^k + \bP^k \bb^{k,T_{\max}}$.
            \State (Optional) Carry over $\bz$ and $\bv$ for the next iteration.
        \EndFor
    \end{algorithmic}
\end{algorithm}

\subsubsection{Numerical Example}
At a high level, we still use MNIST as the underlying data source as in Section \ref{MLP-setting}, but the task here is bilevel data hyper-cleaning rather than direct supervised training. We use $N=2000$ noisy training samples, corrupt 30\% of their labels, and evaluate on a clean validation set of size $500$. The outer variable is the sample-weight vector $\lambda \in \mathbb{R}^{2000}$, while the inner model is a logistic regressor with parameter $\theta \in \mathbb{R}^{10 \times 784}$ (7840 parameters). We report validation loss and validation accuracy for all experiments.

\paragraph{Effect of Subspace Dimension}
We investigate the impact of the subspace dimension $r$ on the convergence of the bilevel optimization process. Specifically, we test SO-SOBA with dimensions $r \in \{10, 50, 200\}$ and compare them with the standard SOBA, which corresponds to the full outer dimension $r=2000$. SO-SOBA updates a low-dimensional vector $\bb \in \mathbb{R}^r$ within a random subspace that is re-sampled every 100 iterations.

\begin{figure}[t]
    \center
    \includegraphics{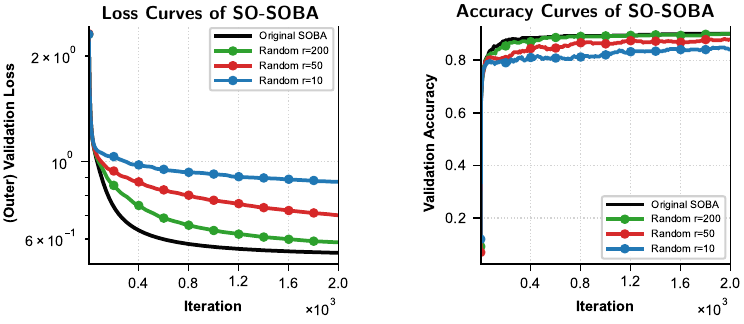} 
    \caption{\small Subspace Optimization with SOBA (SO-SOBA) with different subspace dimensions. Left: Validation Loss. Right: Validation Accuracy. SO-SOBA demonstrates robust convergence even with significantly reduced dimensions.}\label{fig:so-soba-rank}
\end{figure}

As shown in Figure \ref{fig:so-soba-rank}, SO-SOBA exhibits remarkable robustness to dimension reduction. Even with a subspace dimension as low as $r=10$ (only 0.5\% of the full outer dimension), SO-SOBA effectively reduces the validation loss while improving validation accuracy. Increasing the dimension to $r=50$ and $r=200$ further accelerates convergence, closely matching the trajectory of the original algorithm. This result highlights the redundancy in the hyperparameter space for data-cleaning tasks and confirms that the subspace bilevel instantiation can efficiently navigate the optimization landscape through low-dimensional random subspaces, offering a memory-efficient alternative for large-scale bilevel problems without sacrificing solution quality.

\paragraph{Effect of Subspace Selection Strategy}
Beyond reducing the subspace dimension, we investigate how different subspace construction strategies affect the optimization trajectory in bilevel problems. We follow the four-strategy comparison introduced for SO-Newton in \ref{para:so-newton-subspace-selection}: random subspace, random-CD, cyclic-CD, and greedy-CD. Here the subspace dimension is fixed to $r=100$ (5\% of the full outer dimension), and all methods update the same outer variable $\lambda \in \mathbb{R}^{2000}$. Thus, SO-SOBA optimizes a 100-dimensional subspace, while each CD baseline updates exactly 100 coordinates of $\lambda$ per iteration.

The only bilevel-specific adaptation is the oracle used by greedy-CD: instead of ranking coordinates by the full gradient of a single-level objective, it selects the top-$r$ coordinates using the full-dimensional outer-step direction used by SOBA, i.e., the hypergradient surrogate on $\lambda$. This requires computing the full $d$-dimensional direction and therefore remains a computationally expensive oracle baseline.

\begin{figure}[t]
    \center
    \includegraphics{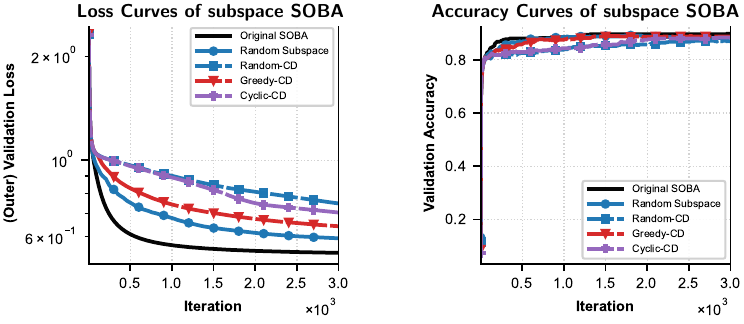}
    \caption{\small Comparison of Subspace Selection Strategy. Left: Validation Loss. Right: Validation Accuracy. Random subspace significantly outperforms random-CD and cyclic-CD, capturing coupled variable interactions and approaching the performance of the full-dimensional baseline and even outperforms the expensive greedy-CD oracle.}\label{fig:so-soba-strategy}
\end{figure}

As illustrated in Figure \ref{fig:so-soba-strategy}, standard axis-aligned methods (random-CD and cyclic-CD) exhibit noticeably slower convergence, higher loss plateaus, and weaker validation accuracy. In tasks like Data Hyper-Cleaning, the outer variables (sample weights) are highly coupled through the inner model's training dynamics. CD methods struggle to capture these complex interactions because they only update isolated coordinates. 

Greedy-CD improves upon random-CD and cyclic-CD but still underperforms the random-subspace strategy. Moreover, greedy-CD strictly violates memory constraints because its selection step requires the full-dimensional hypergradient surrogate. In contrast, the Gaussian random subspace naturally encompasses dense linear combinations of all variables. This allows SO-SOBA to explore coupled update directions efficiently, achieving a convergence profile that approaches the original full-dimensional SOBA while outperforming all coordinate-descent variants, all while remaining strictly memory-efficient.

\FloatBarrier
\section{Large Language Model Pre-training and Fine-tuning}
\subsection{LLM Training and the Memory Wall}
The key memory components involved in the training process consist of four primary elements: model parameters, optimizer states, gradients, and activations. The model parameters store the trainable weights required for optimization. In the case of the Adam optimizer, the optimizer states are given by the first- and second-moment estimates, denoted as $\bM$ and $\bV$ respectively. The gradient term corresponds to the memory cost associated with $\nabla F(\bW; \xi)$. With the Adam optimizer, both the optimizer states and the gradients are determined by the number of trainable parameters.

Another significant memory cost arises from activations, which represent the intermediate values computed during forward propagation. Unlike the optimizer states and gradients, activation memory depends on multiple factors, including model size, batch size, and sequence length. For small batch sizes, the optimizer states constitute a substantial portion of the memory usage; in contrast, for large batch sizes, activations dominate and can account for nearly the entire memory footprint. This ``memory wall'' limits the scale of models and batch sizes that can be trained in practice.

\subsection{How can subspace optimization help with memory-efficient LLM training?}

\begin{table}[t]
\centering
\small
\caption{\small Validation perplexity ($\downarrow$) and estimated BF16 optimizer-state memory in GB ($\downarrow$) during LLaMA pre-training. SO-Adam consistently outperforms LoRA and achieves performance comparable to the full-rank Adam.}\label{tab:llama-pretrain}
\vspace{0.5em}
\begin{threeparttable}
\begin{tabular}{lcccc}
\toprule
\textbf{Algorithm} & \textbf{LLaMA-60M} & \textbf{LLaMA-130M} & \textbf{LLaMA-350M} & \textbf{LLaMA-1B} \\ 
\midrule
Adam  & 34.06 (0.22G)   & 25.08 (0.50G)  & 18.80 (1.37G) & 15.56 (4.99G) \\ 
\midrule
GaLore      & 34.88 (0.14G)   & 25.36 (0.27G)  & 18.95 (0.49G) & \textbf{15.64} (1.46G) \\ 
LoRA       & 34.99 (0.16G)   & 33.92 (0.35G)  & 25.58 (0.69G) & 19.21 (2.27G)  \\ 
\textbf{SO-Adam}         & \textbf{34.55} (0.14G)  & \textbf{25.34} (0.27G) & \textbf{18.86} (0.49G) & 15.86 (1.46G) \\ 
\midrule
\textbf{Rank $r$ / Dim $m$} & 128 / 256  & 256 / 768  & 256 / 1024 & 512 / 2048 \\ 
\bottomrule
\end{tabular}
\end{threeparttable}
\end{table}

\begin{table}[t]
\centering
\small
\caption{\small Pre-training results for the LLaMA-7B model on the C4 dataset (evaluated at 50K steps). SO-Adam and GaLore are configured with a batch size of 16, while Adam uses a batch size of 8 due to memory constraints.}\label{tab:pretrain-llama-7b}
\vspace{0.5em}
\begin{threeparttable}
\begin{tabular}{l c c c}
\toprule
\textbf{Method}  & \textbf{Memory (GB)}  &  \textbf{Training Time (h)} &  \textbf{Perplexity} \\
\midrule
Adam  & 78.92  & 216 & \textbf{15.43} \\
GaLore & 75.33  & 134 & 15.59 \\
\textbf{SO-Adam} & \textbf{54.81}  & \textbf{64} & 15.99 \\
\bottomrule
\end{tabular}
\end{threeparttable}
\end{table}

\begin{table}[t]
    \centering
    \small
    \caption{\small Evaluation of fine-tuning methods on the GLUE benchmark using RoBERTa-Base.}\label{tab:glue-finetuning}
    \vspace{0.5em}
    \begin{threeparttable}
    \resizebox{\textwidth}{!}{
    \begin{tabular}{l c c c c c c c c c} 
    \toprule
    \textbf{Method} & \textbf{CoLA} & \textbf{STS-B} & \textbf{MRPC} & \textbf{RTE} & \textbf{SST2} & \textbf{MNLI} & \textbf{QNLI} & \textbf{QQP} & \textbf{Avg} \\
    \midrule
    Adam  & 62.24 & 90.92 & 91.30 & 79.42 & 94.57 & 87.18 & 92.33 & 92.28 & 86.28 \\
    \midrule
    \multicolumn{10}{l}{\textit{Subspace Rank = 4}} \\
    \quad GaLore  & 60.35 & \textbf{90.73} & \textbf{92.25} & \textbf{79.42} & 94.04 & \textbf{87.00} & 92.24 & 91.06 & 85.89 \\
    \quad LoRA  & 61.38 & 90.57 & 91.07 & 78.70  & 92.89 & 86.82 & 92.18 & \textbf{91.29} & 85.61 \\
    \quad \textbf{SO-Adam}  & \textbf{62.47} & 90.62 & \textbf{92.25} & 78.70 & \textbf{94.84} & 86.67 & \textbf{92.29} & 90.94 & \textbf{86.10} \\
    \midrule
    \multicolumn{10}{l}{\textit{Subspace Rank = 8}} \\
    \quad GaLore   & 60.06 & \textbf{90.82} & 92.01 & \textbf{79.78} & 94.38 & \textbf{87.17} & 92.20 & 91.11 & 85.94 \\
    \quad LoRA  & 61.83 & 90.80 & 91.90 & 79.06  & 93.46 & 86.94 & 92.25 & 91.22 & 85.93 \\
    \quad \textbf{SO-Adam}  & \textbf{64.62} & 90.71 & \textbf{93.56} & 79.42 & \textbf{95.18} & 86.96 & \textbf{92.44} & \textbf{91.26} & \textbf{86.77} \\
    \bottomrule
    \end{tabular}
    }
    \end{threeparttable}
\end{table}

\begin{table}[t]
\centering
\small
\caption{\small Fine-tuning test accuracy ($\uparrow$) of LLaMA and OPT models on WinoGrande and COPA datasets.}\label{tab:finetune-wincopa}
\vspace{0.5em}
\begin{threeparttable}
\begin{tabular}{l cccc cccc}
\toprule
\multirow{2}{*}{\textbf{Model}} & \multicolumn{4}{c}{\textbf{WinoGrande}} & \multicolumn{4}{c}{\textbf{COPA}} \\
\cmidrule(lr){2-5} \cmidrule(lr){6-9}
 & Adam & LoRA & GaLore & \textbf{SO-Adam} & Adam & LoRA & GaLore & \textbf{SO-Adam} \\
\midrule
LLaMA-7B    & 64.4 & 70.9 & 70.9 & \textbf{71.0} & 84.0    & 84.0    & 85.0    & \textbf{86.0}    \\
LLaMA-13B   & 73.3 & \textbf{76.6} & 74.6 & 74.7 & 90.0    & \textbf{92.0}    & \textbf{92.0}    & \textbf{92.0}    \\
\midrule
OPT-1.3B    & 60.4    & 57.3    & 58.3    & \textbf{58.9}    & 76.0 & 73.0 & 72.0 & \textbf{74.0} \\
OPT-6.7B    & 62.2    & 64.7    & 66.8    & \textbf{69.2}    & 78.0 & 80.0 & 80.0 & \textbf{82.0} \\
\bottomrule
\end{tabular}
\end{threeparttable}
\end{table}

Unlike other memory-efficient methods (e.g., GaLore, Adam-mini, Apollo) that primarily focus on reducing the memory usage of optimizer states, the subspace optimization method additionally achieves substantial savings in gradient and activation memory requirements.

\textbf{Memory for optimizer states.} When solving the matrix subspace subproblem, the reduced dimensionality of the subproblem substantially decreases the memory requirements for optimizer states. For instance, the memory required for both the first- and second-moment estimates in each subproblem is $r_\ell n_\ell$ parameters for the $\ell$-th layer, which is significantly lower than the $m_\ell n_\ell$ memory overhead in the standard Adam optimizer.

\textbf{Memory for gradients.} For the matrix subspace subproblems, it is sufficient to compute the gradient with respect to $\bB$, i.e., $\nabla_\bB F(\bW^{k} + \bP^{k} \bB; \xi)$, rather than the full-dimensional gradient $\nabla_\bW F(\bW^k; \xi)$ with respect to the original weight matrix $\bW$. This leads to considerable memory savings in gradient storage.

\textbf{Memory for activations.} The subspace optimization method reduces not only the memory usage for gradients but also the memory required to store activations. Consider the $\ell$-th layer of a neural network. In standard training with Adam, the forward pass is $Z_\ell = Y_\ell W_\ell$ and the loss depends on $Z_\ell$; in the subspace formulation, the forward pass becomes $Z_\ell = Y_\ell (W_\ell + P_\ell B_\ell)$. Here, $Y_\ell \in \mathbb{R}^{s_\ell \times m_\ell}$ denotes the output of the previous layer, i.e., the activation. In the backward pass, the weight gradient in Adam is computed as $\frac{\partial y}{\partial W_\ell} = Y_\ell^\top \frac{\partial y}{\partial Z_\ell}$ for Adam, which requires storing the full activation $Y_\ell$. In contrast, the subspace formulation computes $\frac{\partial y}{\partial B_\ell} = (Y_\ell P_\ell)^\top \frac{\partial y}{\partial Z_\ell}$, which only requires storing $Y_\ell P_\ell \in \mathbb{R}^{s_\ell \times r_\ell}$. Since $r_\ell \ll m_\ell$, this yields substantial memory savings for activations in such layers.

\subsubsection{Numerical Example}
We empirically validate the effectiveness of SO-Adam on both LLM pre-training and downstream fine-tuning tasks across various model architectures (LLaMA, OPT, RoBERTa). The proposed method is compared against the standard Adam, as well as state-of-the-art memory-efficient baselines including LoRA and GaLore.

\paragraph{Pre-training Performance}
We first evaluate SO-Adam on pre-training LLaMA architectures of varying sizes (from 60M up to 1B parameters) on the C4 dataset. Table \ref{tab:llama-pretrain} compares the validation perplexity and the estimated memory usage for optimizer states.

To push the limits of memory efficiency, we further pre-train the \textbf{LLaMA-7B} model on $8 \times$ A800 GPUs. As shown in Table \ref{tab:pretrain-llama-7b}, SO-Adam exhibits performance comparable to GaLore and Adam. Notably, the proposed method achieves this while reducing the total memory footprint by nearly 25~GB compared to Adam, and taking less than half the training time of GaLore and Adam.

\paragraph{Fine-Tuning on Downstream Tasks}
We extend the evaluation to fine-tuning tasks. Table \ref{tab:glue-finetuning} shows the results of fine-tuning a pre-trained RoBERTa-Base model on the GLUE benchmark across different subspace ranks. SO-Adam systematically surpasses LoRA and GaLore at both rank=4 and rank=8, and often achieves equal or better average scores than standard Adam.

To demonstrate generalizability to larger LLM families, we also fine-tune LLaMA (7B, 13B) and OPT (1.3B, 6.7B) models on the WinoGrande and COPA datasets. As shown in Table \ref{tab:finetune-wincopa}, the proposed method consistently outperforms other memory-efficient methods in terms of task accuracy.

\begin{table}[t]
    \centering
    \small
    \caption{\small Comparison of iteration time (in seconds) for different sequence lengths. LLaMA-1B runs on $4\times$ A800 GPUs, while LLaMA-7B uses $8\times$ A800 GPUs. For LLaMA-7B at Seq 256, Adam runs out of memory (OOM), whereas SO-Adam continues to operate efficiently.}\label{tab:time}
    \vspace{0.5em}
    \begin{threeparttable}
    \begin{tabular}{l ccc c ccc}
    \toprule
        \multirow{2}{*}{\textbf{Method}} & \multicolumn{3}{c}{\textbf{LLaMA-1B (Seconds)}} & \phantom{-} & \multicolumn{3}{c}{\textbf{LLaMA-7B (Seconds)}} \\
        \cmidrule{2-4} \cmidrule{6-8} 
        & \textbf{Seq 64} & \textbf{Seq 128} & \textbf{Seq 256} & & \textbf{Seq 64} & \textbf{Seq 128} & \textbf{Seq 256} \\
    \midrule
        \textbf{SO-Adam}    & \textbf{0.94}  & \textbf{1.70}  & \textbf{3.29}  & & \textbf{2.40}  & \textbf{2.94}  & \textbf{4.60}  \\
        GaLore & 1.12  & 1.84  & 3.35  & & 7.86  & 8.26  & 9.12  \\
        Adam   & 1.11  & 1.81  & 3.32  & & 7.84  & 8.23  & OOM  \\
    \bottomrule
    \end{tabular}
    \end{threeparttable}
\end{table}

\begin{table}[t]
\centering
\small
\caption{\small Experimental results on RoBERTa-large (350M). Results are reported as mean accuracy
(standard deviation). SO-ZO outperforms MeZO and MeZO-LoRA by a considerable margin and approaches FT performance.}
\label{appen-tab-roberta}
\vspace{0.5em}
\begin{threeparttable}
\resizebox{\textwidth}{!}{%
\begin{tabular}{@{}lcccccc@{}}
\toprule
Task & \textbf{SST-2} & \textbf{SST-5} & \textbf{SNLI} & \textbf{MNLI} & \textbf{RTE} & \textbf{TREC} \\ 
\cmidrule(lr){2-3} \cmidrule(lr){4-6} \cmidrule(lr){7-7}
Type & \multicolumn{2}{c}{\text{sentiment}} & \multicolumn{3}{c}{\text{natural language inference}} & topic \\ 
\midrule
Zero-shot & 79.0 & 35.5 & 50.2 & 48.8 & 51.4 & 32.0 \\ 
\midrule
\multicolumn{7}{l}{Gradient-free methods: $k = 16$} \\
\quad MeZO & 86.3 (8.0) & 40.8 (2.5) & 68.5 (4.2) & 56.7 (3.4) & 58.6 (9.5) & 62.4 (10.2) \\ 
\quad MeZO-LoRA & \textbf{88.4 (1.6)} & 38.9 (2.0) & 67.0 (3.3) & 56.8 (1.0) & 59.7 (4.5) & 37.0 (5.1) \\ 
\quad SO-ZO & 88.0 (5.8) & \textbf{41.1 (2.8)} & \textbf{73.4 (4.0)} & \textbf{61.6 (4.6)} & \textbf{61.2 (9.1)} & \textbf{77.9 (7.4)} \\ 
\midrule
\multicolumn{7}{l}{Gradient-based methods: $k = 16$} \\
\quad FT & 89.3 (5.3) & 44.0 (1.7) & 72.7 (5.7) & 63.4 (4.3) & 61.3 (5.3) & 83.7 (4.7) \\ 
\quad FT-LoRA & 92.7 (0.9) & 45.5 (1.4) & 71.6 (4.5) & 59.4 (4.6) & 61.4 (6.9) & 75.8 (7.8) \\ 
\midrule
\multicolumn{7}{l}{Gradient-free methods: $k = 512$} \\
\quad MeZO & 93.7 (0.4) & \textbf{53.9 (1.9)} & 84.8 (1.1) & 76.6 (0.8) & 76.8 (3.1) & 95.0 (0.4) \\ 
\quad MeZO-LoRA & 91.7 (0.2) & 45.1 (1.4) & 73.1 (1.1) & 65.5 (0.9) & 72.7 (0.8) & 50.8 (1.9) \\ 
\quad SO-ZO & \textbf{94.1 (0.7)} & 53.0 (0.4) & \textbf{85.4 (0.8)} & \textbf{80.4 (1.0)} & \textbf{79.7 (2.0)} & \textbf{95.5 (0.4)} \\
\midrule
\multicolumn{7}{l}{Gradient-based methods: $k = 512$} \\
\quad FT & 94.4 (0.6) & 55.7 (1.6) & 88.3 (0.8) & 84.8 (0.7) & 82.7 (1.1) & 97.2 (0.3) \\ 
\quad FT-LoRA & 91.9 (2.1) & 52.4 (1.2) & 84.8 (0.6) & 74.8 (3.4) & 81.2 (1.6) & 96.1 (0.6) \\ 
\bottomrule
\end{tabular}%
}
\end{threeparttable}
\end{table}

\paragraph{Memory and Communication Efficiency Breakdown}
Unlike GaLore or LoRA, which primarily reduce optimizer states, the subspace optimization method further reduces \textit{activation memory} by caching lower-dimensional inputs and \textit{communication overhead} by synchronizing low-dimensional gradients. 

Table \ref{tab:time} compares the time per iteration across different sequence lengths. As the model size increases or the sequence length decreases, the advantage of the proposed method becomes more pronounced. For LLaMA-7B with a sequence length of 64, SO-Adam requires only 2.40 seconds per iteration, making it nearly $3\times$ faster than GaLore and Adam. Furthermore, when the sequence length increases to 256, the standard Adam optimizer encounters an out-of-memory (OOM) error due to excessive activation caching, whereas SO-Adam remains feasible.

\paragraph{Fine-Tuning on Large Language Models}
To further demonstrate the scalability and effectiveness of our method in realistic, high-dimensional scenarios, we present the complete results for fine-tuning the RoBERTa-large (350M) model. As shown in Table \ref{appen-tab-roberta}, our SO-ZO method outperforms other gradient-free baselines on almost all datasets.

\FloatBarrier
\begingroup
\setlength{\abovedisplayskip}{4pt}
\setlength{\belowdisplayskip}{4pt}
\setlength{\abovedisplayshortskip}{2pt}
\setlength{\belowdisplayshortskip}{2pt}
\setlength{\jot}{2pt}

\section{Convergence of Matrix Randomized Subspace Optimization}
After surveying how matrix subspace optimization can be instantiated across
stochastic, zeroth-order, second-order, constrained, federated, decentralized,
bilevel, and large-model training settings, we now provide a representative
convergence result. Rather than analyzing every possible pairing of subspace
selection rules and downstream solvers, we focus on the random-subspace case,
which offers a clean theoretical prototype: it isolates the effect of randomized
dimensionality reduction, fits naturally with stochastic large-scale learning,
and already covers the nonconvex matrix-parameterized regime considered
throughout this primer. The proof presented below follows the convergence
analysis of \emph{A Memory Efficient Randomized Subspace Optimization Method for
Training Large Language Models}~\cite{chen2025memory}.

We consider a layerwise matrix variable $W=\{W_\ell\}_{\ell=1}^L$, where $W_\ell\in\mathbb{R}^{m_\ell\times n_\ell}$. At iteration $k$, Matrix Randomized Subspace Optimization samples $P_k=\{P_{k,\ell}\}_{\ell=1}^L$, with $P_{k,\ell}\in\mathbb{R}^{m_\ell\times r_\ell}$ and $r_\ell\ll m_\ell$. The reduced variable is $B=\{B_\ell\}_{\ell=1}^L$, where $B_\ell\in\mathbb{R}^{r_\ell\times n_\ell}$, and $P_kB$ denotes $\{P_{k,\ell}B_\ell\}_{\ell=1}^L$. The reduced subproblem and the lifted update are
\begin{equation}
    g_k(B):=f(W_k+P_kB)+\frac{1}{2\eta_k}\|B\|^2,
    \qquad
    W_{k+1}=W_k+P_k\widetilde B_k,
    \label{eq:matrix-randomized-subspace-subproblem-update}
\end{equation}
where $\widetilde B_k$ is an approximate solution of the reduced subproblem. Throughout this section, $\|\cdot\|$ denotes the layerwise Frobenius norm.

\begin{assumption}[Smoothness and lower boundedness]
The objective $f$ is lower bounded by $f_\star$ and has $L$-Lipschitz continuous gradient:\label{ass:convergence-smoothness}
\[
    \|\nabla f(U)-\nabla f(V)\|\le L\|U-V\|,
\]
for all layerwise matrix variables $U$ and $V$. No convexity of $f$ is required.
\end{assumption}

\begin{samepage}
\begin{assumption}[Random subspace normalization]
For each layer $\ell$ and iteration $k$, $P_{k,\ell}$ is sampled independently of $W_k$ and satisfies\label{ass:convergence-subspace-normalization}
\[
    \mathbb{E}\!\left[P_{k,\ell}P_{k,\ell}^{\top}\right]=I_{m_\ell},
    \qquad
    \|P_{k,\ell}\|_2^2\le \frac{m_\ell}{r_\ell}.
\]
Consequently, for any layerwise matrix $B$ and any fixed layerwise matrix $G$,
\[
    \|P_kB\|\le
    \sqrt{\max_\ell\frac{m_\ell}{r_\ell}}\|B\|,
    \qquad
    \mathbb{E}_{P_k}\!\left[\|P_k^\top G\|^2\right]=\|G\|^2.
\]
\end{assumption}
\end{samepage}

\begin{lemma}
Let $\widehat L:=L\max_\ell\frac{m_\ell}{r_\ell}$, and the notation will be used throughout the rest of this section. Under Assumptions~\ref{ass:convergence-smoothness} and~\ref{ass:convergence-subspace-normalization}, $g_k(B)$ is
$\left(\frac{1}{\eta_k}-\widehat L\right)$-strongly convex and
$\left(\frac{1}{\eta_k}+\widehat L\right)$-smooth whenever
$0<\eta_k<1/\widehat L$.
\end{lemma}

\begin{proof}
Let $h_k(B):=f(W_k+P_kB)$. By the chain rule,
$\nabla h_k(B)=P_k^\top\nabla f(W_k+P_kB)$. For any $B_1,B_2$, the smoothness of $f$ and the normalization of $P_k$ imply
\begin{align}
    \|\nabla h_k(B_1)-\nabla h_k(B_2)\|^2
    &=
    \sum_{\ell=1}^L
    \left\|
    P_{k,\ell}^\top
    \left(
    \nabla_\ell f(W_k+P_kB_1)
    -
    \nabla_\ell f(W_k+P_kB_2)
    \right)
    \right\|_F^2 \notag\\
    &\le
    \max_\ell\frac{m_\ell}{r_\ell}
    \left\|
    \nabla f(W_k+P_kB_1)
    -
    \nabla f(W_k+P_kB_2)
    \right\|^2 \notag\\
    &\le
    L^2\max_\ell\frac{m_\ell}{r_\ell}
    \|P_k(B_1-B_2)\|^2
    \le
    \widehat L^2\|B_1-B_2\|^2 .
    \label{eq:compact-hk-smooth}
\end{align}
Thus $h_k$ is $\widehat L$-smooth, and adding $\frac{1}{2\eta_k}\|B\|^2$ shows that $g_k$ is
$\left(\frac{1}{\eta_k}+\widehat L\right)$-smooth. Moreover, the smoothness inequality for $h_k$ gives
\begin{align}
    g_k(B_2)
    &\ge
    g_k(B_1)
    +
    \langle \nabla h_k(B_1),B_2-B_1\rangle
    -
    \frac{\widehat L}{2}\|B_2-B_1\|^2
    +
    \frac{1}{2\eta_k}\|B_2\|^2
    -
    \frac{1}{2\eta_k}\|B_1\|^2 \notag\\
    &=
    g_k(B_1)
    +
    \langle \nabla g_k(B_1),B_2-B_1\rangle
    +
    \left(\frac{1}{2\eta_k}-\frac{\widehat L}{2}\right)
    \|B_2-B_1\|^2 .
    \label{eq:compact-gk-strong-convex}
\end{align}
Therefore $g_k$ is
$\left(\frac{1}{\eta_k}-\widehat L\right)$-strongly convex whenever
$0<\eta_k<1/\widehat L$.
\end{proof}

\begin{theorem}
Under Assumptions~\ref{ass:convergence-smoothness} and~\ref{ass:convergence-subspace-normalization}, suppose that each reduced subproblem is solved from $B=0$ to an expected $\epsilon$-inexact solution $\widetilde B_k$, meaning
\[
    \mathbb{E}\!\left[g_k(\widetilde B_k)\right]\le g_k(B_k^\star)+\epsilon,
\]
where $B_k^\star$ is the exact minimizer of $g_k$. If $\eta_k=1/(2\widehat L)$, then the sequence generated by Matrix Randomized Subspace Optimization satisfies
\begin{equation}
    \frac{1}{K}
    \sum_{k=0}^{K-1}
    \mathbb{E}\!\left[\|\nabla f(W_k)\|^2\right]
    \le
    \frac{18\widehat L\bigl(f(W_0)-f_\star\bigr)}{K}
    +
    18\widehat L\epsilon .
    \label{eq:matrix-randomized-subspace-convergence}
\end{equation}
\end{theorem}

\begin{proof}
Denote $\mu:=1/\eta_k-\widehat L$. Since $g_k$ is $\mu$-strongly convex, for any $B$,
\[
    g_k(B_k^\star)\le g_k(B)-\frac{\mu}{2}\|B_k^\star-B\|^2.
\]
Taking $B=0$ and using the expected $\epsilon$-inexactness of $\widetilde B_k$, we obtain
\[
    \mathbb{E}\!\left[g_k(\widetilde B_k)\right]
    \le
    \mathbb{E}\!\left[g_k(0)\right]
    -
    \frac{\mu}{2}\mathbb{E}\!\left[\|B_k^\star\|^2\right]
    +
    \epsilon .
\]
Summing over $k=0,\ldots,K-1$ gives
\begin{align}
    \sum_{k=0}^{K-1}
    \frac{\mu}{2}
    \mathbb{E}\!\left[\|B_k^\star\|^2\right]
    &\le
    \sum_{k=0}^{K-1}
    \mathbb{E}\!\left[g_k(0)-g_k(\widetilde B_k)\right]
    +
    K\epsilon .
    \label{eq:compact-first-sum}
\end{align}
By $g_k(0)=f(W_k)$ and $g_k(\widetilde B_k)
    =
    f(W_k+P_k\widetilde B_k)
    +
    \frac{1}{2\eta_k}\|\widetilde B_k\|^2
    =
    f(W_{k+1})
    +
    \frac{1}{2\eta_k}\|\widetilde B_k\|^2,
$
the right-hand side of \eqref{eq:compact-first-sum} telescopes:
\begin{align}
    \sum_{k=0}^{K-1}
    \frac{\mu}{2}
    \mathbb{E}\!\left[\|B_k^\star\|^2\right]
    &\le
    f(W_0)-\mathbb{E}\!\left[f(W_K)\right]+K\epsilon
    \le
    f(W_0)-f_\star+K\epsilon .
    \label{eq:compact-bstar-bound}
\end{align}
Dividing by $K$ yields
\begin{equation}
    \frac{1}{K}
    \sum_{k=0}^{K-1}
    \frac{\mu}{2}
    \mathbb{E}\!\left[\|B_k^\star\|^2\right]
    \le
    \frac{f(W_0)-f_\star}{K}+\epsilon .
    \label{eq:compact-average-bstar-bound}
\end{equation}

It remains to connect $\|B_k^\star\|$ with the stationarity measure
$\|\nabla f(W_k)\|$. Since $g_k$ is
$\left(1/\eta_k+\widehat L\right)$-smooth and $\nabla g_k(B_k^\star)=0$, $\|\nabla g_k(0)\|^2
    =
    \|\nabla g_k(0)-\nabla g_k(B_k^\star)\|^2
    \le
    \left(\frac{1}{\eta_k}+\widehat L\right)^2
    \|B_k^\star\|^2 .$
Combining this inequality with \eqref{eq:compact-average-bstar-bound} gives
\begin{equation}
    \frac{1}{K}
    \sum_{k=0}^{K-1}
    \frac{\mu}{2}
    \left(\frac{1}{\eta_k}+\widehat L\right)^{-2}
    \mathbb{E}\!\left[\|\nabla g_k(0)\|^2\right]
    \le
    \frac{f(W_0)-f_\star}{K}+\epsilon .
    \label{eq:compact-grad-g-bound}
\end{equation}
Furthermore, $\nabla g_k(0)=P_k^\top\nabla f(W_k)$. Since $P_k$ is sampled independently of $W_k$, for any fixed $W_k$,
\begin{align}
    \mathbb{E}_{P_k}\!\left[\|\nabla g_k(0)\|^2\right]
    &=
    \sum_{\ell=1}^L
    \mathbb{E}_{P_{k,\ell}}
    \operatorname{Tr}
    \left(
    \nabla_\ell f(W_k)^\top
    P_{k,\ell}P_{k,\ell}^\top
    \nabla_\ell f(W_k)
    \right) \notag\\
    &=
    \sum_{\ell=1}^L
    \operatorname{Tr}
    \left(
    \mathbb{E}_{P_{k,\ell}}
    \left[P_{k,\ell}P_{k,\ell}^\top\right]
    \nabla_\ell f(W_k)\nabla_\ell f(W_k)^\top
    \right) \notag\\
    &=
    \sum_{\ell=1}^L
    \|\nabla_\ell f(W_k)\|_F^2
    =
    \|\nabla f(W_k)\|^2 .
    \label{eq:compact-random-projection-identity}
\end{align}
Substituting \eqref{eq:compact-random-projection-identity} into
\eqref{eq:compact-grad-g-bound}, and then setting
$\eta_k=1/(2\widehat L)$, we have $1/\eta_k=2\widehat L$ and
$\mu=\widehat L$. Hence\newline
$\frac{2}{\mu}
    \left(\frac{1}{\eta_k}+\widehat L\right)^2
    =
    18\widehat L,$
which gives
\[
    \frac{1}{K}
    \sum_{k=0}^{K-1}
    \mathbb{E}\!\left[\|\nabla f(W_k)\|^2\right]
    \le
    \frac{18\widehat L\bigl(f(W_0)-f_\star\bigr)}{K}
    +
    18\widehat L\epsilon .
\]
This proves the theorem.
\end{proof}

The theorem shows that Matrix Randomized Subspace Optimization reaches an average first-order stationary point up to an error floor determined by the accuracy of the reduced subproblem solver. In particular, if the reduced subproblems are solved increasingly accurately, the average squared gradient norm decays at the standard $O(1/K)$ rate. This result is consistent with the central philosophy of this primer: the convergence guarantee depends on the subspace distribution and the inexactness level of the downstream solver, rather than on the internal details of a particular optimizer. Thus, the analysis of Matrix Randomized Subspace Optimization provides a compact theoretical anchor for the broader matrix subspace optimization framework developed throughout the paper.

\endgroup

\section{Discussion and Future Work}
Subspace Matrix Optimization can be viewed as a meta-procedure that
\begin{enumerate}[label=(\roman*), leftmargin=*, itemsep=0pt, topsep=2pt]
  \item treats subspace selection as a modular component---allowing deterministic, randomized, or data-driven constructions of a low-rank basis $P$;
  \item fixes a specific reduction rule by reparameterizing the high-dimensional variable $W$ as $W+PB$ and optimizing over the low-dimensional coordinates $B$;
  \item constructs subproblems preserving enough of the original problem's structure (geometry, constraints, smooth/nonsmooth composition), in order to utilize the established optimization algorithms within the domain.
\end{enumerate}

Future work can be organized around four tightly coupled questions:
(a) how to build ``good'' subspaces cheaply,
(b) how to redesign the reduced subproblem so that domain solvers remain valid,
(c) how to inherit convergence theory and parallelism advantages from block coordinate descent (BCD),
and (d) when subspace updates are genuinely preferable to BCD.

\paragraph{Subspace design between random sampling and exact spectral selection}
Current subspace constructions span two extremes.
At one end, random Stiefel-type sampling is computationally light and easy to deploy, but it ignores gradient information entirely.
At the other end, choosing $P$ from a truncated SVD of $\nabla f(X)$ produces a subspace containing the most dominant first-order descent information under a rank budget, but it requires access to the full gradient and can erase the intended savings in both computation and memory.
A more compelling direction is therefore to seek not the ``best'' subspace at any cost, but the best subspace affordable under a prescribed resource budget.

One route is to replace exact spectral extraction by approximate spectral construction.
Rather than forming a full SVD of $\nabla f(X)$, randomized range finding can recover an approximate dominant range from compressed sketches, shifting the heavy linear algebra to a much smaller surrogate \cite{halko2011finding}.
A complementary route is gradient-informed sampling:
instead of drawing $P$ uniformly from the Stiefel manifold, one may sample directions, columns, or low-dimensional sketches using optimization-aware scores derived from partial gradient information, predicted decrease, constraint sensitivity, or inexpensive curvature surrogates.
This would retain much of the construction efficiency of random subspaces while biasing the selected basis toward informative descent directions.

Another route is to maintain $P$ incrementally rather than rebuilding it from a fresh full-gradient snapshot.
One may update the current basis through low-rank corrections inferred from iterate/gradient differences, together with lightweight re-orthogonalization to control numerical drift.
From a sketching viewpoint, Liberty's Frequent Directions maintains compact streaming summaries of gradient information with spectral guarantees \cite{liberty2013simple}, while later analysis strengthens its error bounds and mergeability across workers \cite{ghashami2016frequent}.
Viewed this way, subspace construction becomes a dynamic estimation problem:
the goal is to track the persistent low-dimensional geometry of useful updates while keeping the cost of maintaining $P$ well below that of exact SVD-based rules.

\paragraph{Redesigning reduced subproblems to preserve structure} Subspace Matrix Optimization has a crucial limitation: the additive substitution $W \mapsto W + PB$
can destroy the structure that enables specialized solvers.

On matrix manifolds, feasibility is typically violated; hence, future work should develop structure-preserving
subspace reductions built in tangent spaces with retractions (and, if needed, vector transports)
so that the reduced problem remains within the scope of Riemannian optimization toolkits
\cite{absil2008optimization}.
This reframes subspace selection as selecting a low-dimensional subspace of the tangent space---potentially
with geometry-aware metrics---rather than selecting directions in the ambient space.

For nonsmooth or composite problems, the main issue is preserving efficient proximal operators.
A promising direction is to design reduced models where the nonsmooth part remains ``prox-friendly,''
e.g., by constructing local surrogates whose proximal maps are computable
\cite{parikh2014proximal}.
More broadly, subproblem design can borrow from trust-region thinking:
instead of substituting $W+PB$ exactly, build a low-dimensional model with controllable approximation error
and accept/reject steps based on agreement tests
\cite{conn2000trust}.

\paragraph{Learning from BCD: transferring theory and enabling parallel subspace methods}
Since BCD is a special case in which $P$ selects coordinate blocks, it provides a natural reference point for both theory and systems design in subspace methods.

On the theoretical side, BCD analyses justify viewing subspace updates as blockwise steps whose convergence depends on regularity and block quality \cite{beck2013convergence}, with Dykstra--ADMM equivalences clarifying related coordinate-descent structures \cite{tibshirani2017dykstra}.
A key goal is to derive a unified rate framework that interpolates between coordinate descent, where $P$ is cheap but progress per step may be limited, and spectral subspace descent, where $P$ is more expensive but can capture more informative directions. Such a framework should make explicit how convergence depends on subspace dimension, approximation error, and construction cost.

Beyond rate templates, BCD also suggests several acceleration mechanisms that may transfer to subspace methods. 
The coordinate-descent literature highlights complementary mechanisms for selecting and updating informative coordinates, including nonuniform sampling \cite{richtarik2016optimal}, greedy Gauss--Southwell rules \cite{nutini2015coordinate}, accelerated parallel proximal updates \cite{fercoq2015accelerated}, and working-set strategies with extrapolation-based acceleration \cite{bertrand2022beyond}.
These ideas motivate natural subspace analogues: sampling candidate subspaces according to local smoothness or predicted decrease, greedily activating promising subspaces, maintaining active pools of informative subspaces, and combining reduced solves with momentum or extrapolation. A central question is which of these mechanisms remain valid after replacing coordinate blocks by low-rank subspaces, and how their gains depend on subspace overlap and approximation error.

On the systems side, BCD readily supports parallel and asynchronous implementations \cite{bertsekas2015parallel,calafiore2016parallel}.
Subspace methods may seek similar scalability by maintaining multiple low-dimensional subspaces $P=[P_1,\dots,P_m]$, solving the associated reduced subproblems in parallel, and aggregating their updates into a single step. The main challenge is to control interference among concurrent subspaces while remaining robust to stale iterates and outdated subspace information. This also motivates communication-aware subspace design, where both subspace representations and updates admit compact, mergeable summaries in distributed settings.

\paragraph{When low-rank structure truly favors subspace methods} A recurring empirical claim is that subspace methods shine when informative directions live in a low-dimensional spectral subspace.
Future work should formalize this via stable rank, singular value decay of gradient/curvature sketches,
and incoherence of optimal solutions, leading to problem-dependent complexity bounds that trade
subspace dimension, subspace approximation error, and inner-solve accuracy against total cost.
Such theory would yield actionable diagnostics:
estimate spectral decay from inexpensive sketches, select between BCD and subspace updates, and adapt the subspace rank online.

\section{Conclusion}
This primer constructs subspace optimization framework from a matrix-oriented perspective. Rather than flattening modern model parameters into a single vector, we treated matrices as the primitive optimization variables and used low-dimensional subspace to construct reduced problems that preserve the structure of the original optimization task. The central message is that subspace optimization should be viewed not as a single algorithm, but as a flexible interface: one may choose a random, greedy, power-iteration, or coordinate-aligned subspace, and then pair it with an appropriate downstream solver.

Across stochastic optimization, zeroth-order optimization, Newton-type methods, ADMM, PDHG, federated learning, decentralized optimization, bilevel optimization, and large-language-model training, the same principle repeatedly appears: optimizing over \(B\) in \(W+PB\) can substantially reduce memory, computation, and communication while still retaining the dominant update information of the full problem. The numerical examples further suggest that matrix subspace methods are especially effective when useful update directions are coupled across coordinates or concentrated in a low-rank spectral structure. In such cases, subspace updates often outperform coordinate-descent baselines under comparable update budgets and can approach the performance of full-rank methods at a fraction of the resource cost.

Several challenges remain open. Future work should develop cheaper yet informative subspace construction rules, design structure-preserving reduced subproblems for constrained, nonsmooth, and manifold settings, establish convergence theory that interpolates between coordinate descent and spectral subspace descent, and build practical diagnostics for deciding when subspace methods should be preferred. Overall, matrix subspace optimization provides a promising and extensible framework for memory-efficient large-scale optimization, particularly in modern machine-learning systems where parameter dimension, communication cost, and activation memory have become central bottlenecks.

\clearpage

\bibliographystyle{unsrt}
\bibliography{reference}

@inproceedings{chen2025memory,
  title={A Memory Efficient Randomized Subspace Optimization Method for Training Large Language Models},
  author={Chen, Yiming and Zhang, Yuan and Liu, Yin and Yuan, Kun and Wen, Zaiwen},
  booktitle={International Conference on Machine Learning (ICML)},
  pages={9747--9763},
  year={2025},
  organization={PMLR}
}

@inproceedings{nesterov1983method,
  title={A method for solving the convex programming problem with convergence rate {$O(1/k^2)$}},
  author={Nesterov, Yu},
  booktitle={Dokl. Akad. Nauk SSSR},
  volume={269},
  pages={543--547},
  year={1983}
}

@article{fuji2025theoretical,
  title={Theoretical analysis of the randomized subspace regularized {Newton} method for non-convex optimization},
  author={Fuji, Terunari and Poirion, Pierre-Louis and Takeda, Akiko},
  journal={Open Journal of Mathematical Optimization},
  volume={6},
  pages={1--35},
  year={2025}
}

@article{nozawa2023randomized,
  title={Randomized subspace gradient method for constrained optimization},
  author={Nozawa, Ryota and Poirion, Pierre-Louis and Takeda, Akiko},
  journal={arXiv preprint arXiv:2307.03335},
  year={2023}
}

@inproceedings{higuchi2024improving,
  title={Improving convergence guarantees of random subspace second-order algorithm for nonconvex optimization},
  author={Higuchi, Rei and Poirion, Pierre-Louis and Takeda, Akiko},
  booktitle={International Conference on Learning Representations},
  volume={2025},
  pages={3228--3270},
  year={2025}
}

@article{nozawa2025zeroth,
  title={Zeroth-order random subspace algorithm for non-smooth convex optimization},
  author={Nozawa, Ryota and Poirion, Pierre-Louis and Takeda, Akiko},
  journal={Journal of Optimization Theory and Applications},
  volume={204},
  number={3},
  pages={53},
  year={2025},
  publisher={Springer}
}

@article{chen2024sequential,
  title={Sequential subspace methods on {Stiefel} manifold optimization problems},
  author={Chen, Pengwen and Cheng, Chung-Kuan and Holtz, Chester},
  journal={arXiv preprint arXiv:2404.13301},
  year={2024}
}

@article{liu2021subspace,
  title={Subspace methods for nonlinear optimization},
  author={Liu, Xin and Wen, Zaiwen and Yuan, Ya-Xiang},
  journal={CSIAM Transactions on Applied Mathematics},
  volume={2},
  number={4},
  pages={585--651},
  year={2021}
}

@article{nesterov2017random,
  title={Random gradient-free minimization of convex functions},
  author={Nesterov, Yurii and Spokoiny, Vladimir},
  journal={Foundations of Computational Mathematics},
  volume={17},
  number={2},
  pages={527--566},
  year={2017},
  publisher={Springer}
}

@article{duchi2015optimal,
  title={Optimal rates for zero-order convex optimization: The power of two function evaluations},
  author={Duchi, John C and Jordan, Michael I and Wainwright, Martin J and Wibisono, Andre},
  journal={IEEE Transactions on Information Theory},
  volume={61},
  number={5},
  pages={2788--2806},
  year={2015},
  publisher={IEEE}
}

@article{malladi2023fine,
  title={Fine-tuning language models with just forward passes},
  author={Malladi, Sadhika and Gao, Tianyu and Nichani, Eshaan and Damian, Alex and Lee, Jason D and Chen, Danqi and Arora, Sanjeev},
  journal={Advances in Neural Information Processing Systems},
  volume={36},
  pages={53038--53075},
  year={2023}
}

@article{ghadimi2013stochastic,
  title={Stochastic first-and zeroth-order methods for nonconvex stochastic programming},
  author={Ghadimi, Saeed and Lan, Guanghui},
  journal={SIAM Journal on Optimization},
  volume={23},
  number={4},
  pages={2341--2368},
  year={2013},
  publisher={SIAM}
}

@article{bracken1973mathematical,
  title={Mathematical programs with optimization problems in the constraints},
  author={Bracken, Jerome and McGill, James T},
  journal={Operations Research},
  volume={21},
  number={1},
  pages={37--44},
  year={1973},
  publisher={INFORMS}
}

@article{sinha2017review,
  title={A review on bilevel optimization: From classical to evolutionary approaches and applications},
  author={Sinha, Ankur and Malo, Pekka and Deb, Kalyanmoy},
  journal={IEEE transactions on evolutionary computation},
  volume={22},
  number={2},
  pages={276--295},
  year={2017},
  publisher={IEEE}
}

@inproceedings{franceschi2018bilevel,
  title={Bilevel programming for hyperparameter optimization and meta-learning},
  author={Franceschi, Luca and Frasconi, Paolo and Salzo, Saverio and Grazzi, Riccardo and Pontil, Massimiliano},
  booktitle={International Conference on Machine Learning (ICML)},
  pages={1568--1577},
  year={2018},
  organization={PMLR}
}

@article{duchi2011adaptive,
  title={Adaptive subgradient methods for online learning and stochastic optimization.},
  author={Duchi, John and Hazan, Elad and Singer, Yoram},
  journal={Journal of Machine Learning Research},
  volume={12},
  number={7},
  year={2011}
}

@inproceedings{finn2017model,
  title={Model-agnostic meta-learning for fast adaptation of deep networks},
  author={Finn, Chelsea and Abbeel, Pieter and Levine, Sergey},
  booktitle={International Conference on Machine Learning (ICML)},
  pages={1126--1135},
  year={2017},
  organization={PMLR}
}

@article{liu2018darts,
  title={Darts: Differentiable architecture search},
  author={Liu, Hanxiao and Simonyan, Karen and Yang, Yiming},
  journal={arXiv preprint arXiv:1806.09055},
  year={2018}
}

@article{dagreou2022framework,
  title={A framework for bilevel optimization that enables stochastic and global variance reduction algorithms},
  author={Dagr{\'e}ou, Mathieu and Ablin, Pierre and Vaiter, Samuel and Moreau, Thomas},
  journal={Advances in Neural Information Processing Systems},
  volume={35},
  pages={26698--26710},
  year={2022}
}

@book{nocedal2006numerical,
  title={Numerical optimization},
  author={Nocedal, Jorge and Wright, Stephen J},
  publisher={Springer}
}

@book{dennis1996numerical,
  title={Numerical methods for unconstrained optimization and nonlinear equations},
  author={Dennis Jr, John E and Schnabel, Robert B},
  year={1996},
  publisher={SIAM}
}

@inproceedings{martens2010deep,
  title={Deep learning via hessian-free optimization.},
  author={Martens, James and others},
  booktitle={Icml},
  volume={27},
  pages={735--742},
  year={2010}
}

@inproceedings{xu2020second,
  title={Second-order optimization for non-convex machine learning: An empirical study},
  author={Xu, Peng and Roosta, Fred and Mahoney, Michael W},
  booktitle={Proceedings of the 2020 SIAM International Conference on Data Mining},
  pages={199--207},
  year={2020},
  organization={SIAM}
}

@inproceedings{mcmahan2017communication,
  title={Communication-efficient learning of deep networks from decentralized data},
  author={McMahan, Brendan and Moore, Eider and Ramage, Daniel and Hampson, Seth and y Arcas, Blaise Aguera},
  booktitle={International Conference on Artificial Intelligence and Statistics (AISTATS)},
  pages={1273--1282},
  year={2017},
  organization={PMLR}
}

@article{hinton2012divide,
  title={Rmsprop: Divide the gradient by a running average of its recent magnitude},
  author={Hinton, Geoffrey and Srivastava, Nitish and Swersky, Kevin},
  journal={Neural networks for machine learning, Coursera lecture 6e},
  volume={13},
  year={2012}
}

@article{glowinski1975approximation,
  title={Sur l'approximation, par {\'e}l{\'e}ments finis d'ordre un, et la r{\'e}solution, par p{\'e}nalisation-dualit{\'e} d'une classe de probl{\`e}mes de Dirichlet non lin{\'e}aires},
  author={Glowinski, Roland and Marroco, Americo},
  journal={Revue fran{\c{c}}aise d'automatique, informatique, recherche op{\'e}rationnelle. Analyse num{\'e}rique},
  volume={9},
  number={R2},
  pages={41--76},
  year={1975},
  publisher={EDP Sciences}
}

@article{gabay1976dual,
  title={A dual algorithm for the solution of nonlinear variational problems via finite element approximation},
  author={Gabay, Daniel and Mercier, Bertrand},
  journal={Computers \& Mathematics with Applications},
  volume={2},
  number={1},
  pages={17--40},
  year={1976},
  publisher={Elsevier}
}

@article{he20121,
  title={On the {O(1/n)} convergence rate of the {Douglas--Rachford} alternating direction method},
  author={He, Bingsheng and Yuan, Xiaoming},
  journal={SIAM Journal on Numerical Analysis},
  volume={50},
  number={2},
  pages={700--709},
  year={2012},
  publisher={SIAM}
}

@article{deng2016global,
  title={On the global and linear convergence of the generalized alternating direction method of multipliers},
  author={Deng, Wei and Yin, Wotao},
  journal={Journal of Scientific Computing},
  volume={66},
  number={3},
  pages={889--916},
  year={2016},
  publisher={Springer}
}

@article{boyd2011distributed,
  title={Distributed optimization and statistical learning via the alternating direction method of multipliers},
  author={Boyd, Stephen and Parikh, Neal and Chu, Eric and Peleato, Borja and Eckstein, Jonathan and others},
  journal={Foundations and Trends{\textregistered} in Machine Learning},
  volume={3},
  number={1},
  pages={1--122},
  year={2011},
  publisher={Now Publishers, Inc.}
}

@article{yuan2016convergence,
  title={On the convergence of decentralized gradient descent},
  author={Yuan, Kun and Ling, Qing and Yin, Wotao},
  journal={SIAM Journal on Optimization},
  volume={26},
  number={3},
  pages={1835--1854},
  year={2016},
  publisher={SIAM}
}

@article{shi2015extra,
  title={{EXTRA}: an exact first-order algorithm for decentralized consensus optimization},
  author={Shi, Wei and Ling, Qing and Wu, Gang and Yin, Wotao},
  journal={SIAM Journal on Optimization},
  volume={25},
  number={2},
  pages={944--966},
  year={2015},
  publisher={SIAM}
}

@article{nedic2017achieving,
  title={Achieving geometric convergence for distributed optimization over time-varying graphs},
  author={Nedic, Angelia and Olshevsky, Alex and Shi, Wei},
  journal={SIAM Journal on Optimization},
  volume={27},
  number={4},
  pages={2597--2633},
  year={2017},
  publisher={SIAM}
}

@article{nedic2009distributed,
  title={Distributed subgradient methods for multi-agent optimization},
  author={Nedic, Angelia and Ozdaglar, Asuman},
  journal={IEEE Transactions on Automatic Control},
  volume={54},
  number={1},
  pages={48--61},
  year={2009},
  publisher={IEEE}
}

@article{halko2011finding,
  title={Finding structure with randomness: Probabilistic algorithms for constructing approximate matrix decompositions},
  author={Halko, Nathan and Martinsson, Per-Gunnar and Tropp, Joel A},
  journal={SIAM Review},
  volume={53},
  number={2},
  pages={217--288},
  year={2011},
  publisher={SIAM}
}

@article{ghashami2016frequent,
  title={Frequent directions: Simple and deterministic matrix sketching},
  author={Ghashami, Mina and Liberty, Edo and Phillips, Jeff M and Woodruff, David P},
  journal={SIAM Journal on Computing},
  volume={45},
  number={5},
  pages={1762--1792},
  year={2016},
  publisher={SIAM}
}

@inproceedings{liberty2013simple,
  title={Simple and deterministic matrix sketching},
  author={Liberty, Edo},
  booktitle={Proceedings of the 19th ACM SIGKDD international conference on Knowledge discovery and data mining},
  pages={581--588},
  year={2013}
}

@book{absil2008optimization,
  title={Optimization algorithms on matrix manifolds},
  author={Absil, P-A and Mahony, Robert and Sepulchre, Rodolphe},
  year={2008},
  publisher={Princeton University Press}
}

@article{parikh2014proximal,
  title={Proximal algorithms},
  author={Parikh, Neal and Boyd, Stephen},
  journal={Foundations and Trends in Optimization},
  volume={1},
  number={3},
  pages={127--239},
  year={2014},
  publisher={Emerald Publishing Limited}
}

@book{conn2000trust,
  title={Trust region methods},
  author={Conn, Andrew R and Gould, Nicholas IM and Toint, Philippe L},
  year={2000},
  publisher={SIAM}
}

@article{beck2013convergence,
  title={On the convergence of block coordinate descent type methods},
  author={Beck, Amir and Tetruashvili, Luba},
  journal={SIAM Journal on Optimization},
  volume={23},
  number={4},
  pages={2037--2060},
  year={2013},
  publisher={SIAM}
}

@article{tibshirani2017dykstra,
  title={Dykstra's algorithm, admm, and coordinate descent: Connections, insights, and extensions},
  author={Tibshirani, Ryan J},
  journal={Advances in Neural Information Processing Systems},
  volume={30},
  year={2017}
}

@book{bertsekas2015parallel,
  title={Parallel and distributed computation: numerical methods},
  author={Bertsekas, Dimitri and Tsitsiklis, John},
  year={2015},
  publisher={Athena Scientific}
}

@article{zhao2024galore,
  title={{GaLore}: Memory-efficient {LLM} training by gradient low-rank projection},
  author={Zhao, Jiawei and Zhang, Zhenyu and Chen, Beidi and Wang, Zhangyang and Anandkumar, Anima and Tian, Yuandong},
  journal={arXiv preprint arXiv:2403.03507},
  year={2024}
}

@article{kingma2014adam,
  title={Adam: A method for stochastic optimization},
  author={Kingma, Diederik P and Ba, Jimmy},
  journal={arXiv preprint arXiv:1412.6980},
  year={2014}
}

@article{robbins1951stochastic,
  title={A stochastic approximation method},
  author={Robbins, Herbert and Monro, Sutton},
  journal={The Annals of Mathematical Statistics},
  pages={400--407},
  year={1951},
  publisher={JSTOR}
}

@article{kiefer1952stochastic,
  title={Stochastic estimation of the maximum of a regression function},
  author={Kiefer, Jack and Wolfowitz, Jacob},
  journal={The Annals of Mathematical Statistics},
  pages={462--466},
  year={1952},
  publisher={JSTOR}
}

@article{polyak1964some,
  title={Some methods of speeding up the convergence of iteration methods},
  author={Polyak, Boris T},
  journal={Ussr Computational Mathematics and Mathematical Physics},
  volume={4},
  number={5},
  pages={1--17},
  year={1964},
  publisher={Elsevier}
}

@article{johnson2013accelerating,
  title={Accelerating stochastic gradient descent using predictive variance reduction},
  author={Johnson, Rie and Zhang, Tong},
  journal={Advances in Neural Information Processing Systems},
  volume={26},
  year={2013}
}

@article{calafiore2016parallel,
  title={Parallel block coordinate minimization with application to group regularized regression},
  author={Calafiore, Giuseppe C},
  journal={Optimization and Engineering},
  volume={17},
  number={4},
  pages={941--964},
  year={2016},
  publisher={Springer}
}

@incollection{sayed2014diffusion,
  title={Diffusion adaptation over networks},
  author={Sayed, Ali H},
  booktitle={Academic Press Library in Signal Processing},
  volume={3},
  pages={323--453},
  year={2014},
  publisher={Elsevier}
}

@article{li2020federated,
  title={Federated optimization in heterogeneous networks},
  author={Li, Tian and Sahu, Anit Kumar and Zaheer, Manzil and Sanjabi, Maziar and Talwalkar, Ameet and Smith, Virginia},
  journal={Proceedings of Machine Learning and Systems (MLSys)},
  volume={2},
  pages={429--450},
  year={2020}
}

@inproceedings{karimireddy2020scaffold,
  title={{SCAFFOLD}: stochastic controlled averaging for federated learning},
  author={Karimireddy, Sai Praneeth and Kale, Satyen and Mohri, Mehryar and Reddi, Sashank and Stich, Sebastian and Suresh, Ananda Theertha},
  booktitle={International Conference on Machine Learning (ICML)},
  pages={5132--5143},
  year={2020},
  organization={PMLR}
}

@inproceedings{reisizadeh2020fedpaq,
  title={{FedPAQ}: a communication-efficient federated learning method with periodic averaging and quantization},
  author={Reisizadeh, Amirhossein and Mokhtari, Aryan and Hassani, Hamed and Jadbabaie, Ali and Pedarsani, Ramtin},
  booktitle={International Conference on Artificial Intelligence and Statistics (AISTATS)},
  pages={2021--2031},
  year={2020},
  organization={PMLR}
}

@inproceedings{lin2017deep,
  title={Deep gradient compression: reducing the communication bandwidth for distributed training},
  author={Lin, Yujun and Han, Song and Mao, Huizi and Wang, Yu and Dally, William J},
  booktitle={International Conference on Learning Representations (ICLR)},
  year={2018},
  url={https://openreview.net/forum?id=SkhQHMW0W}
}

@article{t2020personalized,
  title={Personalized federated learning with moreau envelopes},
  author={T Dinh, Canh and Tran, Nguyen and Nguyen, Josh},
  journal={Advances in neural information processing systems},
  volume={33},
  pages={21394--21405},
  year={2020}
}

@article{arivazhagan2019federated,
  title={Federated learning with personalization layers},
  author={Arivazhagan, Manoj Ghuhan and Aggarwal, Vinay and Singh, Aaditya Kumar and Choudhary, Sunav},
  journal={arXiv preprint arXiv:1912.00818},
  year={2019}
}

@inproceedings{reddi2020adaptive,
  title={Adaptive federated optimization},
  author={Reddi, Sashank and Charles, Zachary and Zaheer, Manzil and Garrett, Zachary and Rush, Keith and Kone{\v{c}}n{\'y}, Jakub and Kumar, Sanjiv and McMahan, H Brendan},
  booktitle={International Conference on Learning Representations (ICLR)},
  year={2021},
  url={https://openreview.net/forum?id=LkFG3lB13U5}
}

@article{lian2017can,
  title={Can decentralized algorithms outperform centralized algorithms? a case study for decentralized parallel stochastic gradient descent},
  author={Lian, Xiangru and Zhang, Ce and Zhang, Huan and Hsieh, Cho-Jui and Zhang, Wei and Liu, Ji},
  journal={Advances in Neural Information Processing Systems (NeurIPS)},
  volume={30},
  year={2017}
}

@article{richtarik2016optimal,
  title={On optimal probabilities in stochastic coordinate descent methods},
  author={Richt{\'a}rik, Peter and Tak{\'a}{\v{c}}, Martin},
  journal={Optimization Letters},
  volume={10},
  number={6},
  pages={1233--1243},
  year={2016},
  publisher={Springer}
}

@inproceedings{nutini2015coordinate,
  title={Coordinate descent converges faster with the gauss-southwell rule than random selection},
  author={Nutini, Julie and Schmidt, Mark and Laradji, Issam and Friedlander, Michael and Koepke, Hoyt},
  booktitle={International Conference on Machine Learning},
  pages={1632--1641},
  year={2015},
  organization={PMLR}
}

@article{fercoq2015accelerated,
  title={Accelerated, parallel, and proximal coordinate descent},
  author={Fercoq, Olivier and Richt{\'a}rik, Peter},
  journal={SIAM Journal on Optimization},
  volume={25},
  number={4},
  pages={1997--2023},
  year={2015},
  publisher={SIAM}
}

@article{bertrand2022beyond,
  title={Beyond l1: Faster and better sparse models with skglm},
  author={Bertrand, Quentin and Klopfenstein, Quentin and Bannier, Pierre-Antoine and Gidel, Gauthier and Massias, Mathurin},
  journal={Advances in Neural Information Processing Systems},
  volume={35},
  pages={38950--38965},
  year={2022}
}

@article{arora2019implicit,
  title={Implicit Regularization in Deep Matrix Factorization},
  author={Arora, Sanjeev and Cohen, Nadav and Hu, Wei and Luo, Yuping},
  journal={Advances in Neural Information Processing Systems},
  volume={32},
  year={2019}
}

@article{huh2021lowrank,
  title={The Low-Rank Simplicity Bias in Deep Networks},
  author={Huh, Minyoung and Mobahi, Hossein and Zhang, Richard and Cheung, Brian and Agrawal, Pulkit and Isola, Phillip},
  journal={Transactions on Machine Learning Research},
  year={2023}
}

@article{ba2022highdim,
  title={High-dimensional Asymptotics of Feature Learning: How One Gradient Step Improves the Representation},
  author={Ba, Jimmy and Erdogdu, Murat A. and Suzuki, Taiji and Wang, Zhichao and Wu, Denny and Yang, Greg},
  journal={Advances in Neural Information Processing Systems},
  volume={35},
  year={2022}
}

@article{sonthalia2025low,
  title={Low Rank Gradients and Where to Find Them},
  author={Sonthalia, Rishabh and Radhakrishnan, Adityanarayanan and Moniri, Behrad and Hassani, Hamed and Dobriban, Edgar},
  journal={Advances in Neural Information Processing Systems},
  year={2025}
}

@inproceedings{jaiswal2025from,
  title={From Low Rank Gradient Subspace Stabilization to Low-Rank Weights: Observations, Theories, and Applications},
  author={Jaiswal, Ajay Kumar and Wang, Yifan and Yin, Lu and Liu, Shiwei and Chen, Runjin and Zhao, Jiawei and Grama, Ananth and Tian, Yuandong and Wang, Zhangyang},
  booktitle={Proceedings of the 42nd International Conference on Machine Learning},
  pages={26740--26756},
  year={2025},
  organization={PMLR}
}

@article{dai2011geometric,
  title={A Geometric Approach to Low-Rank Matrix Completion},
  author={Dai, Wei and Kerman, Ely and Milenkovic, Olgica},
  journal={IEEE Transactions on Information Theory},
  volume={58},
  number={1},
  pages={237--247},
  year={2012},
  publisher={IEEE}
}

@article{keshavan2009gradient,
  title={A Gradient Descent Algorithm on the Grassman Manifold for Matrix Completion},
  author={Keshavan, Raghunandan H. and Oh, Sewoong},
  journal={CoRR},
  volume={abs/0910.5260},
  year={2009}
}

@article{boumal2011rtrmc,
  title={RTRMC: A Riemannian Trust-Region Method for Low-Rank Matrix Completion},
  author={Boumal, Nicolas and Absil, Pierre-Antoine},
  journal={Advances in Neural Information Processing Systems},
  volume={24},
  pages={406--414},
  year={2011}
}

@article{vandereycken2013low,
  title={Low-Rank Matrix Completion by Riemannian Optimization},
  author={Vandereycken, Bart},
  journal={SIAM Journal on Optimization},
  volume={23},
  number={2},
  pages={1214--1236},
  year={2013},
  publisher={SIAM}
}

@article{meyer2011regression,
  title={Regression on Fixed-Rank Positive Semidefinite Matrices: A Riemannian Approach},
  author={Meyer, Gilles and Bonnabel, Silv{\`e}re and Sepulchre, Rodolphe},
  journal={Journal of Machine Learning Research},
  volume={12},
  pages={593--625},
  year={2011}
}

@inproceedings{li2020efficient,
  title={Efficient Riemannian Optimization on the Stiefel Manifold via the Cayley Transform},
  author={Li, Jun and Li, Fuxin and Todorovic, Sinisa},
  booktitle={International Conference on Learning Representations},
  year={2020}
}

@inproceedings{becigneul2019riemannian,
  title={Riemannian Adaptive Optimization Methods},
  author={B{\'e}cigneul, Gary and Ganea, Octavian-Eugen},
  booktitle={International Conference on Learning Representations},
  year={2019}
}

@inproceedings{han2024riemannian,
  title={Riemannian Coordinate Descent Algorithms on Matrix Manifolds},
  author={Han, Andi and Jawanpuria, Pratik and Mishra, Bamdev},
  booktitle={Proceedings of the 41st International Conference on Machine Learning},
  pages={17393--17415},
  year={2024},
  organization={PMLR}
}

@inproceedings{martens2015optimizing,
  title={Optimizing Neural Networks with Kronecker-factored Approximate Curvature},
  author={Martens, James and Grosse, Roger},
  booktitle={Proceedings of the 32nd International Conference on Machine Learning},
  pages={2408--2417},
  year={2015},
  organization={PMLR}
}

@inproceedings{gupta2018shampoo,
  title={Shampoo: Preconditioned Stochastic Tensor Optimization},
  author={Gupta, Vineet and Koren, Tomer and Singer, Yoram},
  booktitle={Proceedings of the 35th International Conference on Machine Learning},
  pages={1842--1850},
  year={2018},
  organization={PMLR}
}

@inproceedings{shazeer2018adafactor,
  title={Adafactor: Adaptive Learning Rates with Sublinear Memory Cost},
  author={Shazeer, Noam and Stern, Mitchell},
  booktitle={Proceedings of the 35th International Conference on Machine Learning},
  pages={4596--4604},
  year={2018},
  organization={PMLR}
}

@misc{jordan2024muon,
  title={Muon: An Optimizer for Hidden Layers in Neural Networks},
  author={Jordan, Keller and Jin, Yuchen and Boza, Vlado and You, Jiacheng and Cesista, Franz and Newhouse, Laker and Bernstein, Jeremy},
  year={2024}
}

@article{cai2010singular,
  title={A Singular Value Thresholding Algorithm for Matrix Completion},
  author={Cai, Jian-Feng and Cand{\`e}s, Emmanuel J. and Shen, Zuowei},
  journal={SIAM Journal on Optimization},
  volume={20},
  number={4},
  pages={1956--1982},
  year={2010},
  publisher={SIAM}
}

@article{jain2010guaranteed,
  title={Guaranteed Rank Minimization via Singular Value Projection},
  author={Jain, Prateek and Meka, Raghu and Dhillon, Inderjit S.},
  journal={Advances in Neural Information Processing Systems},
  volume={23},
  year={2010}
}

@article{hu2022lora,
  title={Lora: Low-rank adaptation of large language models},
  author={Hu, Edward J and Shen, Yelong and Wallis, Phillip and Allen-Zhu, Zeyuan and Li, Yuanzhi and Wang, Shean and Wang, Lu and Chen, Weizhu},
  journal={arXiv preprint arXiv:2106.09685},
  year={2021}
}

@article{han2024sltrain,
  title={{SLTrain}: A sparse plus low rank approach for parameter and memory efficient pretraining},
  author={Han, Andi and Li, Jiaxiang and Huang, Wei and Hong, Mingyi and Takeda, Akiko and Jawanpuria, Pratik and Mishra, Bamdev},
  journal={Advances in Neural Information Processing Systems},
  volume={37},
  pages={118267--118295},
  year={2024}
}

@article{lialin2023relora,
  title={{ReLoRA}: High-rank training through low-rank updates},
  author={Lialin, Vladislav and Shivagunde, Namrata and Muckatira, Sherin and Rumshisky, Anna},
  journal={arXiv preprint arXiv:2307.05695},
  year={2023}
}

@inproceedings{mo2025parameter,
  title={Parameter and memory efficient pretraining via low-rank riemannian optimization},
  author={Mo, Zhanfeng and Huang, Long-Kai and Pan, Sinno Jialin},
  booktitle={The Thirteenth International Conference on Learning Representations(ICLR)},
  year={2025}
}

@article{chen2024fira,
  title={Fira: Can we achieve full-rank training of {LLMs} under low-rank constraint?},
  author={Chen, Xi and Feng, Kaituo and Li, Changsheng and Lai, Xunhao and Yue, Xiangyu and Yuan, Ye and Wang, Guoren},
  journal={arXiv preprint arXiv:2410.01623},
  year={2024}
}

@article{zmushko2024frugal,
  title={{FRUGAL}: Memory-efficient optimization by reducing state overhead for scalable training},
  author={Zmushko, Philip and Beznosikov, Aleksandr and Tak{\'a}{\v{c}}, Martin and Horv{\'a}th, Samuel},
  journal={arXiv preprint arXiv:2411.07837},
  year={2024}
}

@inproceedings{xiao2025coap,
  title={{COAP}: Memory-efficient training with correlation-aware gradient projection},
  author={Xiao, Jinqi and Sang, Shen and Zhi, Tiancheng and Liu, Jing and Yan, Qing and Luo, Linjie and Yuan, Bo},
  booktitle={Proceedings of the IEEE/CVF Conference on Computer Vision and Pattern Recognition(CVPR)},
  pages={30116--30126},
  year={2025}
}

@article{he2024subspace,
  title={Subspace optimization for large language models with convergence guarantees},
  author={He, Yutong and Li, Pengrui and Hu, Yipeng and Chen, Chuyan and Yuan, Kun},
  journal={arXiv preprint arXiv:2410.11289},
  year={2024}
}

@article{luo2024badam,
  title={Badam: A memory efficient full parameter optimization method for large language models},
  author={Luo, Qijun and Yu, Hengxu and Li, Xiao},
  journal={Advances in Neural Information Processing Systems},
  volume={37},
  pages={24926--24958},
  year={2024}
}

@article{luoaccelerating,
  title={Accelerating Block Coordinate Descent for LLM Finetuning via Landscape Expansion},
  author={Luo, Qijun and Shen, Yifei and Peng, Liangzu and Li, Dongsheng and Li, Xiao},
  journal={Advances in Neural Information Processing Systems},
  volume={38},
  pages={56619--56645},
  year={2026}
}

@article{ramesh2024blockllm,
  title={Blockllm: Memory-efficient adaptation of llms by selecting and optimizing the right coordinate blocks},
  author={Ramesh, Amrutha Varshini and Ganapathiraman, Vignesh and Laradji, Issam H and Schmidt, Mark},
  journal={arXiv preprint arXiv:2406.17296},
  year={2024}
}

@inproceedings{liu2025cola,
  title={{CoLA}: Compute-efficient pre-training of {LLMs} via low-rank activation},
  author={Liu, Ziyue and Zhang, Ruijie and Wang, Zhengyang and Yan, Mingsong and Yang, Zi and Hovland, Paul D and Nicolae, Bogdan and Cappello, Franck and Tang, Sui and Zhang, Zheng},
  booktitle={Proceedings of the 2025 Conference on Empirical Methods in Natural Language Processing(EMNLP)},
  pages={4627--4645},
  year={2025}
}

@inproceedings{shamshoum2025compact,
  title={{CompAct}: Compressed activations for memory-efficient {LLM} training},
  author={Shamshoum, Yara and Hodos, Nitzan and Sieradzki, Yuval and Schuster, Assaf},
  booktitle={Proceedings of the 2025 Conference of the Nations of the Americas Chapter of the Association for Computational Linguistics: Human Language Technologies (Volume 1: Long Papers)(NAACL)},
  pages={1511--1524},
  year={2025}
}

@article{miles2024velora,
  title={{VeLoRA}: Memory efficient training using rank-1 sub-token projections},
  author={Miles, Roy and Reddy, Pradyumna and Elezi, Ismail and Deng, Jiankang},
  journal={Advances in Neural Information Processing Systems},
  volume={37},
  pages={42292--42310},
  year={2024}
}

@article{kong2025cr,
  title={{CR-Net}: Scaling Parameter-Efficient Training with Cross-Layer Low-Rank Structure},
  author={Kong, Boao and Liang, Junzhu and Liu, Yuxi and Deng, Renjia and Yuan, Kun},
  journal={The Fourteenth International Conference on Learning Representations (ICLR)},
  year={2026}
}

@inproceedings{wei2025flexora,
  title={Flexora: Flexible low-rank adaptation for large language models},
  author={Wei, Chenxing and Shu, Yao and He, Ying Tiffany and Yu, Fei},
  booktitle={Proceedings of the 63rd Annual Meeting of the Association for Computational Linguistics (Volume 1: Long Papers)(ACL)},
  pages={14643--14682},
  year={2025}
}

@inproceedings{zangrando2025debora,
  title={{dEBORA}: Efficient bilevel optimization-based low-rank adaptation},
  author={Zangrando, Emanuele and Venturini, Sara and Rinaldi, Francesco and Tudisco, Francesco},
  booktitle={The Thirteenth International Conference on Learning Representations(ICLR)},
  year={2025}
}

@article{rajabi2025subtrack++,
  title={{SubTrack++}: Gradient subspace tracking for scalable {LLM} training},
  author={Rajabi, Sahar and Nonta, Nayeema and Rambhatla, Sirisha},
  journal={arXiv preprint arXiv:2502.01586},
  year={2025}
}

@article{zhang2025i3s,
  title={Breaking the frozen subspace: Importance sampling for low-rank optimization in LLM pretraining},
  author={Zhang, Haochen and Yin, Junze and Wang, Guanchu and Liu, Zirui and Yang, Lin and Zhang, Tianyi and Shrivastava, Anshumali and Braverman, Vladimir},
  journal={Advances in Neural Information Processing Systems},
  volume={38},
  pages={10122--10174},
  year={2026}
}

@article{yuan1995subspace,
  title={A subspace study on conjugate gradient algorithms},
  author={Yuan, Y-X and Stoer, J},
  journal={ZAMM-Journal of Applied Mathematics and Mechanics/Zeitschrift f{\"u}r Angewandte Mathematik und Mechanik},
  volume={75},
  number={1},
  pages={69--77},
  year={1995},
  publisher={Wiley Online Library}
}

@article{karimi2012conjugate,
  title={Conjugate gradient with subspace optimization (2012)},
  author={Karimi, S and Vavasis, S},
  journal={arXiv preprint arXiv:1202.1479}
}

@book{narkiss2005sesop,
  title={Sequential subspace optimization method for large-scale unconstrained problems},
  author={Narkiss, Guy and Zibulevsky, Michael},
  year={2005},
  publisher={Technion-IIT, Department of Electrical Engineering}
}

@article{zibulevsky2008sesop,
  title={{SESOP-TN}: Combining sequential subspace optimization with truncated {Newton} method},
  author={Zibulevsky, Michael},
  journal={Israel Institute of Technology, Technical Report},
  year={2008}
}

@article{pulay1980convergence,
  title={Convergence acceleration of iterative sequences. {The} case of {SCF} iteration},
  author={Pulay, P{\'e}ter},
  journal={Chemical Physics Letters},
  volume={73},
  number={2},
  pages={393--398},
  year={1980},
  publisher={Elsevier}
}

@article{csaszar1984geometry,
  title={Geometry optimization by direct inversion in the iterative subspace},
  author={Cs{\'a}sz{\'a}r, P{\'a}l and Pulay, P{\'e}ter},
  journal={Journal of Molecular Structure},
  volume={114},
  pages={31--34},
  year={1984},
  publisher={Elsevier}
}

@article{chen2019convergence,
  title={Convergence of the {EDIIS} algorithm for nonlinear equations},
  author={Chen, Xiaojun and Kelley, Carl T},
  journal={SIAM Journal on Scientific Computing},
  volume={41},
  number={1},
  pages={A365--A379},
  year={2019},
  publisher={SIAM}
}

@article{walker2011anderson,
  title={Anderson acceleration for fixed-point iterations},
  author={Walker, Homer F and Ni, Peng},
  journal={SIAM Journal on Numerical Analysis},
  volume={49},
  number={4},
  pages={1715--1735},
  year={2011},
  publisher={SIAM}
}

@article{conn1996iterated,
  title={On iterated-subspace minimization methods for nonlinear optimization},
  author={Conn, AR and Gould, Nick and Sartenaer, A and Toint, Ph L},
  journal={Linear and Nonlinear Conjugate Gradient-Related Methods},
  pages={50--78},
  year={1996}
}

@article{gower2015randomized,
  title={Randomized iterative methods for linear systems},
  author={Gower, Robert M and Richt{\'a}rik, Peter},
  journal={SIAM Journal on Matrix Analysis and Applications},
  volume={36},
  number={4},
  pages={1660--1690},
  year={2015},
  publisher={SIAM}
}

@article{gower2019adaptive,
  title={Adaptive sketch-and-project methods for solving linear systems},
  author={Gower, Robert and Molitor, Denali and Moorman, Jacob and Needell, Deanna},
  journal={arXiv preprint arXiv:1909.03604},
  year={2019}
}

@article{necoara2021randomized,
  title={Randomized sketch descent methods for non-separable linearly constrained optimization},
  author={Necoara, Ion and Tak{\'a}{\v{c}}, Martin},
  journal={IMA Journal of Numerical Analysis},
  volume={41},
  number={2},
  pages={1056--1092},
  year={2021},
  publisher={Oxford University Press}
}

@article{byrd1988approximate,
  title={Approximate solution of the trust region problem by minimization over two-dimensional subspaces},
  author={Byrd, Richard H and Schnabel, Robert B and Shultz, Gerald A},
  journal={Mathematical Programming},
  volume={40},
  number={1},
  pages={247--263},
  year={1988},
  publisher={Springer}
}

@article{wang2004subspace,
  title={A subspace trust region method for large scale unconstrained optimization},
  author={Wang, ZH and Wen, ZW and Yuan, Y},
  journal={Numerical linear algebra and optimization},
  pages={264--274},
  year={2004},
  publisher={Science Press, Beijing}
}

@article{wang2006subspace,
  title={A subspace implementation of {quasi-Newton} trust region methods for unconstrained optimization},
  author={Wang, Zhou-Hong and Yuan, Ya-Xiang},
  journal={Numerische Mathematik},
  volume={104},
  number={2},
  pages={241--269},
  year={2006},
  publisher={Springer}
}

@article{liu1989limited,
  title={On the limited memory {BFGS} method for large scale optimization},
  author={Liu, Dong C and Nocedal, Jorge},
  journal={Mathematical Programming},
  volume={45},
  number={1},
  pages={503--528},
  year={1989},
  publisher={Springer}
}

@article{gill2001reduced,
  title={Reduced-{Hessian quasi-Newton} methods for unconstrained optimization},
  author={Gill, Philip E and Leonard, Michael W},
  journal={SIAM Journal on Optimization},
  volume={12},
  number={1},
  pages={209--237},
  year={2001},
  publisher={SIAM}
}

@article{gill2003limited,
  title={Limited-memory {reduced-Hessian} methods for large-scale unconstrained optimization},
  author={Gill, Philip E and Leonard, Michael W},
  journal={SIAM Journal on Optimization},
  volume={14},
  number={2},
  pages={380--401},
  year={2003},
  publisher={SIAM}
}

@article{burdakov2017efficiently,
  title={On efficiently combining limited-memory and trust-region techniques},
  author={Burdakov, Oleg and Gong, Lujin and Zikrin, Spartak and Yuan, Ya-xiang},
  journal={Mathematical Programming Computation},
  volume={9},
  number={1},
  pages={101--134},
  year={2017},
  publisher={Springer}
}

@article{branch1999subspace,
  title={A subspace, interior, and conjugate gradient method for large-scale bound-constrained minimization problems},
  author={Branch, Mary Ann and Coleman, Thomas F and Li, Yuying},
  journal={SIAM Journal on Scientific Computing},
  volume={21},
  number={1},
  pages={1--23},
  year={1999},
  publisher={SIAM}
}

@article{byrd2011active,
  title={An active-set algorithm for nonlinear programming using parametric linear programming},
  author={Byrd, Richard H and Waltz, Richard A},
  journal={Optimization Methods \& Software},
  volume={26},
  number={1},
  pages={47--66},
  year={2011},
  publisher={Taylor \& Francis}
}

@inproceedings{chen2024enhancing,
  title={Enhancing zeroth-order fine-tuning for language models with low-rank structures},
  author={Chen, Yiming and Zhang, Yuan and Cao, Liyuan and Yuan, Kun and Wen, Zaiwen},
  booktitle={The Thirteenth International Conference on Learning Representations (ICLR)},
  year={2025}
}

@inproceedings{yu2025zeroth,
  title={Zeroth-order fine-tuning of {LLMs} in random subspaces},
  author={Yu, Ziming and Zhou, Pan and Wang, Sike and Li, Jia and Tian, Mi and Huang, Hua},
  booktitle={Proceedings of the IEEE/CVF International Conference on Computer Vision (ICCV)},
  pages={4475--4485},
  year={2025}
}

@article{grishchenko2021proximal,
  title={Proximal gradient methods with adaptive subspace sampling},
  author={Grishchenko, Dmitry and Iutzeler, Franck and Malick, J{\'e}r{\^o}me},
  journal={Mathematics of Operations Research},
  volume={46},
  number={4},
  pages={1303--1323},
  year={2021},
  publisher={INFORMS}
}

@article{gower2019rsn,
  title={RSN: randomized subspace {Newton}},
  author={Gower, Robert and Kovalev, Dmitry and Lieder, Felix and Richt{\'a}rik, Peter},
  journal={Advances in Neural Information Processing Systems (NeurIPS)},
  volume={32},
  year={2019}
}

@article{cartis2023scalable,
  title={Scalable subspace methods for derivative-free nonlinear least-squares optimization},
  author={Cartis, Coralia and Roberts, Lindon},
  journal={Mathematical Programming},
  volume={199},
  number={1},
  pages={461--524},
  year={2023},
  publisher={Springer}
}

@article{chen2025clarsta,
  title={{CLARSTA}: a random subspace trust-region algorithm for convex-constrained derivative-free optimization},
  author={Chen, Yiwen and Hare, Warren and Wiebe, Amy},
  journal={arXiv preprint arXiv:2506.20335},
  year={2025}
}

@article{dzahini2024stochastic,
  title={Stochastic trust-region algorithm in random subspaces with convergence and expected complexity analyses},
  author={Dzahini, Kwassi Joseph and Wild, Stefan M},
  journal={SIAM Journal on Optimization},
  volume={34},
  number={3},
  pages={2671--2699},
  year={2024},
  publisher={SIAM}
}

@article{dzahini2024direct,
  title={Direct search for stochastic optimization in random subspaces with zeroth-, first-, and second-order convergence and expected complexity},
  author={Dzahini, KJ and Wild, SM},
  journal={arXiv preprint arXiv:2403.13320},
  year={2024}
}

@article{cartis2026randomized,
  title={Randomized subspace derivative-free optimization with quadratic models and second-order convergence},
  author={Cartis, Coralia and Roberts, Lindon},
  journal={Optimization Methods and Software},
  pages={1--28},
  year={2026},
  publisher={Taylor \& Francis}
}

@article{shang2024accelerated,
  title={Accelerated double-sketching subspace {Newton}},
  author={Shang, Jun and Ye, Haishan and Chang, Xiangyu},
  journal={European Journal of Operational Research},
  volume={319},
  number={2},
  pages={484--493},
  year={2024},
  publisher={Elsevier}
}

@article{cartis2025random,
  title={Random subspace cubic-regularization methods, with applications to low-rank functions},
  author={Cartis, Coralia and Shao, Zhen and Tansley, Edward},
  journal={arXiv preprint arXiv:2501.09734},
  year={2025}
}

@article{yang2026fedkrso,
  title={{FedKRSO}: communication and memory efficient federated fine-tuning of large language models},
  author={Yang, Guohao and Wu, Tongle and Guo, Yuanxiong and Sun, Ying and Gong, Yanmin},
  journal={arXiv preprint arXiv:2602.03019},
  year={2026}
}

@article{obara2022sequential,
  title={Sequential quadratic optimization for nonlinear optimization problems on {Riemannian} manifolds},
  author={Obara, Mitsuaki and Okuno, Takayuki and Takeda, Akiko},
  journal={SIAM Journal on Optimization},
  volume={32},
  number={2},
  pages={822--853},
  year={2022},
  publisher={SIAM}
}

@article{gutman2026randomized,
  title={Randomized Subspace Methods for Optimization on Fixed-Rank Matrix Manifolds},
  author={Gutman, David Huckleberry and Ho-Nguyen, Nam}
}

@inproceedings{jiang2024krylov,
  title={Krylov cubic regularized {Newton}: a subspace second-order method with dimension-free convergence rate},
  author={Jiang, Ruichen and Raman, Parameswaran and Sabach, Shoham and Mokhtari, Aryan and Hong, Mingyi and Cevher, Volkan},
  booktitle={International Conference on Artificial Intelligence and Statistics (AISTATS)},
  pages={4411--4419},
  year={2024},
  organization={PMLR}
}

@article{al2025extended,
  title={Extended-Krylov-subspace methods for trust-region and norm-regularization subproblems},
  author={Al Daas, Hussam and Gould, Nicholas IM},
  journal={SIAM Journal on Optimization},
  volume={36},
  number={4},
  pages={2060--2091},
  year={2026},
  publisher={SIAM}
}

@inproceedings{yang2024lancbio,
  title={{LancBiO}: dynamic {Lanczos-aided} bilevel optimization via {Krylov} subspace},
  author={Yang, Yan and Gao, Bin and Yuan, Ya-xiang},
  booktitle={The Thirteenth International Conference on Learning Representations (ICLR)},
  year={2025}
}

@article{ehrhardt2024efficient,
  title={Efficient gradient-based methods for bilevel learning via recycling {Krylov} subspaces},
  author={Ehrhardt, Matthias J and Gazzola, Silvia and Scott, Sebastian J},
  journal={arXiv preprint arXiv:2412.08264},
  year={2024}
}

@article{zhang2022adaptive,
  title={Adaptive stochastic variance reduction for subsampled {Newton} method with cubic regularization},
  author={Zhang, Junyu and Xiao, Lin and Zhang, Shuzhong},
  journal={INFORMS Journal on Optimization},
  volume={4},
  number={1},
  pages={45--64},
  year={2022},
  publisher={INFORMS}
}

@article{roosta2019sub,
  title={Sub-sampled {Newton} methods},
  author={Roosta-Khorasani, Farbod and Mahoney, Michael W},
  journal={Mathematical Programming},
  volume={174},
  number={1},
  pages={293--326},
  year={2019},
  publisher={Springer}
}

@article{nesterov2012efficiency,
  title={Efficiency of coordinate descent methods on huge-scale optimization problems},
  author={Nesterov, Yu},
  journal={SIAM Journal on Optimization},
  volume={22},
  number={2},
  pages={341--362},
  year={2012},
  publisher={SIAM}
}

@article{richtarik2014iteration,
  title={Iteration complexity of randomized block-coordinate descent methods for minimizing a composite function},
  author={Richt{\'a}rik, Peter and Tak{\'a}{\v{c}}, Martin},
  journal={Mathematical Programming},
  volume={144},
  number={1},
  pages={1--38},
  year={2014},
  publisher={Springer}
}

@inproceedings{flamary2015importance,
  title={Importance sampling strategy for non-convex randomized block-coordinate descent},
  author={Flamary, R{\'e}mi and Rakotomamonjy, Alain and Gasso, Gilles},
  booktitle={2015 IEEE 6th International Workshop on Computational Advances in Multi-Sensor Adaptive Processing (CAMSAP)},
  pages={301--304},
  year={2015},
  organization={IEEE}
}

@inproceedings{perekrestenko2017faster,
  title={Faster coordinate descent via adaptive importance sampling},
  author={Perekrestenko, Dmytro and Cevher, Volkan and Jaggi, Martin},
  booktitle={Artificial Intelligence and Statistics},
  pages={869--877},
  year={2017},
  organization={PMLR}
}

@article{rodomanov2020randomized,
  title={A randomized coordinate descent method with volume sampling},
  author={Rodomanov, Anton and Kropotov, Dmitry},
  journal={SIAM Journal on Optimization},
  volume={30},
  number={3},
  pages={1878--1904},
  year={2020},
  publisher={SIAM}
}

@article{stich2017safe,
  title={Safe adaptive importance sampling},
  author={Stich, Sebastian U and Raj, Anant and Jaggi, Martin},
  journal={Advances in Neural Information Processing Systems},
  volume={30},
  year={2017}
}

@article{sun2020markov,
  title={Markov chain block coordinate descent},
  author={Sun, Tao and Sun, Yuejiao and Xu, Yangyang and Yin, Wotao},
  journal={Computational Optimization and Applications},
  volume={75},
  number={1},
  pages={35--61},
  year={2020},
  publisher={Springer}
}

@article{peng2024block,
  title={Block acceleration without momentum: On optimal stepsizes of block gradient descent for least-squares},
  author={Peng, Liangzu and Yin, Wotao},
  journal={arXiv preprint arXiv:2405.16020},
  year={2024}
}

@article{peng2023block,
  title={Block coordinate descent on smooth manifolds: Convergence theory and twenty-one examples},
  author={Peng, Liangzu and Vidal, Ren{\'e}},
  journal={arXiv preprint arXiv:2305.14744},
  year={2023}
}

@article{platt1998sequential,
  title={Sequential minimal optimization: A fast algorithm for training support vector machines},
  author={Platt, John},
  year={1998}
}

@article{friedman2010regularization,
  title={Regularization paths for generalized linear models via coordinate descent},
  author={Friedman, Jerome H and Hastie, Trevor and Tibshirani, Rob},
  journal={Journal of Statistical Software},
  volume={33},
  pages={1--22},
  year={2010}
}

@article{shalev2013stochastic,
  title={Stochastic dual coordinate ascent methods for regularized loss minimization.},
  author={Shalev-Shwartz, Shai and Zhang, Tong},
  journal={Journal of Machine Learning Research},
  volume={14},
  number={1},
  year={2013}
}

@article{xu2013block,
  title={A block coordinate descent method for regularized multiconvex optimization with applications to nonnegative tensor factorization and completion},
  author={Xu, Yangyang and Yin, Wotao},
  journal={SIAM Journal on Imaging Sciences},
  volume={6},
  number={3},
  pages={1758--1789},
  year={2013},
  publisher={SIAM}
}

@article{birgin2022block,
  title={Block coordinate descent for smooth nonconvex constrained minimization},
  author={Birgin, Ernesto G and Mart{\'\i}nez, Jos{\'e} Mario},
  journal={Computational Optimization and Applications},
  volume={83},
  number={1},
  pages={1--27},
  year={2022},
  publisher={Springer}
}

@article{cheung2024randomized,
  title={Randomized submanifold subgradient method for optimization over Stiefel manifolds},
  author={Cheung, Andy Yat-Ming and Wang, Jinxin and Yue, Man-Chung and So, Anthony Man-Cho},
  journal={arXiv preprint arXiv:2409.01770},
  year={2024}
}

@inproceedings{cai2021zeroth,
  title={A zeroth-order block coordinate descent algorithm for huge-scale black-box optimization},
  author={Cai, HanQin and Lou, Yuchen and McKenzie, Daniel and Yin, Wotao},
  booktitle={International Conference on Machine Learning},
  pages={1193--1203},
  year={2021},
  organization={PMLR}
}

@article{park2025unraveling,
  title={Elucidating Subspace Perturbation in Zeroth-Order Optimization: Theory and Practice at Scale},
  author={Park, Sihwan and Yun, Jihun and Kim, SungYub and Kundu, Souvik and Yang, Eunho},
  journal={arXiv preprint arXiv:2501.19099},
  year={2025}
}

@article{jing2024asynchronous,
  title={Asynchronous distributed reinforcement learning for lqr control via zeroth-order block coordinate descent},
  author={Jing, Gangshan and Bai, He and George, Jemin and Chakrabortty, Aranya and Sharma, Piyush K},
  journal={IEEE Transactions on Automatic Control},
  volume={69},
  number={11},
  pages={7524--7539},
  year={2024},
  publisher={IEEE}
}

@article{david2024locally,
  title={Locally convergent bi-level MDO architectures based on the block coordinate descent algorithm},
  author={David, Yann and Gallard, Fran{\c{c}}ois and Rondepierre, Aude},
  journal={Journal of Optimization Theory and Applications},
  volume={211},
  number={2},
  pages={32},
  year={2026},
  publisher={Springer}
}

@article{gupta2023efficient,
  title={Efficient learning of decision-making models: A penalty block coordinate descent algorithm for data-driven inverse optimization},
  author={Gupta, Rishabh and Zhang, Qi},
  journal={Computers \& Chemical Engineering},
  volume={170},
  pages={108123},
  year={2023},
  publisher={Elsevier}
}

@article{lauga2025flexible,
  title={A flexible block-coordinate forward-backward algorithm for non-smooth and non-convex optimization},
  author={Brice{\~n}o-Arias, Luis and Gon{\c{c}}alves, Paulo and Lauga, Guillaume and Pustelnik, Nelly and Riccietti, Elisa},
  journal={arXiv preprint arXiv:2510.26477},
  year={2025}
}

@article{chambolle2011first,
  title={A first-order primal-dual algorithm for convex problems with applications to imaging},
  author={Chambolle, Antonin and Pock, Thomas},
  journal={Journal of mathematical imaging and vision},
  volume={40},
  number={1},
  pages={120--145},
  year={2011},
  publisher={Springer}
}

@article{heinkenschloss2014matrix,
  title={A matrix-free trust-region {SQP} method for equality constrained optimization},
  author={Heinkenschloss, Matthias and Ridzal, Denis},
  journal={SIAM Journal on Optimization},
  volume={24},
  number={3},
  pages={1507--1541},
  year={2014},
  publisher={SIAM}
}

\end{document}